\documentclass[11pt]{article}
\usepackage{amstext,amssymb,amsmath,amsbsy}

\usepackage{hyperref}
\usepackage{amscd}
\usepackage{amsfonts} 
\usepackage{indentfirst}
\usepackage{verbatim}
\usepackage{amsmath}
\usepackage{amsthm}
\usepackage{enumerate}
\usepackage{graphicx}
\usepackage{color}
\usepackage[OT1]{fontenc}
\usepackage[latin1]{inputenc}
\usepackage[english]{babel}
\usepackage{amssymb,MnSymbol}
\usepackage{tikz}
\usepackage{dsfont}

\newtheorem{theorem}{Theorem}
\newtheorem{lemma}{Lemma}

\newtheorem{proposition}{Proposition}
\newtheorem{corollary}{Corollary}
\newtheorem{remark}{Remark}

\newtheorem{question}{Open problem}
\newtheorem{definition}{Definition}
\newtheorem{pathology}{Pathology}

\newcommand{\proofend}{\hfill $\Box$ }
\newcommand{\mint}{\strokedint}

\newcommand{\bb}{\gamma}

\newcommand{\ess}{\mathrm{ess} \,}

\newcommand{\loc}{_{loc}}

\newcommand{\mN}{\mathbb{N}}

\newcommand{\mZ}{\mathbb{Z}}

\newcommand{\R}{{\cal R}}

\newcommand{\mR}{\mathbb{R}}
\newcommand{\mS}{\mathbb{S}}
\newcommand{\eps}{\varepsilon}

\newcommand{\dsp}{\displaystyle}
\newcommand{\mc}{\mathrm{c}}

\newcommand{\C}{\kappa}

\newcommand{\limf}[1]{\mathop{\liminf}_{#1}}
\newcommand{\lims}[1]{\mathop{\limsup}_{#1}}
\newcommand{\lr}{\rightarrow}

\newcommand{\bvarphi}{\tilde {\varphi} }
\newcommand{\hvarphi}{\hat \varphi }

\newcommand{\supp}{\operatorname{supp}}

\numberwithin{equation}{section}

\title{Non-local functionals related to the total variation and connections with Image Processing}
\author{Ha\"im Brezis\footnote{Rutgers University,
Department of Mathematics, Hill Center, Busch Campus,
110 Frelinghuysen Road, Piscataway, NJ 08854, USA, brezis@math.rutgers.edu} \footnote{Department of Mathematics,
Technion, Israel Institute of Technology,
32.000 Haifa, Israel} \footnote{Laboratoire Jacques-Louis Lions
UPMC,  4  place Jussieu, 75005 Paris,
France} \footnote{Research partially supported by NSF grant DMS-1207793 and by ITN "FIRST" of the European Commission, Grant  Number  PITN-GA-2009-238702.} and Hoai-Minh Nguyen\footnote{EPFL SB MATHAA CAMA, Station 8,  CH-1015 Lausanne,  Switzerland, hoai-minh.nguyen@epfl.ch}}

\begin{document}

\maketitle

\begin{abstract} We present new results concerning the approximation of the total variation, $\int_{\Omega} |\nabla u|$,  of a function $u$ by non-local, non-convex functionals of the form 
$$
\Lambda_\delta u = \int_{\Omega} \int_{\Omega} \frac{\delta \varphi \big( |u(x) - u(y)|/ \delta\big)}{|x - y|^{d+1}} \, dx \, dy, 
$$
as $\delta \to 0$, where $\Omega$ is a domain in $\mR^d$ and $\varphi: [0, + \infty) \to [0, + \infty)$ is a non-decreasing function satisfying some appropriate conditions. The mode of convergence is extremely delicate and numerous problems remain open. De Giorgi's concept of Gamma-convergence illuminates the situation, but also introduces mysterious novelties. The original motivation of our work comes from Image Processing. 

\end{abstract}

\noindent {\bf Key words}: Total variation, bounded variation, non-local functional, non-convex functional, Gamma-convergence, Sobolev spaces.

\medskip 
\noindent{\bf Mathematics Subject Classification}:  49Q20, 26B30, 46E35,  28A75. 
\tableofcontents

\section{Introduction}

Throughout this paper,  we assume that  $\varphi:[0, +\infty) \to [0, + \infty)$
is continuous on $[0, +\infty)$ except at a finite number of points in $(0, +\infty)$ where it admits a limit from the left and from the right. We also assume that $\varphi(0) = 0$ and that $\varphi(t) = \min \{ \varphi(t+), \varphi(t-)\}$ for all $t>0$, so that $\varphi$ is lower semi-continuous. We assume that the domain $\Omega \subset \mR^d$ is either bounded and smooth, or that $\Omega=\mR^d$. The case $d = 1$ is already of great interest;  many difficulties  (and open problems!) occur even when $d=1$.

\medskip
Given a measurable function $u$ on $\Omega$, and a small parameter $\delta > 0$, we define the following non-local functionals:
\begin{equation}\label{def-Lambda}
\Lambda (u): = \int_\Omega \int_\Omega \frac{\varphi(|u(x) - u(y)|) }{|x - y|^{ d  +  1}} \, dx \,
dy \quad \mbox{and} \quad \Lambda_{\delta}(u): = \delta \Lambda (u/\delta).
\end{equation}
Sometimes, it is convenient to be more specific and to write $\Lambda_\delta(u, \varphi, \Omega)$ or $\Lambda_\delta(u, \Omega)$ instead of $\Lambda_{\delta} (u)$.

\medskip
Our main goal in this paper is to study the asymptotic behaviour of $\Lambda_{\delta}$ as $\delta \to 0$.  In order to simplify the presentation we make, throughout the paper,  the following four basic assumptions on $\varphi$:
\begin{equation}\label{cond-varphi-0}
\varphi(t) \le a t^{2} \mbox{ in } [0,1] \mbox{ for some positive constant } a,
\end{equation}
\begin{equation}\label{cond-varphi-01}
\varphi(t) \le b  \mbox{ in } \mR_{+} \mbox{ for some positive constant } b,
\end{equation}
\begin{equation}\label{cond-varphi-decreasing}
\varphi \mbox{ is non-decreasing},
\end{equation}
and
\begin{equation}\label{cond-varphi3}
\gamma_{d}\int_0^\infty \varphi(t) t^{-2} \,dt=1, \, \mbox{where } \gamma_{d}: = \int_{\mS^{d-1}} |\sigma \cdot e | \, d \sigma \mbox{ for some } e \in \mS^{d-1}.
\end{equation}
A straightforward computation gives
\begin{equation}\label{def-gamma}
\gamma_d =  \left\{ \begin{array}{cl} \dsp \frac{2}{d-1} |\mS^{d-2}|  = 2 |B^{d-1}|& \mbox{ if } d \ge 3, \\[6pt]
4 & \mbox{ if } d = 2, \\[6pt]
2 & \mbox{ if } d =1, 
\end{array} \right.
\end{equation}
where $\mS^{d-2}$ (resp. $B^{d-1}$) denotes the unit sphere (resp. ball) in $\mR^{d-1}$. 

Condition~\eqref{cond-varphi3} is a normalization condition prescribed in order to  have \eqref{convergence-pointwise} below with constant 1 in front of $\int_{\Omega} |\nabla u|$. Denote
 \begin{equation}\label{def-calA}
{\cal A} = \big\{\varphi;\;\varphi \mbox{ satisfies } \eqref{cond-varphi-0}-\eqref{cond-varphi3} \big\}.
\end{equation}
Note that $\Lambda$ is {\bf never convex} when  $\varphi \in {\cal A}$.

We also mention the following additional condition on $\varphi$ which will  be imposed in Sections \ref{sec-pro} and \ref{sec-image}:
\begin{equation}\label{cond-varphi-positive}
 \varphi(t) > 0 \quad \mbox{ for all } t > 0.
\end{equation}

Note that if \eqref{cond-varphi-0} and  \eqref{cond-varphi-01} hold, then $\Lambda_{\delta}(u)$ is finite for every $u \in H^{1/2}(\Omega)$, and in particular for every $u \in C^{1}(\bar \Omega)$ when $\Omega$ is bounded.

\medskip
Here is  a list of specific examples  of functions $\varphi$ that we have in mind. They all satisfy \eqref{cond-varphi-0}-\eqref{cond-varphi-decreasing}. In order to achieve \eqref{cond-varphi3}, we choose $\varphi = c_i \bvarphi_i$ where $\bvarphi_i$ is taken from the list below and $c_i$ is an appropriate constant.

\medskip
\noindent{\bf Example 1}:
\begin{equation*}
\bvarphi_{1}(t) = \left\{\begin{array}{ll} 0 & \mbox{ if } t \le 1 \\[6pt]
1 & \mbox{ if } t > 1.
\end{array}\right.
\end{equation*}

\medskip
\noindent{\bf Example 2}:
\begin{equation*}
\bvarphi_{2}(t) = \left\{\begin{array}{ll} t^{2} & \mbox{ if } t \le 1 \\[6pt]
1 & \mbox{ if } t > 1.
\end{array}\right.
\end{equation*}

\medskip
\noindent{\bf Example 3}:
\begin{equation*}
\bvarphi_{3}(t) = 1 - e^{- t^2}.
\end{equation*}

%
%
%
%
%


The asymptotic behavior of $\Lambda_\delta$ as $\delta \to 0$ when $\varphi  = c_1 \bvarphi_1$ has been extensively studied in  \cite{BourNg, NgSob1, NgGammaCRAS, NgSob2,NgGamma, NgSob3, NgSob4, Po2}. Examples 2 and 3 are motivated by Image Processing (see Section~\ref{sec-image}).


In Section~\ref{sec-pointwise} we investigate the pointwise limit of $\Lambda_\delta$ as $\delta \to 0$, i.e., the convergence of $\Lambda_\delta(u)$ for fixed $u$.
We first consider the case where $u \in C^1(\bar \Omega)$, with $\Omega$ bounded, and prove (see Proposition~\ref{thm-Sobolev1p=1}) that
\begin{equation}\label{convergence-pointwise}
\Lambda_\delta(u) \mbox{ converges, as $\delta \to 0$, to } TV(u) = \int_{\Omega} |\nabla u|,  \mbox{ the total variation of $u$}.
\end{equation}
One may then be tempted to infer that the same conclusion holds for every $u \in W^{1, 1}(\Omega)$. Surprisingly, this is not true: for every $d \ge 1$ and for every $\varphi \in {\cal A}$, one can construct a function $u \in W^{1, 1}(\Omega)$ such that
\begin{equation*}
\lim_{\delta \to 0} \Lambda_\delta(u) =  + \infty.
\end{equation*}
Such a ``pathology"  was originally discovered by A. Ponce \cite{Po1} and presented in \cite{NgSob1} for $\varphi = c_1  \bvarphi_1$;
in Section~\ref{sect-pathology} we construct a simpler function $u \in W^{1,1}(\Omega)$ with the additional  property that $\Lambda_\delta (u, \varphi) =+ \infty$ for all $\delta > 0$ when $\varphi = c_2 \bvarphi_2$ (see Pathology~\ref{pathology-1}).

 If $u \in W^{1, 1}(\Omega)$, one may only assert (see Proposition~\ref{thm-Sobolev1p=1}) that
\begin{equation*}
\liminf_{\delta \to 0} \Lambda_{\delta} (u) \ge \int_{\Omega} |\nabla u|,
\end{equation*}
for every $\varphi \in {\cal A}$. 

When dealing with functions $u \in BV(\Omega)$, the situation becomes even more intricate as explained in Section~\ref{sect-pathology}. In particular, it may happen (see Pathology~\ref{pathology-3}) that, for some $\varphi \in {\cal A}$  and some $u \in BV(\Omega)$,
\begin{equation*}
\liminf_{\delta \to 0} \Lambda_{\delta} (u) < \int_{\Omega} |\nabla u|.
\end{equation*}
On the other hand, we prove (see \eqref{statement-part2}) that,  for every $\varphi \in {\cal A}$, 
\begin{equation*}
\liminf_{\delta \to 0} \Lambda_{\delta} (u) \ge K \int_{\Omega} |\nabla u| \quad \forall \, u \in L^1(\Omega),
\end{equation*}
for some $K \in (0, 1]$ depending only on $d$ and $\varphi$, and (see Proposition~\ref{thm-inversep=1-1-1})
\begin{equation*}
\limsup_{\delta \to 0} \Lambda_{\delta} (u) \ge  \int_{\Omega} |\nabla u| \quad \forall \, u \in L^1(\Omega).
\end{equation*}
Here and throughout the paper, we set $\int_{\Omega} |\nabla u| = + \infty$ if $u \in L^1(\Omega) \setminus BV(\Omega)$.

All these facts suggest that the mode of convergence of $\Lambda_\delta$ to $TV$ as $\delta \to 0$ is delicate and that pointwise convergence may be deceptive. It turns out that $\Gamma$-convergence (in the sense of E. De Giorgi) is the appropriate framework to analyze the asymptotic behavior of $\Lambda_\delta$ as $\delta \to 0$. (For the convenience of the reader, we recall the definition of $\Gamma$-convergence in Section~\ref{sec-gamma}).

\medskip
Section~\ref{sec-gamma} deals with  the following crucial result whose proof is extremely involved.

\begin{theorem}\label{thm-gamma}  Let $\varphi \in {\cal A}$. There exists a constant $K = K(\varphi) \in (0, 1]$, which is independent of $\Omega$ such that, as $\delta \to 0$,
\begin{equation}\label{gamma-limit}
\Lambda_{\delta}  \; \mbox{   $\Gamma$-converges to }
\Lambda_0  \mbox{ in } L^1(\Omega),
\end{equation}
where
\begin{equation*}
\Lambda_{0}(u): =  K \int_{\Omega} |\nabla u | \quad \mbox{ for } u \in L^1(\Omega).
\end{equation*}
\end{theorem}


Here is a direct consequence of Theorem~\ref{thm-gamma} in the case $d=1$. 
 
\begin{corollary}
Let $u \in L^1(0, 1)$ and $(u_\delta) \subset L^1(0, 1)$ be such that $u_\delta \to u$ in $L^1(0, 1)$. Then 
\begin{equation}\label{G-d=1}
\liminf_{\delta \to 0} \Lambda_{\delta}(u_\delta, (0, 1)) \ge K(\varphi)  |u(t_2)  - u(t_1)|, 
\end{equation}
for all Lebesgue points $t_1, t_2 \in (0, 1)$ of $u$.  
\end{corollary}

Despite its simplicity, we do not know an easy proof for this assertion even when $u_\delta \equiv u$, $\varphi = c_1 \bvarphi_1$,  and  $K(\varphi)$ is replaced by a positive constant independent of $u$; in this case the proof is due to J. Bourgain and H-M. Nguyen  in \cite[Lemma 2]{BourNg}. 

\begin{remark} \fontfamily{m} \selectfont
The constant $K$ may also depend on $d$ (actually we have not investigated whether it really depends on $d$), but for simplicity we omit this (possible) dependence. 

\end{remark}

\begin{remark} \fontfamily{m} \selectfont
When $\varphi = c_1 \bvarphi_{1}$ and $\Omega = \mR^d$ (with $c_1$ is chosen such that \eqref{cond-varphi3} holds), Theorem~\ref{thm-gamma} is originally due to H-M.~Nguyen \cite{NgGammaCRAS}, \cite{NgGamma}. The lengthy proof of Theorem~\ref{thm-gamma} borrows many ideas from \cite{NgGamma}. The study of non-local functionals of the type $\Lambda_{\delta}$ was initiated in  \cite{NgSob1} at the suggestion of H. Brezis (see \cite{BrSurvey}).
\end{remark}

\begin{remark} \fontfamily{m} \selectfont Let $(\delta_n)$ be any positive sequence converging to 0 as $n \to + \infty$. Then $\Lambda_{\delta_n}$  $\Gamma$-converges to $\Lambda_0$ in $L^1(\Omega)$. This is a consequence of Theorem~\ref{thm-gamma} and the general fact: if a family of functionals $(F_\delta)$ $\Gamma$-converges to a functional  $F_0$ in $L^1(\Omega)$  as $\delta \to 0$ then $(F_{\delta_n})$ also $\Gamma$-converges to $F_0$ in $L^1(\Omega)$ as $n \to + \infty$.

\end{remark}

It would be interesting to remove the monotonicity assumption \eqref{cond-varphi-decreasing} in the definition of ${\cal A}$. More precisely, we have

\begin{question}
Assume  that \eqref{cond-varphi-0}, \eqref{cond-varphi-01}, and \eqref{cond-varphi3} hold.
Is it true that either the conclusion of Theorem~\ref{thm-gamma} holds,  or $(\Lambda_\delta)$ $\Gamma$-converges in $L^1(\Omega)$ to 0 as $\delta \to 0$?
\end{question}

The answer is positive in the one dimensional case \cite{BNg-Gamma-non}.

\begin{remark} \fontfamily{m} \selectfont
Note that if one removes the monotonicity assumption on $\varphi$ it may happen that $\dsp (\Lambda_\delta) \mathop{\to}^{\Gamma} 0 $ in $L^1(\Omega)$ as $\delta \to 0$.
This occurs e.g.,  when  $\supp \varphi \subset \subset (0, + \infty)$.  Indeed, given $u \in L^1(\Omega)$, let $(\tilde u_\delta)$ be a family of functions converging in $L^1(\Omega)$ to $ u $, as $\delta$ goes to 0, such that $\tilde u_\delta$ takes its values in the set $m \delta \mZ$. Here $m$ is chosen such that   $|t| < m/2$ for $t \in \supp \varphi$. It is clear that
\begin{equation*}
\Lambda_\delta(\tilde u_\delta) = 0, \quad \forall \, \delta > 0.
\end{equation*}
Therefore $\dsp \Lambda_\delta \; \mathop{\to}^\Gamma \; 0$ in $L^{1}(\Omega)$.
\end{remark}


The appearance of the constant $K = K(\varphi)$ in Theorem~\ref{thm-gamma} is mysterious and somewhat counterintuitive. Assume for example that $\Omega$ is bounded and  that $u \in C^1(\bar \Omega)$. We know that $\Lambda_\delta(u) \to \int_{\Omega} |\nabla u|$ as $\delta \to 0$ (see Proposition~\ref{thm-Sobolev1p=1}). On the other hand, it follows from Theorem~\ref{thm-gamma} that there exists a family $(u_\delta)$ in $L^1(\Omega)$ such that $u_\delta \to u$ in $L^1(\Omega)$ and $\Lambda_\delta(u_\delta) \to K \int_{\Omega}|\nabla u|$ as $\delta \to 0$.
The reader may wonder how $K$ is determined.  This is rather easy to explain, e.g., when $d=1$ and $\Omega = (0, 1)$; $K(\varphi)$ is given by
\begin{equation}\label{def-ku-11}
K(\varphi) = \inf \liminf_{\delta \to 0} \Lambda_{\delta} (v_\delta),
\end{equation}
where the infimum is taken over all families of  functions $(v_\delta)_{\delta \in (0,1)} \subset L^1(0, 1)$ such that $v_\delta \to v_0$ in $L^1(0, 1)$ as $\delta \to 0$ with $v_0(x) = x$ in $(0, 1)$. Unfortunately, formula \eqref{def-ku-11} provides very  little information about the constant  $K(\varphi)$ (for example, the explicit value of $K(\varphi)$ is not known even when $\varphi = c_1 \bvarphi_1$).  Taking $v_\delta = v_0$ for all $\delta > 0$,
we obtain $K(\varphi) \le 1$ for all $\varphi$. Indeed, an easy computation using the normalization \eqref{cond-varphi3}  shows  (see Proposition~\ref{thm-Sobolev1p=1}) that
$\lim_{\delta \to 0} \Lambda_{\delta}(v_0) = \int_{0}^1 |v_0'| = 1$. A more sophisticated choice of $(v_\delta)$ in \cite{NgGammaCRAS} yields $K(\varphi)< 1$ when $\varphi = c_1 \bvarphi_1$, for every $d \ge 1$. For the convenience of the reader, we include the proof of this fact in Section~\ref{sect3.6}.  
On the other hand, it is nontrivial that $K(\varphi)> 0$ for all $\varphi \in {\cal A}$ and $d \ge 1$. This is not obvious  even in the simple case where $\varphi = c_1 \bvarphi_1$  (see Section~\ref{sect-property-k}); 
the proof  relies heavily on ideas from \cite{BourNg, NgSob2}. It is even less trivial that $\inf_{\varphi \in {\cal A}} K(\varphi) > 0$ (see Section~\ref{sect-k}).

\bigskip

Here is a challenging question, which is open even when $d=1$.

\begin{question}\label{qG3}
Is it always true that $K(\varphi) < 1$ in Theorem~\ref{thm-gamma}? Or even better: Is it true that $\sup_{\varphi \in {\cal A}} K(\varphi) < 1$?
\end{question}

We believe that indeed $K(\varphi)< 1$ for every $\varphi$. (However, if it turns out that $K(\varphi)=1$ for some $\varphi$'s,
it would be interesting to characterize  such  $\varphi$'s.)

\medskip

In Section~\ref{sec-pro}, we establish the following two compactness results. The first one deals with the level sets of $\Lambda_{\delta}$ for a {\bf fixed} $\delta$, e.g.,  for  $\delta = 1$.

\begin{theorem}\label{thm-compact-p=1} Let $\varphi \in {\cal A}$ satisfy  \eqref{cond-varphi-positive},
and let $(u_n)$ be a bounded sequence in $ L^1(\Omega)$ such that
\begin{equation}\label{compactness-p=1}
\sup_{n} \Lambda(u_n) < + \infty.
\end{equation}
There exists a subsequence
$(u_{n_k})$  of $(u_n)$ and $u \in L^1(\Omega)$ such that $(u_{n_k})$ converges to $u$ in $L^{1}(\Omega)$ if $\Omega$ is bounded, resp. in $L^{1}_{\loc}(\mR^d)$ if $\Omega = \mR^d$.
\end{theorem}

The second result concerns a sequence $\big(\Lambda_{\delta_n}\big)$ with $\delta_n \to 0$;  here \eqref{cond-varphi-positive} is not required.

\begin{theorem}\label{thm-compact2} Let $\varphi \in {\cal A}$, $(\delta_n) \to 0$, and let  $(u_n)$ be a bounded sequence in $ L^1(\Omega)$ such that
\begin{equation}\label{bd}
\sup_{n} \Lambda_{\delta_n}(u_n) < + \infty.
\end{equation}
There exists a subsequence
$(u_{n_k})$  of $(u_n)$ and $u \in L^1(\Omega)$ such that $(u_{n_k})$ converges to $u$ in $L^{1}(\Omega)$ if $\Omega$ is bounded, resp. in $L^{1}_{\loc}(\mR^d)$ if $\Omega = \mR^d$.
\end{theorem}

Theorem~\ref{thm-compact2} is due to H-M. Nguyen \cite{NgSob3} when $\Omega = \mR^d$ and $\varphi = c_1 \bvarphi_1$.

\medskip
In Section~\ref{sec-image}, we consider problems of the form
\begin{equation}\label{minimizing-problem}
\inf_{u \in L^{q}(\Omega)} E_{\delta}(u),
\end{equation}
in the case $\Omega$ bounded,  where
\begin{equation}\label{defE}
E_{\delta}(u) = \lambda \int_{\Omega} |u - f|^{q}  + \Lambda_{\delta}(u),
\end{equation}
$q\geq1$,  $f \in L^{q}(\Omega)$ is given, and  $\lambda $ is a fixed positive constant. Our goal is twofold: investigate the existence of minimizers for $E_{\delta}$ ($\delta$ being fixed) and analyze their behavior as $\delta \to 0$.
The existence of a minimizer in \eqref{minimizing-problem} is not obvious since $\Lambda_{\delta}$ is  {\bf not convex} and one cannot invoke the standard tools of Functional Analysis.  Theorem~\ref{thm-compact-p=1} implies the existence of a minimizer in \eqref{minimizing-problem}. Next we study the behavior of these minimizers as $\delta \to 0$. More precisely, we prove

\begin{theorem}\label{thm-image} Assume that $\Omega$ is bounded, and that $\varphi \in {\cal A}$  satisfies  \eqref{cond-varphi-positive}.  Let $q\geq1$,  $f \in L^{q}(\Omega)$, and let $u_\delta$ be a minimizer of \eqref{defE}. Then $u_{\delta} \to u_{0}$ in $L^{q}(\Omega)$ as $\delta \to 0$,  where $u_{0}$ is the unique minimizer of the functional $E_{0}$ defined on $L^{q}(\Omega)\cap BV(\Omega)$ by
\begin{equation*}
E_{0}(u) : = \lambda \int_{\Omega} |u - f|^{q} + K \int_{\Omega} |\nabla u|,
\end{equation*}
and $0 < K \le 1$ is the constant coming from Theorem~\ref{thm-gamma}.
\end{theorem}


Basic ingredients in the proof are the $\Gamma$-convergence result (Theorem~\ref{thm-gamma}) and the compactness result (Theorem~\ref{thm-compact2}).

\medskip
As explained in Section~\ref{sec-image}, $E_{\delta}$ and $E_{0}$ are closely related to  functionals used in Image Processing for the purpose of denoising the image $f$. In fact, $E_{0}$ corresponds to the celebrated ROF filter  originally introduced by L. I. Rudin, S. Osher and E. Fatemi in \cite{RudinOsherFatemi}.  While $E_{\delta}$ (with $\varphi$ as in Examples 2-3) is reminiscent of filters introduced by L.~S.~Lee~\cite{Lee} and  L. P.~ Yaroslavsky  (see \cite{Yaroslavsky1, YaroslavskyEden}).  More details can be found in  the expository paper by A. ~Buades, B. Coll, and J.~M.~Morel
\cite{BuadesCollMorelReview1}; see also \cite{BuadesCollMorel2,  Morel1, SmithBrady, PKTD} where various terms, such as  ``neighbourhood filters",``non-local means" and ``bilateral filters",  are used.  Some of these filters admit a variational formulation,  as explained  by S.~Kindermann, S.~Osher and P.~W. Jones in \cite{KindermannOsherJones}. Theorem~\ref{thm-image} says that such  filters ``converge" to the ROF filter,  as $\delta \to 0$,   a fact which seems to be new to the experts in Image Processing.

\medskip
In a forthcoming paper \cite{BrezisNguyen-p>1} we investigate similar questions where the functional $\Lambda_{\delta}(u)$ is replaced by
\begin{equation}\label{p}
\Lambda_{\delta} (p, u) =  \delta^{p} \int_{\Omega} \int_{\Omega} \frac{\varphi(|u(x) - u(y)|/ \delta)}{|x - y|^{d+ p}} \, dx \, dy,
\end{equation}
and the total variation is replaced by $\int_{\Omega} |\nabla u|^{p}$ (with $p>1$).

In recent years there has been much interest in the convergence of {\bf convex} non-local functionals to the total variation, going  back to the work of  J.~Bourgain, H.~Brezis and P.~Mironescu  \cite{BBM} (see Remark~\ref{rem-BBM} below).
Related works may be found  in  \cite{B1, Davila, BBM2, Po, VanWi, CRS, LeoniSpector, Fusco, Figalli, Bucur, BN-Two, BN-BBM, BrMirbook}. For the convergence of non-local functionals to the perimeter, we mention in particular \cite{CV, ADM},   the two surveys \cite[Section 5]{Bucur}, \cite[Section 5.6]{Fusco},  and  the references therein.
 As one can see, there is a ``family resemblance'' with questions studied in our paper.
We warn the reader that the non-convexity of $\Lambda_\delta$ is a source of major difficulties. Moreover, new and surprising phenomena emerged  over the past fifteen years, in particular the discovery in \cite{NgGammaCRAS, NgGamma} that the $\Gamma$-limit and the pointwise limit of $(\Lambda_\delta)$ do not coincide; we refer to  \cite{BrSurvey} for some historical comments.
We also mention that a different type of approximation of the BV-norm of a function $u$,
 especially suited when $u$ is the characteristic function  of a set $A$,   so that its $BV$-norm is the perimeter of $A$,  has been recently developed in \cite{ABBF1} and \cite{ABBF2} (with roots in \cite{BBM2-1}).

\medskip 
Part of the results in this paper are announced in \cite{BrSurvey, NgSob4, BN-CRAS-long}.

\section{Pointwise convergence of  $\Lambda_\delta$ as $\delta \to 0$}\label{sec-pointwise}

\subsection{Some positive results}

The first result in this section is

\begin{proposition} \label{thm-Sobolev1p=1}
Assume that $\varphi \in {\cal A}$. Then
\begin{equation}\label{concl-char1-p=1}
\lim_{\delta \to 0} \Lambda_{\delta}(u)  =  \int_{\Omega} |\nabla u |,
\end{equation}
for all $u \in C^1(\overline{\Omega})$ if $\Omega$ is bounded, resp. for all $u \in C^1_{\mc}(\mR^d)$ if $\Omega = \mR^d$. 
However, if $u \in W^{1,1}(\Omega)$ (with  $\Omega$ bounded  or $\Omega = \mR^d$), we can only assert that
\begin{equation}\label{ine-Fatou}
\liminf_{\delta \to 0} \Lambda_\delta (u) \ge \int_{\Omega} |\nabla u|,
\end{equation}
and strict inequality may happen (see Pathology~\ref{pathology-1} in Section~\ref{sect-pathology}).
\end{proposition}

\begin{remark}\label{rem-BBM}  \fontfamily{m} \selectfont  The convergence of a special sequence of {\bf convex} non-local functionals to the total variation was originally analyzed by J.~Bourgain, H.~Brezis and P.~Mironescu  \cite{BBM} and further  investigated in  \cite{B1, Davila, BBM2, Po, VanWi, LeoniSpector, BrMirbook, BN-Two, BN-BBM}. More precisely, it has been shown that, for every $u \in L^{1}(\Omega)$,
\begin{equation}\label{limitBBM}
\lim_{\eps \to 0} J_{\eps} (u) = \gamma_d \int_{\Omega} |\nabla u|,
\end{equation}
where
\begin{equation}\label{def-I}
J_{\eps}(u) = \int_{\Omega} \int_{\Omega} \frac{|u(x) - u(y)|}{|x - y|} \rho_{\eps} (|x -y|) \, dx \, dy.
\end{equation}
Here $\gamma_d$ is defined in \eqref{cond-varphi3},  $\rho_{\eps}$ is an arbitrary sequence of radial mollifiers (normalized by the condition $\int_0^\infty \rho_\eps(r) r^{d-1} \, dr = 1$). As the reader can see, \eqref{concl-char1-p=1} and \eqref{limitBBM} look somewhat similar. However, the asymptotic analysis of $\Lambda_\delta$ is much more delicate because two basic properties satisfied by $J_\eps$ are {\bf not} fulfilled by $\Lambda_\delta$:
\begin{enumerate}
\item[i)] there is {\bf no} constant $C$ such that, e.g.,  with $\Omega$  bounded,
\begin{equation}\label{pro1-JL}
\Lambda_\delta(u) \le  C \int_{\Omega} |\nabla u| \quad \forall \, u \in C^1(\bar \Omega), \, \forall \, \delta > 0,
\end{equation}
despite the fact $\lim_{\delta \to 0} \Lambda_\delta(u) = \int_{\Omega} |\nabla u|$ for all $u \in C^1(\bar \Omega)$. Indeed, if \eqref{pro1-JL} held, we would deduce by density the same estimate for every $u \in W^{1, 1}(\Omega)$ and this fact contradicts Pathology~\ref{pathology-1}  in Section~\ref{sect-pathology}. 


\item[ii)] $\Lambda_\delta(u)$ is {\bf not} a convex functional.

\end{enumerate}
\end{remark}

It is known (see \cite{Po}) that the $\Gamma$-limit and the pointwise limit of $(J_\eps)$ coincide and are equal to $\gamma_d \int_{\Omega} |\nabla \cdot |$. By contrast, this is not true  for $\Lambda_{\delta}$ since the constant $K$ in Theorem~\ref{thm-gamma} might be less than 1.

\medskip 
Here and in what follows in this paper, given a function  $\varphi: [0, + \infty) \to \mR$ and $\delta > 0$, we denote $\varphi_\delta$ the function 
$$
\varphi_\delta (t) =\delta \varphi (t / \delta) \mbox{ for } t \ge 0. 
$$
With this notation, one has   
$$
\Lambda_{\delta}(u, \varphi) = \Lambda(u, \varphi_\delta).
$$ 

\medskip

\noindent {\bf Proof of  Proposition~\ref{thm-Sobolev1p=1}.} We first consider the case $\Omega  = \mR^d$ and $u \in C^1_{\mc}(\mR^d)$.  
Fix  $M > 1$  such that $u(x) = 0$ if $|x| \ge M-1$.   We have
\begin{equation*}
\Lambda_\delta(u) =  \int_{|x| > M} \,
dx \int_{\mR^d} \frac{\varphi_\delta (|u(x) - u(y)|) }{|x - y|^{ d  +  1}} \, dy  +  \int_{|x| \le M}  \,
dx \int_{\mR^d} \frac{\varphi_\delta (|u(x) - u(y)|) }{|x - y|^{ d  +  1}} \, dy.
\end{equation*}
Since $\varphi$ is bounded and
\begin{equation*}
 \int_{|x| > M} \, dx  \int_{|y| < M-1} \frac{1}{|x - y|^{ d  +  1}} \, dy  < + \infty,
\end{equation*}
it follows from the choice of $M$ that
\begin{equation}\label{cl1}
 \lim_{\delta \to 0} \int_{|x| > M} \, dx  \int_{\mR^d} \frac{\varphi_\delta (|u(x) - u(y)|) }{|x - y|^{ d  +  1}} \, dy = 0.
\end{equation}
Replacing $y$ by $x + z$ and using polar coordinates in the $z$ variable, we find
\begin{equation}\label{tt1}
\int_{|x| \le  M} \, dx  \int_{\mR^d} \frac{\varphi_\delta(|u(x) - u(y)|)}{|x -y|^{d + 1}} \, dy = \int_{|x| \le M} \, dx 
\int_0^{+\infty}   \, dh  \int_{\mS^{d-1}} \frac{\varphi_\delta(|u(x + h \sigma) - u(x)|)}{h^{2}} \, d \sigma .
\end{equation}
We have
\begin{multline}\label{tt2}
\int_{|x| \le M} \, dx \int_0^{+\infty}\, d h \int_{\mS^{d-1}} \frac{\varphi_\delta(|u(x + h \sigma) - u(x)|)}{h^{2}} \, d
\sigma \\[6pt]
= \int_{|x| \le M} \, dx  \int_0^{+\infty} \, d h \int_{\mS^{d-1}} \frac{\delta \varphi \Big(|u(x + h \sigma) - u(x)| \big/
\delta \Big)}{h^{2}} \, d \sigma.
\end{multline}
Rescaling the variable $h$ gives
\begin{multline}\label{tt3}
\int_{|x| \le M} \, dx \int_0^{+\infty} \, dh  \int_{\mS^{d-1}}  \frac{\delta \varphi \Big(|u(x + h \sigma) - u(x)| \big/
\delta \Big)}{h^{2}}  \, d \sigma\\[6pt]
= \int_{|x| \le M} \, dx  \int_0^{+\infty} \, d h \int_{\mS^{d-1}}\frac{ \varphi \Big(|u(x + \delta h \sigma) - u(x)| \big/
\delta \Big)}{h^{2}} \, d \sigma  .
\end{multline}
Combining \eqref{tt1}, \eqref{tt2}, and \eqref{tt3} yields
\begin{equation}\label{tt4}
\int_{|x| \le M} \ dx  \int_{\mR^d} \frac{\varphi_\delta(|u(x) - u(y)|)}{|x -y|^{d+1}} \, dy = \int_{|x| \le M} \, dx
\int_0^{+\infty} \, d h \int_{\mS^{d-1}} \frac{ \varphi \Big(|u(x + \delta h \sigma) - u(x)| \big/
\delta \Big)}{h^{2}}  \, d \sigma .
\end{equation}
Note that
\begin{equation}\label{tt5}
\lim_{\delta \to 0}\frac{|u(x + \delta h \sigma) - u(x)|}{\delta}  = |\nabla u(x) \cdot \sigma| h  \mbox{ for
} (x, \, h, \, \sigma) \in \mR^d \times [0, + \infty) \times \mS^{d-1}.
\end{equation}
Since $\varphi$ is continuous at $0$ and on $(0, +\infty)$ except at a finite number of points, it follows that
\begin{multline}\label{mono1}
\lim_{\delta \to 0}\frac{1}{h^{2}} \varphi \Big(|u(x + \delta h \sigma) - u(x)|\big/\delta \Big) =
\frac{1}{h^{2}}\varphi \Big(|\nabla u(x) \cdot \sigma| h \Big) \\[6pt]
\mbox{ for  a.e. } (x, \, h, \, \sigma) \in \mR^d \times (0, + \infty) \times \mS^{d-1}
\end{multline}
(if $|\nabla u (x) \cdot \sigma| h$ is a point of discontinuity of $\varphi$, we may change a little bit $h$).
Rescaling once more the variable $h$ gives
\begin{equation}\label{ttt}
 \int_{0}^{\infty}  \, d h \int_{\mS^{d-1}} \frac{1}{h^{2}} \varphi\Big( |\nabla u (x) \cdot \sigma| h\Big)  \, d \sigma = |\nabla u(x)| \int_{0}^{\infty} \varphi(t) t^{-2}  \, dt \; \int_{\mS^{d-1}} |\sigma \cdot e| \, d \sigma ;
\end{equation}
here we have also used the obvious fact that, for every $V \in \mR^d$, and for any fixed $e \in \mS^{d-1}$,
\begin{equation}\label{V}
\int_{\mS^{d-1}} |V \cdot \sigma| \, d \sigma = |V| \int_{\mS^{d-1}} |\sigma \cdot e| \, d \sigma. 
\end{equation}
Thus, by the normalization condition \eqref{cond-varphi3}, we obtain
\begin{equation}\label{limit}
\int_{|x| \le M} \, dx  \int_{0}^\infty  \, d h \int_{\mS^{d-1}} \frac{1}{h^{2}} \varphi\Big( |\nabla u (x) \cdot \sigma| h\Big)  \, d \sigma
 =  \int_{|x| \le M} |\nabla u| \, dx.
\end{equation}
Define $ \hvarphi: [0, \infty) \to \mR$ as follows
\begin{equation*}
\hvarphi (t) = \left\{\begin{array}{cl} (a +b) t^{2} & \mbox{ if }  0 \le t \le 1, \\[6pt]
a + b & \mbox{ if } t > 1,
\end{array}\right.
\end{equation*}
where $a$ and $b$ are the constants  in \eqref{cond-varphi-0} and \eqref{cond-varphi-01}.
Then
\begin{equation}\label{choice1}
\varphi \le \hvarphi \mbox{ on } [0, + \infty).
\end{equation}
We note that
\begin{equation}\label{int1}
\int_{0}^{\infty} \hvarphi(t) t^{-2} \, dt < + \infty.
\end{equation}
Since $u \in C^{1}_{\mc}(\mR^{d})$, it is clear that
\begin{equation}\label{int2}
\frac{|u(x + \delta h \sigma) - u(x)|}{\delta}  \le C h  \mbox{ for
} (x, \, h, \, \sigma) \in \mR^d \times [0, + \infty) \times \mS^{d-1},
\end{equation}
for some positive constant $C$.
On the other hand,  by \eqref{int1},
\begin{equation}\label{int3}
\int_{|x| \le M} \, dx  \int_{0}^{\infty} \, d h  \int_{\mS^{d-1}} \frac{1}{h^{2}}\hvarphi(C h) \, d \sigma   < + \infty.
\end{equation}
Applying the dominated convergence theorem, and using \eqref{tt4}, \eqref{mono1}, \eqref{limit}, \eqref{choice1},  \eqref{int2} and \eqref{int3}, we find
\begin{equation}\label{cl2}
\lim_{\delta \to 0} \int_{|x| \le M} \, dx  \int_{\mR^d} \frac{\varphi_\delta(|u(x) - u(y)|)}{|x -y|^{d + 1}} \, dy   = \int_{|x| \le M} |\nabla u| \, dx.
\end{equation}
Assertion \eqref{concl-char1-p=1} now follows from \eqref{cl1} and \eqref{cl2}.

\medskip
The proof of \eqref{ine-Fatou} is almost identical, even simpler. In fact \eqref{ine-Fatou} is an immediate consequence of \eqref{mono1} and \eqref{ttt},  and Fatou's lemma.

\medskip

We next consider the case where $\Omega$ is bounded.   Let $D \subset \subset \Omega$ and fix $t > 0$ small enough such that $B(x, t)= \{ y \in \mR^d; \; |y -x| < t \} \subset \subset \Omega $ for every $x \in D$.   We have, for every $u \in W^{1, 1}(\Omega)$,
\begin{equation*}
\Lambda_{\delta}(u) \ge  \int_{D} \, dx \int_{B(x, t)} \frac{\varphi_\delta(|u(x) - u(y)|)}{|x -y|^{d+1}} \, dy = \int_{D} \, dx
\int_0^{t/\delta} \int_{\mS^{d-1}} \frac{\varphi \Big(|u(x + \delta h \sigma) - u(x)| \big/ \delta \Big)}{h^{2}} \, d \sigma  \, dh.
\end{equation*}
By the same method as above, we deduce that
\begin{equation}\label{p2-S2}
\liminf_{\delta \to 0}\Lambda_\delta(u)  \ge  \int_{D} |\nabla u| \quad \forall \, u \in W^{1, 1}(\Omega);
\end{equation}
which implies \eqref{ine-Fatou} since $D \subset \Omega$ is arbitrary.

In order to prove \eqref{concl-char1-p=1} for every $u \in C^1(\bar \Omega)$, we write
\begin{equation*}
\Lambda_{\delta} (u)= A_\delta + B_\delta + C_\delta,
\end{equation*}
where
\begin{equation*}
A_\delta = \int_D \,dx \int_{B(x, t)} \frac{\varphi_\delta(|u(x) - u(y)| )}{|x - y|^{d+1}} \, dy,
\end{equation*}
\begin{equation*}
B_\delta = \int_D \,dx \int_{\Omega \setminus B(x, t)} \frac{\varphi_\delta(|u(x) - u(y)| )}{|x - y|^{d+1}} \, dy,
\end{equation*}
and
\begin{equation*}
C_\delta = \int_{\Omega \setminus D} \,dx \int_{\Omega} \frac{\varphi_\delta(|u(x) - u(y)| )}{|x - y|^{d+1}} \, dy.
\end{equation*}
By the same method as above, we find
\begin{equation}\label{p11-P1}
\lim_{\delta \to 0} A_\delta = \int_{D} |\nabla u|.
\end{equation}
On the other hand, we have
\begin{equation}\label{p12-P1}
B_\delta \le \delta b |\Omega|^2/ t^{d+1},
\end{equation}
and, as above,
\begin{equation*}
C_\delta \le \delta \int_{\Omega \setminus D}\, dx \int_{\Omega} \frac{\hvarphi(L |x - y|/ \delta)}{|x - y|^{d+1}} \, dy,
\end{equation*}
where $L$ is the Lipschitz constant of $u$ on $\Omega$.  An immediate computation gives
\begin{equation}\label{p13-P1}
C_\delta \le C |\Omega \setminus D|,
\end{equation}
where $C$ depends only on $L$, $a$, $b$, and $d$. It is clear that
\begin{equation*}
\big|\Lambda_{\delta}(u) - \int_{\Omega} |\nabla u| \big| \le \big|A_\delta - \int_{D} |\nabla u| \big| + B_\delta + C_\delta + \int_{\Omega \setminus D} |\nabla u|.
\end{equation*}
Using \eqref{p11-P1}, \eqref{p12-P1}, and \eqref{p13-P1}, we conclude that
\begin{equation*}
\limsup_{\delta \to 0} |\Lambda_{\delta } (u) - \int_{\Omega} |\nabla u|| \le C |\Omega \setminus D|;
\end{equation*}
which implies \eqref{concl-char1-p=1} since $D$ is arbitrary. The proof is complete. \proofend

\begin{remark} \fontfamily{m} \selectfont  We call the attention of the reader that the monotonicity assumption~\eqref{cond-varphi-decreasing} on $\varphi$ has {\bf not} been used in the proof of Proposition~\ref{thm-Sobolev1p=1}.
\end{remark}

\begin{remark}\label{RemarkW1p}  \fontfamily{m} \selectfont  The condition $u \in C^1(\overline{\Omega})$ if $\Omega$ is bounded (resp. $u \in C^{1}_{\mc}(\mR^{d})$  if $\Omega = \mR^d$) in \eqref{concl-char1-p=1} is much too strong.  In fact, the same conclusion  holds under the assumption that $\Omega$ is bounded and $u$ is Lipschitz (with an identical proof).   More generally, equality~\eqref{concl-char1-p=1} holds  e.g.,   when $u \in W^{1, p}(\Omega)$ for some $p>1$, and $\Omega$ is bounded (see Proposition~\ref{pro-C} in Appendix~\ref{C}).  It would be interesting to characterize the set
$$
\Big\{u \in W^{1, 1}(\Omega); \; \lim_{\delta \to 0} \Lambda_{\delta}(u) = \int_{\Omega} |\nabla u| \Big\}.
$$
\end{remark}

So far we have been dealing with the pointwise convergence of $\Lambda_\delta(u)$ when $u \in W^{1, 1}(\Omega)$, but it is natural to ask similar questions when $u \in BV(\Omega)$. As a consequence of Theorem~\ref{thm-gamma}, we know that  for every $\varphi \in {\cal A}$, there exists  a constant $K = K(\varphi) \in(0, 1]$ such that
\begin{equation}\label{statement-part2}
 \liminf_{\delta \to 0} \Lambda_{\delta}(u) \ge K  \int_{\Omega} |\nabla u|  \quad \forall \, u \in L^1(\Omega).
\end{equation}
On the other hand, we also have

\begin{proposition}\label{thm-inversep=1-1-1} Assume that $\varphi \in {\cal A}$. Then
\begin{equation}
 \limsup_{\delta \to 0} \Lambda_{\delta}(u) \ge  \int_{\Omega} |\nabla u|  \quad \forall \, u \in L^1(\Omega).
\end{equation}
\end{proposition}

\noindent {\bf Proof of Proposition~\ref{thm-inversep=1-1-1}.}  The proof relies on  ideas  from  \cite{NgSob1}.  It suffices to consider the case
\begin{equation}\label{limsup-P1}
F: = \limsup_{\delta \to 0} \Lambda_\delta(u) < + \infty.
\end{equation}
We first assume that $u \in L^\infty(\Omega)$. Set
\begin{equation}\label{def-A}
A = 2 \| u\|_{L^\infty}.
\end{equation}
Fix $0< \delta_0 < 1$.  Set, for  $0< \eps < 1/2$,
\begin{equation}\label{Tepsdelta}
T(\eps, \delta_0): =  \int_0^{\delta_{0}}   \eps \delta^{\eps -1 } \Lambda_\delta(u) \, d \delta =  \int_0^{\delta_{0}}   \eps \delta^{\eps -1 } \, d\delta \int_{\Omega} \int_{\Omega} \frac{\delta \varphi(|u(x) - u(y)| / \delta)}{|x-y|^{d + 1}}  \, dx \, dy .
\end{equation}
Using Fubini's theorem and integrating first with respect to $\delta$,  we have
\begin{equation*}
T(\eps, \delta_0) = \int_{\Omega} \int_{\Omega} \frac{\eps |u(x) - u(y)|^{1 + \eps}}{|x - y|^{d+1}} \, dx \, dy \int_{|u(x) - u(y)|/ \delta_0}^\infty \varphi(t) t^{-2 - \eps}  \, dt.
\end{equation*}
This implies
\begin{equation*}
T(\eps, \delta_0) \ge c(\eps, \delta_0) \mathop{\int_{\Omega} \int_{\Omega}}_{|u(x) - u(y)| < \delta_0^2} \frac{\eps |u(x) - u(y)|^{1 + \eps}}{|x - y|^{d+1}}\, dx \, dy,
\end{equation*}
where
\begin{equation*}
c(\eps, \delta_0) = \int_{\delta_0}^\infty \varphi(t) t^{-2 - \eps} \, dt.
\end{equation*}
It follows from \eqref{def-A} that
\begin{equation}\label{p0-Pro3}
T(\eps, \delta_0) \ge c(\eps, \delta_0) \int_{\Omega} \int_{\Omega} \frac{\eps |u(x) - u(y)|^{1 + \eps}}{|x - y|^{d+1}}\, dx \, dy -  c(\eps, \delta_0) \mathop{\int_{\Omega} \int_{\Omega}}_{|u(x) - u(y)| \ge  \delta_0^2} \frac{\eps A^{1+  \eps}}{|x - y|^{d+1}}\, dx \, dy.
\end{equation}
Let $\tau > 0$ be arbitrary small. First choose $\delta_0$ small enough such that
\begin{equation}\label{p0-Pro3-1}
 \int_{\delta_0}^\infty \varphi(t) t^{-2} \, dt  \ge \gamma_d^{-1} ( 1 - \tau)
\end{equation}
and
\begin{equation}\label{limsup-Pro2}
\Lambda_\delta(u) \le  F + \tau  \quad \forall \, 0< \delta < \delta_0.
\end{equation}
We next observe that
\begin{equation}\label{alpha-finite-Pro2}
\mathop{\int_{\Omega} \int_{\Omega}}_{|u(x) - u(y)| \ge \alpha} \frac{1}{|x - y|^{d+1}} \, d x \, dy < + \infty \quad \forall \, \alpha > 0.
\end{equation}
Indeed, fix $t_0 > 0$  such that $\varphi(t_0) > 0$ and note
\begin{equation}
\Lambda_\delta(u)  \ge   \mathop{\int_{\Omega} \int_{\Omega}}_{|u(x) - u(y)| \ge \alpha} \frac{\varphi_\delta (|u(x) - u(y)|)}{|x - y|^{d+1}} \, dx \, dy \ge \delta \varphi(\alpha/ \delta) \mathop{\int_{\Omega} \int_{\Omega}}_{|u(x) - u(y)| \ge \alpha} \frac{1}{|x - y|^{d+1}} \,dx \, dy.
\end{equation}
Choosing $0< \delta < \min\{ \delta_0, \alpha/ t_0\}$ and using \eqref{limsup-Pro2}, we obtain \eqref{alpha-finite-Pro2}.
We deduce from \eqref{alpha-finite-Pro2} that
\begin{equation}\label{p1-Pro3}
\lim_{\eps \to 0} c(\eps, \delta_0) \mathop{\int_{\Omega} \int_{\Omega}}_{|u(x) - u(y)| >  \delta_0^2} \frac{\eps A^{1+  \eps}}{|x - y|^{d+1}}\, dx \, dy = 0.
\end{equation}
We next invoke the following lemma which is an immediate consequence of the BBM formula \eqref{limitBBM} applied with  $\rho_\eps(t) = \eps t^{\eps - d} \mathds{1}_{(0, 1)}$ (see \cite[Proposition 1]{BN-Two}).

\begin{lemma}\label{lem-limit} We have
\begin{equation}\label{p2-Pro3}
\liminf_{\eps \to 0} \gamma_d^{-1}\int_{\Omega} \int_{\Omega} \frac{\eps |u(x) - u(y)|^{1 + \eps}}{|x - y|^{d+1}}\, dx \, dy \ge \int_{\Omega} |\nabla u|  \quad \forall \, u \in L^1(\Omega).
\end{equation}
\end{lemma}

Combining  \eqref{p0-Pro3}, \eqref{p0-Pro3-1}, \eqref{p1-Pro3}, and \eqref{p2-Pro3} yields
\begin{equation}\label{p3-Pro3}
\liminf_{\eps \to 0} T(\eps, \delta_0) \ge (1 - \tau) \int_{\Omega} |\nabla u|.
\end{equation}
On the other hand, using \eqref{Tepsdelta} and \eqref{limsup-Pro2}, we find
\begin{equation*}
T(\eps, \delta_0) \le \int_0^{\delta_0} \eps \delta^{\eps-1} (F  + \tau) \, d \delta = (F + \tau ) \delta_0^\eps,
\end{equation*}
so that
\begin{equation}\label{p4-Pro3}
\limsup_{\eps \to 0} T(\eps, \delta_0) \le F + \tau.
\end{equation}
From \eqref{p3-Pro3} and \eqref{p4-Pro3}, we deduce that
\begin{equation*}
F + \tau \ge (1 - \tau) \int_{\Omega} |\nabla u|.
\end{equation*}
Since $\tau > 0$ is arbitrary, we obtain
\begin{equation*}
\limsup_{\delta \to 0} \Lambda_\delta (u) \ge \int_{\Omega} |\nabla u|.
\end{equation*}
The proof is complete in the case $u \in L^\infty(\Omega)$. In the general case, we proceed  as follows. Set, for $A > 0$,
\begin{equation}\label{def-TA}
T_A(s) = \left\{\begin{array}{cl} s & \mbox{ if } |s| \le A, \\[6pt]
A & \mbox{ if } s  > A, \\[6pt]
-A & \mbox{ if } s <- A,
\end{array}\right.
\end{equation}
and
\begin{equation*}
u_A = T_A(u).
\end{equation*}
Since $\varphi$ is non decreasing,
\begin{equation*}
\Lambda_\delta(u_A) \le \Lambda_\delta(u).
\end{equation*}
It follows that
\begin{equation*}
\int_{\Omega} |\nabla u_A| \le \limsup_{\delta \to 0} \Lambda_\delta(u_A) \le \limsup_{\delta \to 0} \Lambda_\delta(u).
\end{equation*}
By letting $A \to + \infty$, we obtain
\begin{equation*}
\int_{\Omega} |\nabla u| \le \limsup_{\delta \to 0} \Lambda_\delta(u).
\end{equation*}
The proof is complete. \proofend

\subsection{Some pathologies} \label{sect-pathology}

Our first example is related to Proposition~\ref{thm-Sobolev1p=1} and shows that inequality \eqref{ine-Fatou} can be strict.



\begin{pathology}\label{pathology-1} Let $d \ge 1$. There exists $u \in W^{1, 1}(\Omega)$ such that
$$
\lim_{\delta \to 0} \Lambda_{\delta} (u) = + \infty \mbox{ for all } \varphi \in {\cal A};
$$
moreover,
$$
 \Lambda_{\delta} (u) = + \infty \quad  \forall \,  \delta > 0  \mbox{ for } \varphi = c_2 \bvarphi_2.
$$
\end{pathology}

\noindent{\bf Proof.} For simplicity, we present only the case $d=1$ and choose $\Omega = (-1/2, 1/2)$.  Define, for $\alpha > 0$,
\begin{equation*}
u(x) = \left\{\begin{array}{cl} 0 & \mbox{ if } -1/2 < x < 0, \\[6pt]
|\ln x|^{-\alpha} & \mbox{ if } 0 < x < 1/2.
\end{array}\right.
\end{equation*}
Clearly, $u \in W^{1, 1}(\Omega)$.  We claim that, for $0 < \alpha < 1$,
$$
\lim_{\delta \to 0} \Lambda_{\delta} (u) = + \infty \mbox{ for } \varphi = c_1 \bvarphi_1
$$
and, for $0  < \alpha < 1/2$,
$$
 \Lambda_{\delta} (u) = + \infty \quad  \forall \,  \delta > 0  \mbox{ for } \varphi = c_2 \bvarphi_2.
$$
It is clear that the conclusion follows from the claim  since for all $\varphi \in {\cal A}$ there exist $\alpha, \beta > 0$ such that $\varphi (t) \ge \alpha c_1 \bvarphi_1(\beta t)$ for all  $t > 0$.

\medskip 
It remains to prove the claim. For $\varphi = c_1 \bvarphi_1$, we have
$$
\Lambda_{\delta}(u) \ge c_1 \mathop{\int_0^{1/2}}_{|u(x)| > \delta} \, dx  \int_{-1/2}^0  \frac{\delta}{|x - y|^2} \, dy.
$$
For $\delta$ sufficiently small, let $x_\delta \in (0, 1/2)$ be the unique solution of  $|\ln x |^{-\alpha} = \delta$. A straightforward computation yields
$$
\Lambda_{\delta} (u) \ge  c_1 \delta \int_{x_\delta}^{1/2}  \Big(  \frac{1}{x} - \frac{1}{x+ 1/2}\Big) \, dx \sim  \delta |\ln x_\delta | =  \delta^{1 - 1/\alpha} \to + \infty \mbox{ as } \delta \to 0,
$$
if $\alpha < 1$. We now consider the case $\varphi = c_2 \bvarphi_2$. We have, since $|u| \le 1$,
$$
\Lambda_{\delta} (u) \ge C_\delta   \int_0^{1/2} \, dx \int_{-1/2}^0   \frac{|u(x)|^2}{|x - y|^2} \, dy =  C_\delta \int_0^{1/2} |\ln x |^{-2 \alpha} \Big(  \frac{1}{x} - \frac{1}{x+ 1/2}\Big)  \, d x = + \infty,
$$
if $2 \alpha < 1$.  \proofend

\medskip
Next, we mention an example of $\varphi \in {\cal A}$ and $u \in W^{1,1}$ such that $\lim_{\delta \to 0} \Lambda_{\delta}(u)$ does {\bf not} exist and the gap between $\liminf_{\delta \to 0} \Lambda_{\delta} (u)$ and $\limsup_{\delta \to 0} \Lambda_{\delta} (u)$
is ``maximal".

\begin{pathology}\label{pathology-2} Let $\Omega = (0, 1)$. There exists a function $\varphi \in {\cal A}$ and a function $u \in W^{1, 1}(\Omega)$ such that
\begin{equation}\label{ine-Fatou1}
\liminf_{\delta \to 0} \Lambda_{\delta} (u) = \int_{\Omega} |\nabla u| \quad \mbox{ and } \quad  \limsup_{\delta \to 0} \Lambda_{\delta} (u) = + \infty.
\end{equation}
\end{pathology}

The construction is presented in Appendix~\ref{A}. Our next example shows that assertion~\eqref{ine-Fatou} in Proposition~\ref{thm-Sobolev1p=1} may fail for $u \in BV(\Omega) \setminus W^{1,1}(\Omega)$.

\begin{pathology}\label{pathology-3} Let $\Omega = (0, 1)$. There exists a continuous function $\varphi \in {\cal A}$  and a function $u \in BV(\Omega) \cap C(\bar \Omega)$  such that
\begin{equation}\label{ine-Fatou2}
\liminf_{\delta \to 0} \Lambda_{\delta} (u) < \int_{\Omega} |\nabla u|.
\end{equation}
\end{pathology}

The construction is presented in Appendix~\ref{B}.

\medskip
\noindent{\bf Concluding remark}: the abundance of pathologies is quite mystifying and  a reasonable theory of pointwise convergence of $\Lambda_{\delta}$ seems out of reach.
Fortunately, $\Gamma$-convergence saves the situation!

\section{$\Gamma$-convergence of $\Lambda_\delta$ as $\delta \to 0$.}\label{sec-gamma}

\subsection{Structure  of the proof of Theorem~\ref{thm-gamma}}

Recall that (see  e.g., \cite{Braides, Dal}), by definition, a family of functionals $(\Lambda_{\delta})_{\delta \in (0, 1)}$ defined on $L^1(\Omega)$ (with values in $\mR \cup \{+\infty\}$), $\Gamma$-converges to $\Lambda_0$ in $L^1(\Omega)$ as $\delta \to 0$ if the following two properties hold:
\begin{enumerate}
\item[(G1)] For every $u \in L^1(\Omega)$ and for every family
$(u_\delta)_{\delta \in (0,1)} \subset L^1(\Omega)$ such that
$u_\delta \to u$ in $L^1(\Omega)$ as $\delta \to 0$,
one has
\begin{equation*}
\liminf_{\delta \to 0} \Lambda_\delta(u_\delta) \ge \Lambda_0(u).
\end{equation*}
\item[(G2)] For every $u \in L^1(\Omega)$, there exists a family
$(\tilde u_\delta)_{\delta \in (0,1)} \subset L^1(\Omega)$ such that
$\tilde u_\delta \to u$ in $L^1(\Omega)$ as $\delta \to 0$,
and
\begin{equation*}
\limsup_{\delta \to 0}  \Lambda_\delta(\tilde u_{\delta}) \le  \Lambda_0(u).
\end{equation*}
\end{enumerate}

The constant $K$ which occurs in Theorem~\ref{thm-gamma} will be defined via a ``semi-explicit" construction. More precisely, {\bf fix} any (smooth)  function $u \in {\cal B}: =  \Big\{u \in BV(\Omega); \; \int_{\Omega} |\nabla u| = 1  \Big\}$; given any $\varphi \in {\cal A} =  \{\varphi;\; \varphi \mbox{ satisfies } \eqref{cond-varphi-0}-\eqref{cond-varphi3} \}$, set
\begin{equation}\label{def-ku}
K(u, \varphi, \Omega) = \inf \liminf_{\delta \to 0} \Lambda_{\delta} (v_\delta),
\end{equation}
where the infimum is taken over all families of  functions $(v_\delta)_{\delta \in (0,1)} \subset L^1(\Omega)$ such that $v_\delta \to u$ in $L^1(\Omega)$ as $\delta \to 0$.

\medskip
We will eventually establish that
\begin{equation}
\mbox{$K(u, \varphi, \Omega)$ is independent of $u$ and $\Omega$; it depends only on $\varphi$ and $d$,}
\end{equation}
and
\begin{equation}
\mbox{Theorem~\ref{thm-gamma} holds with  $K = K(u, \varphi, \Omega)$.}
\end{equation}

 A priori, it is very surprising that $K(u, \varphi, \Omega)$ is independent of $u \in {\cal B}$.
 However, a posteriori,  if one believes Theorem~\ref{thm-gamma}, this becomes natural.
Indeed,
\begin{equation*}
\liminf_{\delta \to 0} \Lambda_{\delta} (u_\delta) \ge K \int_{\Omega} |\nabla u| = K,
\end{equation*}
for every family $\dsp (u_\delta) \; \mathop{\to}^{L^1} \; u \in {\cal B}$ by $(G1)$, and thus $K(u, \varphi, \Omega) \ge K$. On the other hand, by $(G2)$, there exists a family $\dsp (\tilde u_\delta) \; \mathop{\to}^{L^1} \; u \in {\cal B}$ such that
\begin{equation*}
\limsup_{\delta \to 0} \Lambda_{\delta} (\tilde u_\delta) \le K \int_{\Omega} |\nabla u| = K,
\end{equation*}
and hence $K(u, \varphi, \Omega) \le K$.
      
\medskip 
In view of what we just said, the special choice of $u$ and $\Omega$ is irrelevant. For convenience, we define, for $\varphi \in {\cal A}$,
\begin{equation}\label{k-Gamma}
\C(\varphi) = K(U, \varphi, Q),
\end{equation} 
where 
$$
Q = [0, 1]^d \quad \mbox{ and } \quad U(x) : = (x_1 + \dots + x_d) \big/ \sqrt{d} \mbox{ in } Q, 
$$ 
so that $\int_{Q} |\nabla U| = 1$. 


\medskip
Here is a comment about property (G2). From Property (G2), it follows easily that a stronger form of (G2) holds: 
\begin{enumerate}
\item[(G2')] For every $u \in L^1(\Omega)$, there exists a family
$(\hat  u_\delta)_{\delta \in (0,1)} \subset L^1(\Omega) \cap  L^\infty(\Omega)$ such that
$\hat u_\delta \to u$ in $L^1(\Omega)$ as $\delta \to 0$,
and
\begin{equation*}
\limsup_{\delta \to 0}  \Lambda_\delta(\hat u_{\delta}) \le  \Lambda_0(u).
\end{equation*}
\end{enumerate}
Indeed, it suffices to take 
$$
\hat u_\delta = T_{A_\delta}(\tilde u_\delta), 
$$ 
where $T_A$ denotes the truncation at the level $A$ (see \eqref{def-TA})  and $A_\delta \to \infty$.  This leads naturally to the following 

\begin{question}\label{OP3} Given $u \in L^1(\Omega)$, is it possible to find 
$(\hat  u_\delta)_{\delta \in (0,1)} \subset L^1(\Omega) \cap C^0(\bar \Omega)$ (resp. $W^{1,1}(\Omega)$, resp. $L^1(\Omega) \cap C^\infty(\bar \Omega)$) such that
$\hat u_\delta \to u$ in $L^1(\Omega)$ as $\delta \to 0$,
and
\begin{equation*}
\limsup_{\delta \to 0}  \Lambda_\delta(\hat u_{\delta}) \le  \Lambda_0(u)? 
\end{equation*}
The question is open even if $\Omega = (0, 1)$, $u(x) = x$, and $\varphi =c_1 \bvarphi$.
\end{question}
The heart of the matter is the non-convexity of $\varphi$, so that one {\bf cannot}  use  convolution.  If the answer to Open problem \ref{OP3} is negative,  this would be a kind of  Lavrentiev gap phenomenon. In that case, it would  be very interesting to study the asymptotics as $\delta \to 0$ of ${\Lambda_\delta}_{| L^1(\Omega) \cap C^0(\bar \Omega)}$ (with numerous possible variants). 

\medskip

In Section~\ref{sect-property-k}, we prove  that
\begin{equation}\label{property-k}
0< \C(\varphi) \le 1 \mbox{ for all } \varphi \in {\cal A}.
\end{equation}


In Section~\ref{sect-proof-G2}, we prove  Property (G2) in Theorem~\ref{thm-gamma}.


In Section~\ref{sect-proof-G1}, we prove Property  (G1) in Theorem~\ref{thm-gamma}.




In Section~\ref{sect-k}, we discuss further properties  of    $\C(\varphi)$. In particular, we show that $\inf_{\varphi \in {\cal A}} \C(\varphi)  > 0$.

In Section~\ref{sect3.6}, we prove that $\C(c_1 \bvarphi_1) < 1$. 

\subsection{Proof of \eqref{property-k}} \label{sect-property-k}

By \eqref{concl-char1-p=1} in Proposition~\ref{thm-Sobolev1p=1}, we have
\begin{equation*}
\lim_{\delta \to 0}\Lambda_\delta(U, \varphi, Q) = \int_{Q} |\nabla U| = 1
\end{equation*}
(the reader may be concerned that $Q$ is not smooth, but the conclusion of Proposition~\ref{thm-Sobolev1p=1} can be easily extended to this case). Hence $\C (\varphi) \le 1$ by the definition \eqref{def-ku}  applied with $U$ and $Q$.

We next claim that   $\C(\varphi)>0$. Recall that,  by \cite[Theorem 2, formulas (1.2) and (1.3)]{NgGamma}  
$$
\liminf_{\delta \to 0} \Lambda_{\delta} (v_\delta, c_1 \bvarphi_1, Q) \ge K_1 \int_{Q} |\nabla U| = K_1,
$$
for every sequence $v_\delta \to U$ in $L^1(Q)$ and for some positive constant $K_1$ (here we also use the fact that convergence in $L^1(Q)$ implies convergence in measure in  $Q$).
On the other hand, it is easy to check that for every $\varphi \in {\cal A}$ there exist  $\alpha, \beta > 0$ such that
\begin{equation*}
\varphi (t) \ge \alpha c_1 \bvarphi_1(\beta t) \quad \forall \, t > 0.
\end{equation*}
Thus, for every sequence $(v_\delta) \to U$ in $L^1(Q)$,
\begin{equation*}
\liminf_{\delta \to 0} \Lambda_{\delta} (v_\delta, \varphi, Q) \ge \alpha \beta K_1 > 0.
\end{equation*}
Consequently,
$$
\C(\varphi) > 0.
$$
\proofend

\subsection{Proof of Property (G2)} 
\label{sect-proof-G2}

\subsubsection{Preliminaries} \label{sect-preliminaries}

This section is devoted to several lemmas  which  are used in the proof of Property (G2) (some of them are also used in the proof of Property (G1)) and are
in the spirit of \cite[Sections 2 and 3]{NgGamma}.  
In this section, $\varphi \in {\cal A}$ is fixed (arbitrary) and $0 < \delta < 1$. We recall that 
\begin{equation} \label{varphi-delta}
\varphi_\delta (t) = \delta \varphi(t/ \delta) \quad \mbox{ for } t \ge 0. 
\end{equation}
All subsets $A$ of $\mR^d$ are assumed to be measurable  and $C$ denotes a positive constant depending only on $d$ unless stated otherwise. For $A \subset \mR^d$ and $f: A \mapsto \mR$, we denote  $\mathrm{Lip}(f, A)$  the Lipschitz constant of $f$ on $A$.

\medskip 
We begin with

\begin{lemma} \label{lemtechnical1} Let  $A \subset \mR^d$ and  $f, g$ be measurable functions on $A$.  Define $h_1 = \min (f, g)$ and $h_2 = \max (f, g)$. We have
\begin{equation}\label{estLemfond1}
\Lambda_\delta (h_1, A) \le \Lambda_\delta(f, A)
+ \iint_{A^2 \setminus B_1^2} \frac{ \varphi_\delta (|g(x) - g(y)|)}{|x - y |^{d + 1}} \ dx \, dy 
 \end{equation}
and
\begin{equation}\label{estLemfond1*}
\Lambda_\delta (h_2, A) \le \Lambda_\delta(f, A)
+ \iint_{A^2 \setminus B_2^2} \frac{ \varphi_\delta (|g(x) - g(y)|)}{|x - y |^{d + 1}} \ dx \, dy, 
\end{equation}
where
$$
B_1 = \big\{ x \in A; f(x) \le g(x) \big\} \quad \mbox{and} \quad B_2 = \big\{ x \in A;\, f(x) \ge  g(x) \big\}.
$$
Assume in addition that  $g$ is Lipschitz on $A$  and $L = \mathrm{Lip}(g, A)$. Then
\begin{equation}\label{estLemfond2}
\Lambda_\delta (h_1, A) \le \Lambda_\delta(f, A) + C L |A \setminus B_1|
\end{equation}
and
\begin{equation}\label{estLemfond3}
\Lambda_\delta (h_2, A) \le \Lambda_\delta(f, A)+ C L |A \setminus B_2|. 
\end{equation}
\end{lemma}

\noindent{\bf Proof.} It suffices to prove \eqref{estLemfond1} and \eqref{estLemfond2} since \eqref{estLemfond1*} and \eqref{estLemfond3} are consequences of \eqref{estLemfond1} and \eqref{estLemfond2} by considering $-f$ and $-g$. We first prove \eqref{estLemfond1}. One can easily verify that   
$$
|h_1(x) - h_1(y)| \le \max (|f(x) - f(y)|, |g(x) - g(y)|).
$$
This implies \eqref{estLemfond1} since $\varphi \ge 0$ and $\varphi$ is non-decreasing.

To obtain \eqref{estLemfond2} from \eqref{estLemfond1}, one just notes that, since 
 $|g(x) - g(y)| \le L |x -y|$ and $\varphi$ is non-decreasing, 
 \begin{equation*}
\iint_{A^2 \setminus B_1^2} \frac{\varphi_\delta (|g(x) - g(y)|)}{|x - y |^{d + 1}} \ dx \, dy \le
\iint_{A^2 \setminus B_1^2} \frac{\varphi_\delta(L|x-y|)}{|x - y |^{d + 1}} \ dx \, dy, 
\end{equation*}
and, by a change of variables and the normalization condition of $\varphi$,  
$$
\iint_{A^2 \setminus B_1^2} \frac{\varphi_\delta(L|x-y|)}{|x - y |^{d + 1}} \ dx \, dy \le  2 \int_{A \setminus B_1}  \, dx \int_{\mS^{d-1}}  \, d \sigma \int_0^\infty\frac{\varphi_\delta (L r)}{r^2} \, d r  \le C L  |A \setminus B_1|. 
$$
\proofend


\medskip
Here is an immediate consequence of Lemma~\ref{lemtechnical1}.

\begin{corollary}\label{corCor00}
Let   $-\infty \le m_1 < m_2 \le + \infty$, $A \subset \mR^d$, and $f$ be  a measurable function on $A$. Set 
$$
h= \min\big(\max(f,
m_1), m_2\big). 
$$
We have
\begin{equation}\label{Cor00}
\Lambda_\delta(h, A) \le \Lambda_\delta(f, A). 
\end{equation}
\end{corollary}

Another useful consequence of Lemma~\ref{lemtechnical1} is
\begin{corollary}\label{cortechnical2}
Let $c>0$, $A \subset \mR^d$,  and $f, g$ be measurable functions on $A$.  Set $ B= \big\{x \in A; \, |f(x)-g(x)| >c \big\}$, 
$$
h= \min \big(\max(f, g-c), g+c \big). 
$$
Assume  that $g$ is
Lipschitz on $A$ with  $L = \mathrm{Lip}(g, A)$. We have 
\begin{equation*}
\Lambda_{\delta}(h, A) \le \Lambda_{\delta}(f, A) +   C L|B|. 
\end{equation*}
\end{corollary}

 An important consequence of Corollary~\ref{cortechnical2} is
\begin{corollary}\label{cortechnical3}
Let $A \subset \mR^d$, $g \in L^\infty(A)$, $(\delta_k) \subset \mR$, and $(g_k) \subset L^1(A)$ be  such
that $A$ is bounded,   $g$ is Lipschitz,  and  $g_k \to g$ in  $L^1(A)$. There exists $(h_k) \subset L^\infty(A)$  such that
$\| h_k - g\|_{L^\infty(A)} \to 0$ and
\begin{equation*}
\limsup_{k \to \infty} \Lambda_{\delta_k}(h_k, A) \le
\limsup_{k \to \infty} \Lambda_{\delta_k}(g_k, A).
\end{equation*}
Similarly, if $g_\delta \to g$ in $L^1(A)$ as $\delta \to 0$, there exists $(h_\delta) \subset L^\infty(A)$ such that $\| h_\delta - g \|_{L^\infty(A)} \to 0$ and 
\begin{equation*}
\limsup_{\delta  \to 0} \Lambda_{\delta}(h_\delta, A) \le
\limsup_{\delta \to 0} \Lambda_{\delta}(g_\delta, A).
\end{equation*}
\end{corollary}

\noindent{\bf Proof}: Set $c_k = \| g_k - g \|_{L^1(A)}^{1/2}$. Then $c_k |A_k| \le \| g_k - g\|_{L^1(A)} = c_k^2$ where $A_k = \{x \in A; |g_k(x) - g(x)| > c_k \}$, so that 
\begin{equation}\label{pro-pro-k}
\lim_{k \to + \infty} c_k = 0 \quad \mbox{ and } \quad  \lim_{k \to + \infty} |A_k| = 0. 
\end{equation}
Define $h_k = \min( \max(g_k, g - c_k), g + c_k)$ in $A$. Clearly $\|h_k - g \|_{L^\infty(A)} \le c_k$. Applying
Corollary~\ref{cortechnical2}, we have
\begin{equation}\label{pro-pro-k-1}
 \Lambda_{\delta_k}(h_k, A)  \le
 \Lambda_{\delta_k}(g_k, A) + C L |A_k|,  
\end{equation}
where $L$ is the Lipschitz constant of $g$. Letting $k \to + \infty$ in \eqref{pro-pro-k-1} and using \eqref{pro-pro-k}, one reaches  the conclusion for $(g_k)$. The argument for $(g_\delta)$ is the same.   \proofend

\medskip 
We now introduce some notations used  later. We denote 
\begin{enumerate}
\item[i)] for $x, y \in \mR^d$, 
$$|x-y|_\infty =  \sup_{i=1, \cdots, d} |x_i - y_i|.$$ 

\item[ii)] for $c>0$ and $A \subset \mR^d$,  
\begin{equation}\label{defAc}
A_{c}  = \big\{x \in A ; \; \mbox{dist}_{\infty}(x, \partial A) \le c \big\},
\end{equation}
where 
\begin{equation*}
\mbox{dist}_{\infty}(x, \partial A)  : = \inf_{y \in \partial A} |x-y|_\infty. 
\end{equation*}

\item[iii)]  for $c \in \mR$ and for $A, B \subset \Omega$, 
$$
cA = \{ ca \; ;  a \in A \}
$$
and
$$
A+ B : = \{a + b \; ; \, a \in A \mbox{ and } b \in B \}.
$$
We write $A + v$ instead of $A+ \{ v\}$ for $v \in \mR^d$.  
\end{enumerate}

We now present an estimate which will be used repeatedly in this paper. 

\begin{lemma} \label{lem-computation} Let $c>0$, $g \in L^1(\mR^d)$, and let $D$  be a  Lipschitz, bounded open subset of $\mR^d$. Assume that  $g$ is Lipschitz in $D_c$ with $L = \mathrm{Lip}(g,D_c)$. We have 
\begin{equation*}
\int_{D} \int_{\mR^d} \frac{\varphi_{\delta}(|g(x) - g(y)|)}{|x - y|^{d+1}} \, dx \, dy \le  \Lambda_{\delta} (g, D \setminus D_{c/2}) + C_D \Big(L c + b\delta/ c \Big) \mbox{ for } \delta > 0, 
\end{equation*}
for some positive constant $C_D$ depending only on $D$ where $b$ is the constant in \eqref{cond-varphi-01}. 
\end{lemma}

\noindent{\bf Proof.} Set 
$$
A_1 = (D \setminus D_{3c / 4}) \times (D \setminus D_{c/2}), \quad A_2 = D_{3 c/ 4} \times  D_c, 
$$
and
$$
A_3 = \Big((D \setminus D_{3c/4}) \times \big(\mR^d \setminus (D \setminus D_{c/2}) \big) \Big) \cup \Big( D_{3 c/ 4}  \times (\mR^d \setminus D_c) \Big). 
$$
It is clear that $D \times \mR^d \subset A_1 \cup A_2 \cup A_3$ and $A_1 \subset (D \setminus D_{c/2}) \times (D \setminus D_{c/2})$. A straightforward computation yields 
$$
\iint_{A_2} \frac{\varphi_{\delta}(|g(x) - g(y)|)}{|x - y|^{d+1}} \, dx \, dy \le \iint_{A_2} \frac{\varphi_{\delta}(L|x-y|)}{|x - y|^{d+1}} \, dx \, dy \le C_D L c. 
$$
and, since $\varphi \le b$ and  if $(x, y) \in A_3$ then $x \in D$ and  $|x - y| \ge C_D c$,  
$$
\iint_{A_3} \frac{\varphi_{\delta}(|g(x) - g(y)|)}{|x - y|^{d+1}} \, dx \, dy \le \iint_{A_3} \frac{\delta b}{|x - y|^{d+1}} \, dx \, dy \le \int_{D} \, dx \int_{\mS^{d-1}} \, d \sigma  \int_{C_D c}^\infty \frac{\delta b}{h^2} \, dh  \le C_D \delta b / c.  
$$
Therefore, 
\begin{equation*}
\int_{D} \int_{\mR^d} \frac{\varphi_{\delta}(|g(x) - g(y)|)}{|x - y|^{d+1}} \, dx \, dy \le  \Lambda_{\delta} (g, D \setminus D_{c/2}) + C_D \Big( Lc + b\delta/ c \Big). 
\end{equation*}
\proofend 

\medskip 

We have

\begin{lemma}\label{lemcorlemfondN-111} Let $(\delta_k)  \subset \mR_+$  and  $(g_k) \subset L^1(Q)$ be such that $\delta_k \to 0$ and   $g_k \to U$ in $L^1(Q)$. There exist  $(c_k) \subset \mR_+$ and $(h_k) \subset L^\infty(Q)$ such that 
$$
c_k \ge  \sqrt{\delta_k},  \quad  \lim_{k \to + \infty} c_k = 0, 
$$
$$
\| h_k - U \|_{L^\infty(Q)} \le 2 d c_k, \quad  \mathrm{Lip}(h_k, Q_{c_k}) \le 1,
$$
and 
\begin{equation*}
\lims{k \to + \infty} \Lambda_{\delta_k} (h_k, Q) \le \lims{k \to + \infty} \Lambda_{\delta_k} (g_k, Q).  
\end{equation*}
Similarly, if $(g_\delta) \subset L^1(Q)$ is such that $g_\delta \to U$ in $L^1(Q)$, there exist $(c_\delta) \subset \mR_+$ and $(h_\delta) \subset L^\infty(Q)$ such that 
$$
c_\delta \ge  \sqrt{\delta_\delta},  \quad  \lim_{\delta \to 0} c_\delta = 0, 
$$
$$
\| h_\delta - U \|_{L^\infty(Q)} \le 2 d c_\delta, \quad  \mathrm{Lip}(h_\delta, Q_{c_\delta}) \le 1,
$$
and 
\begin{equation*}
\lims{\delta \to 0} \Lambda_{\delta} (h_\delta, Q) \le \lims{\delta \to 0} \Lambda_{\delta} (g_\delta, Q).  
\end{equation*}
\end{lemma}

\noindent {\bf Proof.} We only give the proof for the sequence $(g_k)$. The proof for the family $(g_\delta)$ is the same. By Corollary~\ref{cortechnical3}, 
one may assume that $\|g_k - U\|_{L^\infty(Q)} \to 0$. Set
\begin{equation}\label{defck}
c_k = \max \Big(\|g_k - U\|_{L^\infty(Q)}, \sqrt{\delta_k} \Big), 
\end{equation}
denote $g_{0, k} = g_k$, and define
\begin{equation}\label{defgk} \begin{array}{l}
g_{1, k}(x) =  \min \Big(\max \Big(g_{0,k}(x), U(0, x_2, \dots, x_d)
+ 2 c_k \Big),
U(1, x_2, \dots, x_d) - 2 c_k \Big), \\[6pt]
g_{2, k} (x) =  \min \Big(\max \Big(g_{1,k}(x), U(x_1, 0, \dots, x_d)
+ 4 c_k \Big),
U(x_1, 1, \dots, x_d) - 4 c_k \Big),\\[6pt]
\dots \\[6pt]
g_{d, k}(x) =  \min \Big(\max \Big(g_{d-1,k}(x), U(x_1, \dots,
x_{d-1},0) + 2d c_k \Big),  U(x_1,  \dots, x_{d-1}, 1) - 2 d c_k \Big). 
\end{array}
\end{equation}
From the definition of $U$, we have
\begin{equation*}
\left\{\begin{array}{ll} U(x_1, \dots, x_{i-1},0, x_{i+1}, \dots,
x_d) + 2 i c_k \le U(x) +
2 i c_k \\[6pt]
U(x_1, \dots, x_{i-1},1, x_{i+1}, \dots, x_d) - 2 i c_k \ge U(x) -
2 i c_k,
\end{array}\right. \quad \mbox{ for } 1 \le i \le d.  
\end{equation*}
It follows from \eqref{defgk} that, for $1 \le i \le d$, 
$$
\min \Big(g_{i-1, k}(x), U(x)-2 i c_k \Big) \le g_{i,k}(x) \le \max \Big(g_{i-1, k}(x), U(x)+2 i c_k \Big).
$$
Using the fact $U(x) - c_k \le g_{0,k}(x) \le U(x) + c_k$ by \eqref{defck}, we obtain, for $1 \le i \le d$, 
\begin{equation}\label{TT1}
U(x)-2 i c_k \le g_{i,k}(x) \le U(x)+2 i c_k. 
\end{equation}
Since $\lim_{k \to + \infty} c_k = 0$, it follows from \eqref{defgk} and \eqref{TT1} that, for large $k$, 
\begin{equation}\label{propertygk}
g_{i, k}(x) = \left \{\begin{array}{ll} U(x_1, \dots, x_{i-1},0,
x_{i+1}, \dots, x_d) + 2 i c_k & \mbox{if } 0 \le x_i \le c_k,\\[6pt]
U(x_1, \dots, x_{i-1},1, x_{i+1}, \dots, x_d) - 2 i c_k, &
\mbox{if } 1 - c_k \le x_i \le 1. 
\end{array}\right.
\end{equation}
We derive from \eqref{defgk} and \eqref{propertygk} that $g_{d,k}$ is Lipschitz on $Q_{c_k}$ with a
Lipschitz constant $1$ ($=| \nabla U|$). We claim that 
\begin{equation}\label{claim-B}
\limsup_{k \to \infty} \Lambda_{\delta_k} (g_{i, k}, Q)
\le \limsup_{k \to \infty} \Lambda_{\delta_k} (g_{i-1, k}, Q)
\end{equation}
for all $1 \le i \le d$. 

We establish \eqref{claim-B} for $i=1$ (the argument is the same for every $i$). We first apply Lemma~\ref{lemtechnical1} with $A = Q$,  $f(x) = \max\big(g_{0, k}(x), U(0, x_2, \cdots, x_d) + 2 c_k\big)$, and $g(x) = U(1, x_2, \cdots, x_d) - 2 c_k$. Recall that 
$$
Q \setminus B_1 = \big\{ x \in Q; f(x) > g(x) \big\}. 
$$
Note that 
$$
f(x) \le \max \Big( U(x) + c_k, U(0,  x_2, \cdots, x_d) + 2 c_k\Big)  \le U(x) + 2 c_k. 
$$
It follows that if $x \in Q \setminus B_1$ then $U(x) + 2 c_k > U(1, x_2, \cdots, x_d) - 2 c_k$; this implies $1 - x_1 < 4 \sqrt{d} c_k $. Hence $|Q \setminus B_1| \le C c_k$ and it follows from Lemma~\ref{lemtechnical1} that
\begin{equation}\label{claim-B1}
\Lambda_{\delta_k} (g_{1, k}, Q) \le \Lambda_{\delta_k} \Big(\max\big(g_{0, k}(x), U(0, x_2, \cdots, x_d) + 2 c_k\big), Q \Big)  + C c_k. 
\end{equation}
We next apply Lemma~\ref{lemtechnical1} with $A = Q$,  $f(x) = g_{0, k}(x)$, $g(x) = U(0, x_2, \cdots, x_d) + 2 c_k$, and $B_2 = \{ x \in Q; f(x) > g(x) \}$. If $x \in Q \setminus B_2$ we have $U(0, x_2, \cdots, x_d) + 2 c_k < g_{0, k} $ so that $U(0, x_2, \cdots, x_d)  + 2 c_k <  U(x)  - c_k$; this implies $x_1 <  3 \sqrt{d} c_k $.  Hence $|Q \setminus B_2| \le C c_k$ and it follows from Lemma~\ref{lemtechnical1} that
\begin{equation}\label{claim-B2}
 \Lambda_{\delta_k} \Big(\max\big(g_{0, k}(x), U(0, x_2, \cdots, x_d) + 2 c_k\big), Q \Big)  \le \Lambda_{\delta_k} (g_{0, k},  Q \Big)  + C c_k. 
\end{equation}
Claim~\eqref{claim-B} now follows from \eqref{claim-B1} and \eqref{claim-B2} since $\lim_{k \to + \infty} c_k = 0$.

From \eqref{claim-B}, we deduce that 
\begin{equation*}
\limsup_{k \to \infty} \Lambda_{\delta_k} (g_{d, k}, Q)
\le \limsup_{k \to \infty} \Lambda_{\delta_k} (g_{k}, Q). 
\end{equation*}
The conclusion follows by choosing $h_k = g_{d, k}$.  \proofend

\medskip
We next establish the following lemma which plays an important role in the proof of Property (G2).

\begin{lemma}\label{lemcorlemfondN} There exist  
$(c_\delta)\subset \mR_+$ and $(g_\delta)  \subset L^\infty(Q)$  such that 
$$
c_\delta \ge \sqrt{\delta}, \quad  \lim_{\delta \to 0 } c_\delta = 0,
$$ 
$$
\|g_\delta  - U \|_{L^\infty(Q)} \le 2 d c_\delta,  \quad \mathrm{Lip}(g_\delta, Q_{c_\delta}) \le 1, 
$$
and
\begin{equation*}
\lim_{\delta \to 0} \Lambda_\delta(g_\delta, Q) \le
\C.
\end{equation*}
\end{lemma}

\noindent{\bf Proof.}  Applying Lemma~\ref{lemcorlemfondN-111}, we derive  from the definition of $\C$ that there exist a sequence $(g_k) \subset L^\infty(Q)$ and two sequences $(\delta_k), (c_k) \subset \mR_+$ such that 
\begin{equation}\label{delta-c-k}
 \lim_{k \to + \infty} \delta_k =  \lim_{k \to + \infty} c_k = 0, \quad c_k \ge  \sqrt{\delta_k},
\end{equation}
\begin{equation}\label{gk-infty}
\| g_k - U \|_{L^\infty(Q)} \le 2 d c_k, \quad \mathrm{Lip}(g_k, Q_{c_k}) \le 1,
\end{equation}
and 
\begin{equation}\label{estlemfondN.1}
 \lims{k \to + \infty} \Lambda_{\delta_k} (g_k, Q) \le \C.  
\end{equation}
We next  construct a family $(h_\delta) \subset L^\infty(Q)$ such that
\begin{equation}\label{haha-B}
\|h_\delta - U \|_{L^\infty(Q)} \to 0 \quad \mbox{ and } \quad \limsup_{\delta \to 0} \Lambda_\delta(h_\delta, Q) \le \C.
\end{equation}
Let $(\tau_k)$ be a strictly decreasing positive sequence such
that $\tau_k \le c_k \delta_{k}$.
For each $\delta$ small, let $k$ be such that $\tau_{k+1} < \delta \le \tau_{k}$ and  define $
m_1=\delta_{k}/\delta \ge 1/ c_k$  and $m = [m_1]$. As usual, for $a > 0$,  $[a]$ denotes the largest
integer $ \le a$. Define $h^{(1)}_{\delta}: [0, m]^d \to \mR$  as follows
\begin{equation}\label{defh1delta}
h^{(1)}_{\delta} (y) =   \frac{\sum_{i=1}^{d} [y_i]}{\sqrt{d}} + g_{k}(x) \mbox{ with } x= (y_1 - [y_1], \cdots, y_d - [y_d]).  
\end{equation}
For $\alpha \in \mN^d$ and $c \ge 0$, set
\begin{equation*}
Q_{+\alpha} := Q + (\alpha_1, \cdots ,\alpha_d), \quad Q_{+\alpha, c} := Q_c + (\alpha_1, \cdots ,\alpha_d),
\quad \mbox{and} \quad D_{+\alpha, c} := Q_{+\alpha} \setminus Q_{+\alpha, c}. 
\end{equation*}
Define
$$
Y = \mN^d \cap [0, m-1]^d \quad \mbox{ and } \quad B = \bigcup_{\alpha \in Y} \Big(Q_{+\alpha, c_k} \setminus Q_{+\alpha, c_k/2} \Big).
$$
We claim that 
\begin{equation}\label{Liph1delta}
\mathrm{Lip}(h^{(1)}_\delta, B) \le C. 
\end{equation}
Indeed, it is clear from \eqref{gk-infty} and  \eqref{defh1delta} that 
$$
\mathrm{Lip}(h^{(1)}_\delta, Q_{+\alpha, c_k} \setminus Q_{+\alpha, c_k/2} ) \le 1 \mbox{ for } \alpha \in Y. 
$$
On the other hand, if $y \in Q_{+\alpha, c_k} \setminus Q_{+\alpha, c_k/2}$
and $y' \in Q_{+\alpha', c_k} \setminus Q_{+\alpha', c_k/2}$ with $\alpha \neq \alpha'$ then $c_k \le C |y - y'|$ so that 
\begin{align*}
|h_\delta^{(1)}(y) - h_\delta^{(1)}(y')|  & \le  |h_\delta^{(1)}(y) - U(y)| +  |h_\delta^{(1)}(y') - U(y')| + |U(y) - U(y')| \\[6pt]
& \mathop{\le}^{ \mbox{ by } \eqref{gk-infty}}  |U(y) - U(y')| + 4 d c_k \le  |y-y'|  + 4 d c_k \le C |y -y'| .
\end{align*}
Claim \eqref{Liph1delta} follows. 

By  classical  Lipschitz extension it follows from \eqref{Liph1delta} that  there exists $h^{(2)}_{\delta}: \mR^d \mapsto \mR$ such that $h^{(2)}_{\delta} = h^{(1)}_{\delta}$ on $B$ and
\begin{equation}\label{Liph2delta}
\mathrm{Lip}(h^{(2)}_{\delta}, \mR^d) \le C.
\end{equation}
Define, for $x \in \mR^d$, 
\begin{equation}\label{defh3delta}
h^{(3)}_{\delta} (x) =
\left\{ \begin{array}{ll}
h^{(1)}_{\delta}(x) &\mbox{ if } x \in D_{+\alpha, c_k/2} \mbox{ for some } \alpha \in Y, \\[6pt]
h^{(2)}_{\delta}(x) & \mbox{ otherwise},
\end{array}\right.
\end{equation}
and set
\begin{equation*}
h_\delta(x) = \frac{1}{m_1} h^{(3)}_{\delta} (m x) \mbox{ in } [0, 1]^d. 
\end{equation*}
Since $\delta = \delta_k/ m_1$, by a change of variables, we obtain
\begin{equation}\label{estlemfondN.2}
\Lambda_{\delta} (h_\delta, Q) = \frac{\delta}{\delta_k} \Lambda_{\delta_k} \big(h_\delta^{(3)} (m \, \cdot  \, ), Q \big)
= \frac{m^{1-d}}{
m_1} \Lambda_{\delta_k} \big(h_\delta^{(3)}, [0, m]^d \big). 
\end{equation} 
We next estimate $\Lambda_{\delta_k} \big(h_\delta^{(3)}, [0, m]^d \big)$. 
For $\alpha \in Y$,  applying   Lemma~\ref{lem-computation} with $c = c_k$, $D = Q_{+\alpha}$ and $g = h^{(3)}_\delta$, we have
\begin{equation}\label{mama-1}
\mathop{\iint}_{Q_{+\alpha} \times [0, m]^d} \frac{\varphi_\delta(|h^{(3)}_{\delta} (x) -  h^{(3)}_{\delta} (y)|) }{|x-y|^{d+1}} \, dx \, dy  \le \Lambda_{\delta} \big(h^{(3)}_{\delta}, D_{+\alpha, c_k/2}\big) + C(c_k + b \delta / c_k ).  
\end{equation}
From \eqref{defh1delta} and \eqref{defh3delta},  we obtain
\begin{equation}\label{h3-delta-11}
\Lambda_{\delta}(h_\delta^{(3)}, D_{+\alpha, c_k/2}) =  \Lambda_{\delta} (g_k, D_{+(0, \cdots, 0), c_k/2}) \le \Lambda_{\delta}(g_k, Q). 
\end{equation}
Since 
$$
\int_{[0, m]^d} \int_{[0, m]^d} \cdots= \sum_{\alpha \in Y} \int_{Q_\alpha} \int_{[0, m]^d} \cdots, 
$$
it follows from \eqref{mama-1}  and \eqref{h3-delta-11} that  
\begin{equation}\label{estlemfondN.3}
\Lambda_{\delta_k} (h_\delta^{(3)}, [0, m]^d)  \le   m^d \Lambda_{\delta_k} (g_{k}, Q)
+  Cm^d (c_k + b \delta_k/ c_k).
\end{equation}
Since $m \le m_1$ and $c_k \ge \delta_k^{1/2}$ by \eqref{delta-c-k}, we deduce from \eqref{estlemfondN.2} and
\eqref{estlemfondN.3} that
\begin{equation}\label{To1}
\Lambda_{\delta}(h_\delta, Q)   \le
\Lambda_{\delta} (g_{k}, Q)+ C (c_k +
b \delta_{k}^\frac{1}{2}).
\end{equation}
Combining \eqref{delta-c-k}, \eqref{estlemfondN.1},   and \eqref{To1} yields 
\begin{equation*}
\limsup_{\delta \to 0}\Lambda_{\delta}(h_\delta, Q)   \le \C. 
\end{equation*}

We next claim that 
\begin{equation}\label{claim-mm1}
\| h_3^{(\delta)} - U \|_{L^\infty([0, m]^d)} \le C c_k. 
\end{equation}
Indeed, for $y \in [0, m]^d$, we have, by \eqref{defh1delta}, 
$$
|h_1^{(\delta)}(y) - U(y) |= |g_k(x) - U(x)|  \mbox{ where } x= (y_1 - [y_1], \cdots, y_d - [y_d]). 
$$ 
It follows from \eqref{gk-infty} that 
\begin{equation} \label{claim-mm2}
\|h_1^{(\delta)} - U \|_{L^\infty([0, m]^d)} \le C c_k. 
\end{equation}
On the other hand, for  $y \in [0, m]^d \setminus  \dsp \bigcup_{\alpha \in Y} D_{+\alpha, c_k/2}$, let $\dsp \hat y \in B$ such that $|y - \hat y| \le c_k$. Since $h_\delta^{(2)} (\hat y) = h_\delta^{(1)} (\hat y)$,  it follows from \eqref{Liph2delta} and \eqref{claim-mm1} that 
\begin{equation} \label{claim-mm3}
|h_2^{(\delta)} (y) - U(y)| \le |h_2^{(\delta)}(y) - h_2^{(\delta)} (\hat y)|+  |h_1^{(\delta)} (\hat y) - U(\hat y)| + |U(\hat y) - U(y)| \le C c_k  
\end{equation}
Claim~\eqref{claim-mm1} now follows from \eqref{claim-mm2} and  \eqref{claim-mm3}.

Using \eqref{claim-mm1}, we derive from the definition of $h_\delta$ and the facts that $m_1 \ge 1/ c_k$ and  $c_k \to 0$ that $\| h_\delta - U\|_{L^\infty(Q)} \to 0$. Hence \eqref{haha-B} is established. 

\medskip
The conclusion now follows from \eqref{haha-B} and Lemma~\ref{lemcorlemfondN-111}.  \proofend

\medskip 

We next establish 
\begin{lemma}\label{lemlimsupaffine1}
Let $S$ be an open bounded subset of $\mR^d$ with  Lipschitz boundary and let  $g$ be an affine function defined on $S$. Then
\begin{equation}\label{k-S}
\inf  \limf{\delta \lr 0}\Lambda_\delta(g_\delta, S)  =  \C |\nabla g| |S|,
\end{equation}
where the infimum is taken over all families  $(g_\delta)_{\delta \in (0,1)} \subset L^1(S)$  such
that $g_\delta \to g$ in $L^1(S)$. Moreover, there
exists a family $(h_\delta) \subset L^\infty(S)$ such that $\|h_\delta - g\|_{L^\infty(S)} \to 0$ and
\begin{equation*}
\lim_{\delta \lr 0} \Lambda_\delta (h_\delta, S) = \C |\nabla g| |S|.
\end{equation*}
\end{lemma}

\noindent {\bf Proof.} Note that if $T: \mR^d \to \mR^d$ is an affine conformal transformation, i.e., 
$$
T(x) = a R x + b \mbox{ in } \mR^d 
$$
for some $a > 0$, some linear unitary operator $R: \mR^d \to \mR^d$, and for some $b \in \mR^d$,  then, for a measurable subset $D$ of $\mR^d$ and $f \in L^1(D)$,   
\begin{equation*}
\Lambda_{\delta} (f, D) = a^{1 - d} \Lambda_{\delta} (f \circ T^{-1}, T(D)), 
\end{equation*}
by a change of variables. 

Using a transformation $T$ as above, we may write $
g \circ T^{-1} = U$.  Then 
$$
\Lambda_{\delta} (g_\delta, S) = a^{1-d} \Lambda_\delta(g_\delta \circ T^{-1}, T(S))
$$
and 
$$
|\nabla g| |S| = a^{1 - d} |T(S)|
$$
Hence, it suffices to prove Lemma~\ref{lemlimsupaffine1} for  $g = U$. 

Denote $m$ the LHS of \eqref{k-S}. Since $|\nabla g| = |\nabla U| =1$, \eqref{k-S} becomes
\begin{equation}\label{mS}
m = \C  |S|. 
\end{equation}
The proof  of  \eqref{mS} is based on a covering lemma \cite[Lemma 3]{NgGamma} (applied first with $\Omega = S$ and $B = Q$ and then with $\Omega = Q$ and $B = S$) which asserts that 
\begin{enumerate}
\item[$i)$] There exists a
sequence of disjoint sets $(Q_k)_{k \in \mN}$ such that $Q_k$ is the  image
of $Q$ by a dilation and a translation,  $Q_k \subset S$ for all $k$, and
$$
|S| = \sum_{k \in \mN} |Q_k|.
$$
\item[$ii)$] There exists a
sequence of  disjoint sets $(S_{k})_{k \in \mN}$ such that $S_{k}$ is
the image of $S$ by a dilatation and a translation, $S_{k} \subset Q$ for all $k$, and
$$
|Q| = \sum_{k \in \mN} |S_{k}|.
$$
\end{enumerate}
We first claim that 
$$
m \ge \C |S|. 
$$ 
Indeed, let $(Q_k)$ be the sequence of disjoint sets in  $i)$.

Clearly, 
\begin{equation}\label{p1}
 \Lambda_{\delta} (g_\delta, S) \ge  \sum_{k \in \mN} \Lambda_{\delta} (g_\delta, Q_k). 
\end{equation}
Fix $k \in \mN$ and let $a_k  > 0$ and $b_k \in \mR^d$ be such that $Q_k = a_k Q + b_k$. Then $|Q_k| = a_k^d$ and, by a change of variables,  
\begin{equation}\label{f1}
\Lambda_\delta (g_\delta, a_k Q + b_k) = a_k^d \Lambda_{\delta/ a_k} (\hat g_\delta, Q) \mbox{ where } \hat g_\delta (x) = \frac{1}{a_k} g_\delta (a_k x + b_k). 
\end{equation}
From the definition of $\C$, we have
\begin{equation}\label{f2}
\limf{\delta \lr 0} \Lambda_{\delta/ a_k} (\hat g_\delta, Q)
\ge \C. 
\end{equation}
We deduce from \eqref{f1} and \eqref{f2} that 
\begin{equation}\label{p2}
 \limf{\delta \lr 0}  \Lambda_{\delta} (g_\delta, Q_k) \ge \C |Q_k|. 
\end{equation}
Combining \eqref{p1} and \eqref{p2} yields 
$$
\liminf_{\delta \to 0} \Lambda_{\delta} (g_\delta, S) \ge \C |S|; 
$$
which implies $m \ge \C |S|$. Similarly, using $ii)$ one can show that $
\C |S| \ge m.$ We thus obtain   \eqref{mS}. 

It remains to prove that there exists a family $(h_\delta)$  such that $\|h_\delta - g \|_{L^\infty(S)} \to 0$ and  
\begin{equation}\label{lemlimsupaffine1-hoho}
\lim_{\delta \to 0} \Lambda_\delta (h_\delta, S) = \C |\nabla g| |S|.
\end{equation}
As above, we can assume that $g = U$. Let $\hat Q$ be the image of $Q$ by a dilatation and a translation such that  $S \subset \subset \hat Q$. By Lemma~\ref{lemcorlemfondN} and a change of variables, there exists a family $(h_\delta) $ such that $h_\delta \to U$ in $L^1(\hat Q)$  and 
$$
\lim_{\delta \to 0} \Lambda_{\delta}(h_\delta, \hat Q) = \C |\hat Q|. 
$$
On the other hand, we have, by \eqref{k-S},  
$$
\liminf_{\delta \to 0} \Lambda_{\delta}(h_\delta, \hat Q \setminus  S) \ge \C |\hat Q \setminus S|. 
$$
Moreover,  
$$
\Lambda_\delta(h_\delta, \hat Q) \ge \Lambda_\delta(h_\delta, S) + \Lambda_\delta(h_\delta, \hat Q \setminus S). 
$$
It follows that 
$$
\limsup_{\delta \to 0} \Lambda_{\delta}(h_\delta, S) \le \C |S|, 
$$
which implies \eqref{lemlimsupaffine1-hoho} by \eqref{k-S}.  \proofend

\medskip

 
Throughout the rest of Section~\ref{sect-proof-G2}, we 
let $A_1$, $A_2$, \dots, $A_m$ be disjoint open $(d+1)$-simplices
in $\mR^d$ such that every coordinate component of any vertex of
$A_i$ is equal to 0 or 1, 
\begin{equation*}
\bar Q= \bigcup_{i =1}^m \bar A_i,
\end{equation*}
and
\begin{equation*}
A_1 = \Big\{x= (x_1, \dots, x_d) \in \mR^d; \, x_i > 0 \mbox{ for
all } 1 \le i \le d, \mbox{ and }\sum_{i=1}^d x_i < 1 \Big\}.
\end{equation*}
We also denote $A_{\ell, c}$ the set $(A_\ell)_c$ (see \eqref{defAc}).

\medskip 
The following lemma is a variant of Lemma~\ref{lemcorlemfondN} for $\{ A_\ell \}_{\ell =1}^m$.

\begin{lemma}\label{lemfondN.1} Let $ \ell \in \{1, \dots, m \}$ and $g$ be an affine function defined on
$A_\ell$ such that its normal derivative $\frac{\partial g}{\partial n} \neq 0$ along the
boundary of $A_\ell$, where $n$ denotes the inward normal. There exist a family $(g_\delta)  \subset L^\infty(A_\ell)$  and a
family $(c_\delta) \subset \mR_+$ such that 
$$
c_\delta \ge \sqrt{\delta}, \quad \lim_{\delta \to 0} c_{\delta} = 0,
$$
$$
\|g_\delta - g \|_{L^\infty(A_\ell)} \le C_d |\nabla g|  c_\delta,  \quad \mathrm{Lip}(g_\delta, A_{\ell, c_\delta}) \le |\nabla g|, 
$$
and
\begin{equation*}
\lims{\delta \lr 0} \Lambda_{\delta} (g_\delta, A_\ell) \le \C |\nabla g| |A_\ell|.
\end{equation*}
\end{lemma}

\medskip 

\noindent{\bf Proof.} For notational ease,  we assume that  $\ell =1$. The proof is in the spirit of  the one of Lemma~\ref{lemcorlemfondN-111}.
By Lemma~\ref{lemlimsupaffine1}, there
exists a family  $(h_\delta) \subset L^\infty(A_1)$ such that $\| h_\delta - g\|_{L^\infty(A_1)}  \to  0 $ and
\begin{equation}\label{huhuhuhu}
\lim_{\delta \lr 0} \Lambda_\delta (h_\delta, A_1) = \C |\nabla g| |A_1|.
\end{equation}
Set
$$
c_\delta = \max \big(\|h_\delta - g \|_{L^\infty(A_1)}, \sqrt{\delta} \big)
\quad \mbox{ and } \quad l_\delta = 2  |\nabla g| c_\delta. 
$$
Denote $h_{0, \delta} = h_\delta$, and define,  for $i$ = 1, 2, \dots, $d$, and $x \in A_1$, 
\begin{equation}\label{defg-i}
h_{i, \delta}(x) = \left\{\begin{array}{ll} \max
\big(h_{i-1,\delta}(x),  g(x_1, \dots, x_{i-1}, 0,  x_{i+1},
\dots, x_d) +  i l_\delta  \big) \mbox{ if } \frac{\partial g}{\partial x_i} > 0,\\[6pt]
\min \big(h_{i-1,\delta}(x),  g(x_1, \dots, x_{i-1}, 0,  x_{i+1},
\dots, x_d) -  i l_\delta \big) \; \mbox{ if } \frac{\partial
g}{\partial x_i} < 0.
\end{array}\right.
\end{equation}
Set $e = (\frac{1}{\sqrt{d}}, \dots, \frac{1}{\sqrt{d}})$ and
define, for $x \in A_1$, 
\begin{equation}\label{defg-N+1}
h_{d+1, \delta}(x) = \left\{\begin{array}{ll} \max
\big(h_{d,\delta}(x),  g(z(x)) +  (d+1) l_\delta \big) & \mbox{if } \frac{\partial g}{\partial e} < 0,\\[6pt]
\min \big(h_{d,\delta}(x),  g(z(x)) - (d+1) l_\delta \big) &
\mbox{if } \frac{\partial g}{\partial e} > 0.
\end{array}\right.
\end{equation}
Here for each $x \in A_1$, $z(x) := x - \langle x, e \rangle e +
e$ (the projection of $x$ on the hyperplane
$P$ which is orthogonal to $e$ and contains $e$). 

As in the proof of Lemma~\ref{lemcorlemfondN-111}, we have the following three assertions, for $x \in A_1$,  \\
i)  for $1 \le i \le d+1$,    
$$
g(x) -  i \ell_\delta \le h_{i, \delta} (x) \le g(x) +  i \ell_\delta, 
$$
ii) for $1 \le i \le d$ and  $0 \le  x_i \le c_\delta$, 
\begin{equation*}
h_{i, \delta}(x) = \left\{\begin{array}{ll}  g(x_1, \dots, x_{i-1}, 0,  x_{i+1},
\dots, x_d) +  i l_\delta   \mbox{ if } \frac{\partial g}{\partial x_i} > 0,\\[6pt]
g(x_1, \dots, x_{i-1}, 0,  x_{i+1},
\dots, x_d) -  i l_\delta  \; \mbox{ if } \frac{\partial
g}{\partial x_i} < 0, 
\end{array}\right.
\end{equation*}
iii) for $|x - z(x)| \le c_\delta$, 
\begin{equation*}
h_{d+1, \delta}(x) = \left\{\begin{array}{ll}   g(z(x)) +  (d+1) l_\delta  & \mbox{if } \frac{\partial g}{\partial e} < 0,\\[6pt]
g(z(x)) - (d+1) l_\delta   &
\mbox{if } \frac{\partial g}{\partial e} > 0.
\end{array}\right.
\end{equation*}
It follows that $h_{d+1, \delta}$ is Lipschitz on $A_{1, {c_\delta}}$ with Lipschitz constant $|\nabla g|$. As in the proof of Lemma~\ref{lemcorlemfondN-111}, one has, for $0 \le i \le d$, 
$$
\limsup_{\delta \to 0} \Lambda_{\delta} (h_{i+1, \delta}, A_1) \le \limsup_{\delta \to 0} \Lambda_{\delta} (h_{i, \delta}, A_1); 
$$
which implies, by \eqref{huhuhuhu},  
$$
\limsup_{\delta \to 0} \Lambda_{\delta} (h_{d+1, \delta}, A_1) \le \limsup_{\delta \to 0} \Lambda_{\delta} (h_{0, \delta}, A_1)  =  \limsup_{\delta \to 0} \Lambda_{\delta} (h_{ \delta}, A_1)  = \C |\nabla g| |A_1|.  
$$
The conclusion now holds for $g_\delta = h_{d+1, \delta}$.  \proofend

\medskip 
We end this section with the following result which is a consequence of Lemma~\ref{lemfondN.1} by a change of variables. 
\begin{definition} \fontfamily{m} \selectfont
For each $k \in \mN$, a set $K$ is called a $k$-sim of $\mR^d$ if
there exist $z \in \mZ^d$ and $ \ell \in \{ 1, 2, \dots,
m \}$ such that $K = \frac{1}{2^k} A_\ell + \frac{z}{2^k}$.
\end{definition}

We have

\begin{corollary}\label{corfondN.2} Let $K$ be a $k$-sim of $\mR^d$ and $g$ be an affine function defined on
$K$ such that $\frac{\partial g}{\partial n} \neq 0$ along the
boundary of $K$. There exist a family $(g_\delta) \subset L^\infty(K)$ and a
family  $(c_\delta) \subset \mR_+$ such that 
$$
c_\delta \ge C_k \sqrt{\delta}, \quad \lim_{\delta \to 0} c_\delta = 0, 
$$
$$
\|g_\delta - g\|_{L^\infty(K)} \le C_k |\nabla g| c_\delta, \quad \mathrm{Lip}(g_\delta, K_{ c_\delta}) \le C_k|\nabla g|,
$$
and
\begin{equation*}
\lims{\delta \lr 0} \Lambda_\delta (g_\delta, K)  \le \C |\nabla g| |K|.
\end{equation*}
\end{corollary}

In Corollary~\ref{corfondN.2} and Section~\ref{sect-pro-proG2} below, $C_k$ denotes a positive constant depending only on $k$ and $d$ and can be different from one place to another. 

\subsubsection{Proof of Property (G2)}\label{sect-pro-proG2}

Our goal is to show that (G2) holds with $K = \C$, i.e.,  
\begin{enumerate}
\item[(G2)] For every $u \in L^1(\Omega)$, there exists a family
$(\tilde u_\delta)_{\delta \in (0,1)} \subset L^1(\Omega)$ such that
$\tilde u_\delta \to u$ in $L^1(\Omega)$ as $\delta \to 0$,
and
\begin{equation*}
\limsup_{\delta \to 0}  \Lambda_\delta(\tilde u_{\delta}) \le\C \int_{\Omega} |\nabla u|. 
\end{equation*}
\end{enumerate}

 We consider the case $\Omega = \mR^d$ and the case where $\Omega$ is bounded separately. 
  
  \medskip 
\noindent {\bf Case 1: $\Omega = \mR^d$.} The proof is divided into two steps. Given $k \in \mN$, set 
\begin{equation}
R_k: = \left\{ u \in C^0_{c}(\mR^d) 
\left|\begin{array}{l} u \mbox{ is affine on each $k$-sim and $\partial u/ \partial n \neq 0$ along the } \\[6pt] 
\mbox{boundary of each $k$-sim,  unless $u$ is constant on
that $k$-sim}
\end{array} \right. 
 \right\}. 
\end{equation}

{\bf Step 1}.   We prove Property (G2) when $u \in R_k$ and $k \in \mN$ is arbitrary but fixed. Set 
$$
{\cal K} = \{K \mbox{ is a } \mbox{$k$-sim  and }   u \mbox{ is not constant on } K \}. 
$$
From now on in the proof of Step 1, $K$ denotes a $k$-sim.  By  Corollary~\ref{corfondN.2}, for each $K \in {\cal K}$, there exist 
$(u_{K,\delta}) \subset L^\infty(K)$ and $(c_{K,
\delta}) \subset \mR_+$ such that 
\begin{equation}\label{G2-1}
c_{K, \delta } \ge C_k \sqrt{\delta}, \quad \lim_{\delta \to 0} c_{K, \delta} = 0, 
\end{equation}
\begin{equation}\label{G2-2}
\| u_{K,  \delta}- u \|_{L^\infty(K)} \le  C_k \|\nabla u\|_{L^\infty(\mR^d)} c_{K, \delta}, \quad \mathrm{Lip}(u_{K, \delta}, K_{ c_{K,\delta}})\le C_k \|\nabla u\|_{L^\infty(\mR^d)},
\end{equation}
and
\begin{equation}\label{part2-haha}
\lims{\delta \lr 0} \Lambda_{\delta}(u_{K, \delta}, K)
 \le \C  \int_{K} |\nabla u| \, dx.
\end{equation}
For each $\delta$, let  $u_{\delta}$ be a function defined in $\mR^d$ such that
\begin{equation}\label{part3-haha}
u_{\delta} = u_{K, \delta} \mbox{ in } K \setminus K_{c_{K, \delta}/2} \mbox{ for } K \in {\cal K},  
\quad
u_{\delta} =  u \mbox{ in } K \mbox{ for } K \not \in {\cal K},  
\end{equation}
and 
\begin{equation}\label{part4-haha}
|\nabla u_{\delta}(x)| \le C_k \| \nabla u\|_{L^\infty(\mR^d)} \quad \mbox{ for } x \in \mR^d \setminus \bigcup_{K \in {\cal K}} (K \setminus K_{c_{K, \delta}/2}). 
\end{equation}
Such a $u_\delta$ exists  by  \eqref{G2-2} via standard Lipschitz extension.  Applying Lemma~\ref{lem-computation} with $D = K$ and $g = u_\delta$, we have, by \eqref{G2-2},   
\begin{equation}\label{hohoho}
\mathop{\iint}_{K \times \mR^d} \frac{\varphi_{\delta} (|u_\delta(x) - u_\delta(y)|) } { |x - y|^{d+1} } \, dx \, dy \le \Lambda_{\delta} \big(u_\delta, K \setminus K_{c_{K, \delta}/2} \big) + C_k  (\| \nabla u\|_{L^\infty(\mR^d)}c_{K, \delta} + b \delta/ c_{K, \delta}). 
\end{equation}
From the definition of $u_\delta$, there exists $R>1$, independent of $\delta$, such that  $u_\delta = u =  0$ in $\mR^d \setminus B_{R}$.  We have, for some $b > 0$  (see \eqref{cond-varphi-01}),  
\begin{equation}\label{hohoho-1}
\mathop{\iint}_{(\mR^{d} \setminus B_{R+1}) \times \mR^d} \frac{\varphi_{\delta} (|u_\delta(x) - u_\delta(y)|) } { |x - y|^{d+1} } \, dx \, dy \le \mathop{\iint}_{  B_{R} \times (\mR^d \setminus B_{R+1}) } \frac{ \delta b} { |x - y|^{d+1} } \, dy \, dx \le C_d R^d \delta b. 
\end{equation}
Combining  \eqref{part2-haha},  \eqref{part3-haha}, \eqref{hohoho}, and \eqref{hohoho-1} yields  
\begin{equation*}
\lims{\delta \lr 0} \Lambda_{\delta}(u_\delta, \mR^d)
 \le \C  \int_{\mR^d} |\nabla u| \, dx. 
\end{equation*}

We next claim that  $u_{\delta} \to u$ in $L^1(\mR^d)$.
Indeed, for $x \in K_{c_{K, \delta}/2}$ for some $K \in {\cal K}$, let $\hat x \in K_{c_{K, \delta}} \setminus K_{c_{K, \delta}/2}$ be such that $|\hat x - x| \le c_{K, \delta}$. We have
\begin{equation*} 
|u_\delta (x) - u(x)|  \le |u_\delta(x) - u_\delta(\hat x)| + |u_\delta(\hat x) - u (\hat x)| + |u(\hat x) - u(x)| \le C_k \| \nabla u \|_{L^\infty(\mR^d)} c_{K, \delta}. 
\end{equation*}
This implies, for $K \in {\cal K}$,  
$$
\lim_{\delta \to 0}\|u_\delta - u \|_{L^\infty(K)} = 0, 
$$
Since $u_\delta = u$ in $K$ for $K \not \in {\cal K}$, we deduce that 
$$
\lim_{\delta \to 0}\|u_\delta - u \|_{L^1(\mR^d)} = 0. 
$$

The proof of Step 1 is complete.

\medskip 
{\bf Step 2.} We prove Property (G2) for a general $u \in L^1(\mR^d)$. Without loss of generality, one may assume that $u \in BV(\mR^d)$ since there is nothing to prove otherwise.  Let $(u_n) \subset C^\infty_{\mc}(\mR^d)$ be  such
that $u_n$ converges to $u$ in $L^1(\mR^d)$ and $\| \nabla u_n \|_{L^1(\mR^d)} \to \int_{\mR^d} |\nabla u|$ as $n \to + \infty$.  We next use

\begin{lemma}\label{lem-EF} Let $v \in C^1_{\mc}(\mR^d)$ with $\supp v \subset B_R$ for some $R>0$. There exists a sequence  $(v_m) \subset W^{1, \infty}(\mR^d)$ such that  $v_m \subset R_m$, $\supp v_m \subset B_R$ for large $m$ and  $v_m \to v$  in $W^{1, 1}(\mR^d)$ as $m \to + \infty$. 
\end{lemma}

\noindent{\bf Proof of Lemma~\ref{lem-EF}.} There exist a sequence $(k_m) \subset \mN$  and a sequence $(v_m) \subset W^{1, \infty}(\mR^d)$ with $\supp v_m \subset B_R$ for large $m$ such that \\
i) $v_m \to v$  in $W^{1, 1}(\mR^d)$ as $m \to + \infty$. \\
ii) $v_m$ is affine  on each  $m$-sim. \\
This fact is standard in finite element theory, see e.g.,   \cite[Proposition 6.3.16]{Allaire}.  By a small perturbation of $v_m$, one can also assume that $\partial v_m/ \partial n \neq 0$ along the boundary of each  $m$-sim, unless $v_m$ is constant there. \proofend

\medskip 
We now return to the proof of Step 2. By Lemma~\ref{lem-EF}, for each $n \in \mN$, there exists $v_n  \in R_k$ for some $k \in \mN$ such that 
$$
\| v_n - u_n\|_{W^{1, 1}(\mR^d)} \le 1/ n. 
$$
By Step 1, there exists a family 
$(v_{\delta, n}) \subset L^1(\mR^d)$ such
that $v_{\delta, n}$ converges to $v_{n}$ in $L^1(\mR^d)$ as $\delta \to 0$ and 
\begin{equation*}
\lims{\delta \lr 0} \Lambda_{\delta} (v_{\delta, n}) \le \C
\int_{\mR^d} |\nabla  v_{n}| \, dx.
\end{equation*}
Hence there exists $\delta_n > 0$ such that, for $0 < \delta < \delta_n$ 
$$
\Lambda_{\delta}(v_{\delta, n}) \le \C \int_{\mR^d} |\nabla v_n| + 1/ n \quad \mbox{ and } \quad \|v_{\delta, n}  - v_n \|_{L^1(\mR^d)} \le 1/ n. 
$$
Without loss of generality, one may assume that $(\delta_n)$ is decreasing to 0. Set
$$
u_{\delta} = v_{\delta_{n+1}, n} \quad  \mbox{ for } \delta_{n+1} \le \delta < \delta_n. 
$$
Then $(u_\delta)$ satisfies the properties required. 

\medskip 
The proof of Case 1 is complete. \proofend

\bigskip
\noindent{\bf Case 2: $\Omega$ is bounded.}  We prove Property (G2) for a general $u \in L^1(\Omega)$. Without loss of generality, one may assume that $u \in BV(\Omega)$. Let $R> 0$ be such that $\Omega \subset \subset B_R$ and let  $(u_n) \subset C^\infty(\mR^d)$ with $\supp u_n \subset B_R$ such that $u_n \to u$ in $L^1(\Omega)$ and $\| \nabla u_n \|_{L^1(\Omega)} \to \int_{\Omega} |\nabla u|$ as $n \to + \infty$ (the existence of such a sequence $(u_n)$ is standard). Set, for $k \in \mN$, 
$$
\Omega_k = \big\{x \in K \mbox{ for some $k$-sim $K$ such that $K \cap \Omega \neq \O$} \big\}. 
$$
It is clear that, for each $n$,  
\begin{equation}\label{hahahaha}
\lim_{k \to + \infty} \int_{\Omega_k} |\nabla u_n| = \int_{\Omega} |\nabla u_n|. 
\end{equation}
By Lemma~\ref{lem-EF} (applied with $v = u_n$) and \eqref{hahahaha}, for each $n$,  there exist $k = k_n \in \mN$ and  $v_n \in R_k$ such that
\begin{equation}\label{hahahaha-1}
\|v_n - u_n \|_{W^{1, 1}(\mR^d)} \le  1/ n \quad \mbox{ and }
\int_{\Omega_k} |\nabla v_n| \le \int_{\Omega} |\nabla v_n| + 1/ n. 
\end{equation}
In what follows (except in the last two sentences), $n$ is fixed. 
By Case 1 (applied with $u = v_n$), there exists a family $(v_{\delta, n}) \subset L^1(\mR^d)$ such that $v_{\delta, n} \to v_n$ in $L^1(\mR^d)$ as $\delta \to 0$ and 
\begin{equation}\label{part7-haha}
\lims{\delta \to 0} \Lambda_{\delta}(v_{\delta, n}, \mR^d) \le \C \int_{\mR^d} |\nabla v_{n}|.  
\end{equation}
Applying  Lemma~\ref{lemlimsupaffine1}, we have 
\begin{equation*}
\limf{\delta \to 0} \Lambda_{\delta}(v_{\delta, n}, K) \ge \C \int_{K} |\nabla v_n| \mbox{ for each $k$-sim $K$}.  
\end{equation*}
Since 
$$
\mR^d \setminus \Omega_k = \mathop{\bigcup_{K \mbox{ is a $k$-sim}}}_{K \subset \mR^d \setminus \Omega_k} \bar K, 
$$
it follows that 
\begin{equation}\label{part8-haha}
\liminf_{\delta \to 0} \Lambda_{\delta}(v_{\delta, n}, \mR^d \setminus \Omega_k) \ge \C \int_{\mR^d \setminus \Omega_k} |\nabla v_n|. 
\end{equation}
Clearly
\begin{equation}\label{ll}
\Lambda_{\delta}(v_{\delta, n}, \mR^d \setminus \Omega_k)  + \Lambda_{\delta}(v_{\delta, n}, \Omega_k)  \le \Lambda_{\delta}(v_{\delta, n}, \mR^d).   
\end{equation}
We derive from  \eqref{part7-haha}, \eqref{part8-haha}, and \eqref{ll} that 
\begin{equation*}
 \lims{\delta \to 0} \Lambda_{\delta}(v_{\delta, n}, \Omega_k) \le \C \int_{\Omega_k} |\nabla v_n|,  
\end{equation*}
which implies, by \eqref{hahahaha-1},  
$$
\lims{\delta \to 0} \Lambda_{\delta}(v_{\delta, n}, \Omega) \le   \C \int_{\Omega} |\nabla v_n| + \C /n. 
$$
Hence there exists $\delta_n > 0$ such that, for $0 < \delta < \delta_n$ 
$$
\Lambda_{\delta}(v_{\delta, n}, \Omega) \le \C \int_{\Omega} |\nabla v_n| + \C/ n + 1/ n \quad \mbox{ and } \quad \|v_{\delta, n}  - v_n \|_{L^1(\Omega)} \le 1/ n. 
$$
Without loss of generality, one may assume that $(\delta_n)$ is decreasing to 0. Set
$$
u_{\delta} = v_{\delta_{n+1}, n} \mbox{ in } \Omega \quad  \mbox{ for } \delta_{n+1} \le \delta < \delta_n. 
$$
Then $(u_\delta)$ satisfies the required properties. \proofend 

\medskip

\subsection{Proof of Property (G1)}\label{sect-proof-G1}



Our goal is to show that (G1) holds with $K = \C$, i.e.,  
\begin{enumerate}
\item[(G1)] For every $u \in L^1(\Omega)$ and for every family
$(u_\delta)_{\delta \in (0,1)} \subset L^1(\Omega)$ such that
$u_\delta \to u$ in $L^1(\Omega)$ as $\delta \to 0$,
one has
\begin{equation*}
\liminf_{\delta \to 0} \Lambda_\delta(u_\delta) \ge \C \int_\Omega |\nabla u|. 
\end{equation*}
\end{enumerate}
The proof is in the spirit of \cite[Sections 5 and 6]{NgGamma} where Property (G1) is proved for $\Omega = \mR^d$ and $\varphi = c_1 \bvarphi_1$. 

\medskip 
Set 
\begin{equation*}
H(x) = \left\{\begin{array}{ll} 0 & \mbox{if } x_1 < 0, \\[6pt]
1 & \mbox{otherwise}, 
\end{array}\right.
\end{equation*}
and denote $H_c(x) : = H(x_1- c, x')$ for $(x_1, x') \in \mR \times \mR^{d-1}$ and  $c \in \mR$.  Define
\begin{equation}\label{defbN}
\bb := \inf \limf{\delta \lr 0}
\Lambda_{\delta} (g_\delta, Q), 
\end{equation}
where the infimum is taken over all families of functions
$(g_\delta)_{\delta \in (0,1)} \subset L^1(Q)$ such that $g_\delta \to H_\frac{1}{2}$ in $L^1(Q)$. Note that $\int_Q |\nabla H_{1/2}| = 1$. 
It follows from Property (G2) that 
\begin{equation}\label{gamma-K}
\gamma \le \C. 
\end{equation}

In the next section, we prove
\begin{proposition}\label{probC} We have
$$
\bb = \C.
$$
\end{proposition}

Clearly, this is consistent with Theorem~\ref{thm-gamma}. 

\subsubsection{Proof of Proposition~\ref{probC}}

The proof of Proposition~\ref{probC} is based on two lemmas. The first one in the spirit of Lemma~\ref{lemcorlemfondN-111} is: 

\begin{lemma}\label{lembC2-1} There exist a sequence $(h_k) \subset L^1(Q)$ and two sequences
$(\delta_k), (c_k) \subset \mR_+$  such that 
$$
\lim_{k \to + \infty} \delta_k = \lim_{k \to + \infty} c_k = 0, \quad \lim_{k \to + \infty} h_k = H_{1/2} \mbox{ in } L^1(Q), 
$$
$$
h_k(x) = 0 \mbox{ for  } x_1 < 1/2 - c_k, \quad h_k(x) = 1 \mbox{ for } x_1 > 1/2 + c_k, \quad 0 \le h_k(x) \le 1 \mbox{ in } Q, 
$$
and
\begin{equation*}
\lim_{k \lr \infty} \Lambda_{\delta_k} (h_k, Q)
= \bb.
\end{equation*}
\end{lemma}

\noindent{\bf Proof.} From the definition of $\bb$, there exist a
sequence $(\tau_k) \subset \mR_+$ and a sequence $(g_k) \subset L^1(Q)$ such that $\tau_k \to 0$, $g_k  \to H_\frac{1}{2}$  in $L^1(Q)$, and 
\begin{equation}\label{eqgk1}
\lim_{k \lr \infty} \Lambda_{\tau_k} (g_k, Q) = \bb.
\end{equation}
Set $c_k  = \|g_k - H_{1/2} \|_{L^1(Q)}^{1/4}$ so that 
\begin{equation}\label{condck1}
\lim_{k \to + \infty} c_k = 0 \quad  \mbox{ and } \quad \lim_{k \lr \infty} \frac{\big|\{ x \in Q; |g_k(x) - H_\frac{1}{2}(x)|
\ge c_k \}\big|}{c_k} \le \lim_{k \lr \infty} \frac{\|g_k - H_{1/2} \|_{L^1(Q)}}{c_k^2}   =0.
\end{equation}
Define two continuous functions $h_{1,k}, \, h_{2,k}: Q \to \mR$ which depend only on $x_1$ as follows
\begin{equation*}
\dsp h_{1,k}(x) = \left\{\begin{array}{cl} c_k & \mbox{if }
x_1 < \frac{1}{2} - c_k, \\[6pt]
1 + c_k & \mbox{if } x_1 > \frac{1}{2}, \\[6pt]
\mbox{ affine w.r.t. $x_1$} &
\mbox{otherwise},
\end{array} \right. \quad 
\dsp h_{2,k}(x) = \left\{\begin{array}{cl} - c_k & \mbox{if }
x_1 < \frac{1}{2}, \\[6pt]
1 - c_k & \mbox{if } x_1 > \frac{1}{2} + c_k, \\[6pt]
\mbox{ affine w.r.t. $x_1$} &
\mbox{otherwise}.
\end{array} \right.
\end{equation*}
Set
$$g_{1, k} = \min \big(\max(g_k, h_{2,k}), h_{1,k} \big) \quad \mbox{and} \quad g_{2, k} = \min \big( \max(g_{1,k}, c_k ), 1- c_k \big).
$$
It is clear that, in $Q$, 
\begin{equation}\label{TTT1}
g_{2,k}(x) = c_k \mbox{ for } x_1 < 1/2 - c_k,  \quad g_{2,k}(x) = 1 - c_k \mbox{ for } x_1 > 1/2
+ c_k, \quad c_k \le g_{2, k} \le 1 -c_k. 
\end{equation}
We claim that 
\begin{equation}\label{ineqg2kBV}
\lims{k \lr \infty} \Lambda_{\tau_k} (g_{2, k}, Q)  =
\bb.
\end{equation}
Indeed,  by 
Corollary~\ref{corCor00}, 
we have
\begin{equation}\label{toto1.1}
\Lambda_{\tau_k} (g_{2, k}, Q)  \le \Lambda_{\tau_k} (g_{1, k}, Q). 
\end{equation}
Note that 
$$
\|\nabla h_{1, k}\|_{L^\infty(Q)} \le 1/ c_k  \quad \mbox{ and } \quad \|\nabla h_{2, k}\|_{L^\infty(Q)} \le 1/ c_k. 
$$
Using \eqref{condck1} and applying  Lemma~\ref{lemtechnical1}, we obtain 
\begin{equation}\label{toto1}
\lims{k \lr \infty} \Lambda_{\tau_k} (g_{1, k}, Q) 
\le \lims{k \lr \infty}
\Lambda_{\tau_k} (g_{k}, Q) .
\end{equation}
Combining \eqref{eqgk1}, \eqref{toto1.1}, and \eqref{toto1} yields
\begin{equation*}
\lims{k \lr \infty} \Lambda_{\tau_k} (g_{2, k}, Q)
\le \bb; 
\end{equation*}
which is \eqref{ineqg2kBV}. 
One can now verify that the conclusion holds for $h_k : = (1- 2c_k)^{-1}(g_{2,k} - c_k)
$ and $ \delta_k := (1- 2 c_k)^{-1}\tau_k$   by  \eqref{TTT1} and \eqref{ineqg2kBV}.   \proofend

\medskip 
We next prove 

\begin{lemma}\label{lembC2} Set $g(x) = x_1$ in $Q$. There exist a sequence $(g_k) \subset L^1(Q)$ and a sequence
$(\delta_k) \subset \mR_+$  such that 
$$
\lim_{k \to + \infty} \delta_k = 0, \quad \lim_{k \to + \infty} g_k = g  \mbox{ in } L^1(Q), 
$$
and
\begin{equation*}
\lims{k \lr \infty} \Lambda_{\delta_k} (g_k, Q)
\le \bb.
\end{equation*}
\end{lemma}

\noindent{\bf Proof.} By Lemma~\ref{lembC2-1}, there exist a sequence $(h_k) \subset L^1(Q)$ and two sequences
$(\delta_k), (c_k) \subset \mR_+$  such that 
\begin{equation}\label{H-1-part1}
\lim_{k \to + \infty} \delta_k = \lim_{k \to + \infty} c_k = 0, 
\end{equation}
\begin{equation}\label{H-1-part2}
h_k(x) = 0 \mbox{ for  } x_1 < 1/2 - c_k, \quad h_k(x) = 1 \mbox{ for } x_1 > 1/2 + c_k, \quad 0 \le h_k(x) \le 1 \mbox{ in } Q, 
\end{equation}
and
\begin{equation}\label{H-1-part3}
\lim_{k \lr \infty} \Lambda_{\delta_k} (h_k, Q)
= \bb.
\end{equation}
{\bf Fix} $n \in \mN$ and consider the sequence
$(f_k): Q \mapsto \mR$  defined as follows
\begin{equation}\label{def-fk}
f_k (x) =  \frac{1}{n} h_k \Big(x_1 - \frac{j}{n} + \frac{1}{2} - \frac{1}{2n}, x' \Big) + \frac{j}{n} \mbox{ for } x \in Q_j, \, 0 \le j \le n-1,
\end{equation}
where $Q_j=  [j/n, (j+1)/ n] \times [0, 1]^{d-1}$.  We deduce from \eqref{H-1-part2} that 
\begin{equation}\label{limfk2-11}
\int_{Q} |f_k(x) - x_1| \, dx \le \frac{1}{n}.
\end{equation}
We claim that
\begin{equation}\label{limfk1--1}
\limsup_{k \lr \infty} \Lambda_{\delta_k/n} (f_k, Q) \le \limsup_{k \lr \infty} \Lambda_{\delta_k} (h_k, Q)  = \gamma.
\end{equation}
It is clear that 
\begin{equation}\label{TTT2}
\Lambda_{\delta_k/n} (f_k, Q)
\le   \sum_{j=0}^{n-1} \Lambda_{\delta_k/n} (f_k, Q_j) +  \sum_{j=0}^{n-1} \mathop{\iint}_{Q_j \times (Q \setminus Q_j)}
 \frac{ \varphi_{\delta_k/n} (|f_k(x) - f_k(y)|)}{|x-y|^{d+1}} \, dx \, dy.
\end{equation}
Set $\hat Q = [\frac{1}{2} - \frac{1}{2n}, \frac{1}{2} +
\frac{1}{2n}] \times [0, 1]^{d-1}$.  We have, by the definition of $f_k$, 
\begin{equation}\label{TTT3}
\Lambda_{\delta_k/n} (f_k, Q_j) 
= \frac{1}{n} \Lambda_{\delta_k} (h_k, \hat Q) \le \frac{1}{n}\Lambda_{\delta_k} (h_k, Q). 
\end{equation}
If $(x, y) \in Q_j \times (Q \setminus Q_j)$ then  $f_k(x) = f_k(y)$ if $|x_1 - y_1| < 1/ (2n) - c_k$ by \eqref{H-1-part2}. It follows from  \eqref{cond-varphi-01} and \eqref{H-1-part1} that
\begin{equation}\label{TTT4}
\limsup_{k \to + \infty }\mathop{\iint}_{Q_i \times (Q \setminus Q_i)}
 \frac{ \varphi_{\delta_k/n} (|f_k(x) - f_k(y)|)}{|x-y|^{d+1}} \, dx \, dy \le  \limsup_{k \to + \infty} \mathop{\mathop{\iint}_{Q_i \times (Q \setminus Q_i) }}_{|x_1 - y_1| \ge 1/ (2n) - c_k }  \frac{  b \delta_k / n }{ |x-y|^{d+1}} \, dx \, dy = 0 
\end{equation}
(recall that $n$ is fixed). 
Combining \eqref{H-1-part3}, \eqref{TTT2},  \eqref{TTT3}, and \eqref{TTT4} yields  \eqref{limfk1--1}. 

We now reintroduce the dependence on $n$. By the above, there exists $f_{k, n}$,  defined for $k, n \ge 1$,  such that 
\begin{equation*}
\int_{Q} |f_{k, n} (x) - x_1| \, dx \le \frac{1}{n}
\end{equation*}
and
\begin{equation*}
\limsup_{k \lr \infty} \Lambda_{\delta_k/n} (f_{k, n} , Q) \le  \gamma \mbox{ for each } n. 
\end{equation*}
Thus for each $n$, there exists  $k_n$ such that 
$\Lambda_{\delta_{k_n}/n} (f_{k_n, n} , Q) \le  \gamma + 1/n$. The desired conclusions hold for $(f_{k_n, n})$ and $(\delta_{k_n}/ n)$. 
\proofend

\medskip
%
%
%

\medskip 
\noindent{\bf Proof of Proposition~\ref{probC}}: We have $\C \le \bb$ by Lemmas \ref{lemlimsupaffine1} and \ref{lembC2}; and $\C \ge \bb$ by \eqref{gamma-K}. Hence $\bb= \C$.  \proofend

\subsubsection{Some useful lemmas}

We begin with a consequence of the definition of $\gamma$ and Proposition~\ref{probC}.  

\begin{lemma}\label{lem-Def} For any $\eps > 0$ there exists $\delta_\eps > 0$ such that if $\delta < \delta_\eps$ and $g \in L^1(Q)$ with $\|g - H_{1/2} \|_{L^1(Q)} < \delta_\eps$ then 
$$
\Lambda_\delta(g, Q) \ge \C - \eps. 
$$
\end{lemma}


We now prove 

\begin{lemma}\label{lemp=1.1}
Let $c, \tau > 0$ and $(g_\delta) \subset 
L^1(\R)$ with $\R = (a_1,b_1) \times (a, b)^{d-1} $ for some $a_1 < b_1$ and $a < b$ be such that $\tau < (b_1 - a_1)/ 8$. Assume that, for small $\delta$,  
$$
g_\delta(x) = 0 \mbox{ for } x_1 <  a_1
+ \tau, \quad  g_\delta(x) = c \mbox{ for } x_1 >  b_1 -
\tau, \quad \mbox{ and } \quad 0 \le g_\delta(x)  \le c \mbox{ for } x \in \R. 
$$  
We have, with $\R' = (a, b)^{d-1}$, 
\begin{equation*}
\limf{\delta \lr 0}  \Lambda_{\delta} (g_\delta, \R)
\ge
c \C |\R'|.
\end{equation*}
\end{lemma}

Here and in what follows, for a subset in $\mR^{d-1}$, $| \, \cdot \,  |$ denotes its  $(d-1)$-dimensional Hausdorff measure unless stated otherwise.

\medskip 

\noindent {\bf Proof.}  We only present the proof in two dimensions for simplicity of notations. Let $d = 2$. For $s > 0$, set 
$$
\R^s = (a_1, b_1) \times [(a, a+ s) \cup (b - s, b)]. 
$$
We first prove that, for every $s > 0$,   
\begin{equation}\label{M-1}
\liminf_{\delta \to 0} [ \Lambda_\delta(g_\delta, \R) + \Lambda_\delta(g_\delta, \R^s)] \ge c \C |\R'|. 
\end{equation}
Without loss of generality, one may assume that 
\begin{equation}\label{SS}
\R = (0, b_1) \times (0, b_2) \quad \mbox{ and } \quad c = 1. 
\end{equation}
Let $g_{1, \delta}: (0, b_1) \times \mR$ be such that 
\begin{equation}\label{def-g1}
g_{1, \delta}(x) = g_\delta(x) \mbox{ for } x \in \R, \quad g_{1, \delta}(x) = g_\delta(x_1, -x_2) \mbox{ for } x \in (0, b_1) \times (-b_2, 0),  
\end{equation}
and $g_{1, \delta}$ is a periodic function in $x_2$ with period  $2b_2$. Set 
$$
\R_j = (0, b_1) \times (j b_2, jb_2 + b_2) \mbox{ for } j \ge 0 \quad \mbox{ and } \quad \R(m)  =(0, b_1) \times (0, 2m  b_2) \mbox{ for } m \ge 0. 
$$
It is clear that, for $m \in \mN$,  
\begin{equation}\label{H-part0}
\Lambda_\delta (g_{1, \delta}, \R(m)) = \sum_{j = 0}^{2m -1} \Lambda_{\delta} (g_{1, \delta},  \R_j) + \sum_{j = 0}^{2m -1} \mathop{\iint}_{\R_j \times (\R(m) \setminus \R_j)}  \frac{ \varphi_{\delta} (|g_{1, \delta}(x) - g_{1, \delta}(y)|)}{|x-y|^{3}} \, dx \, dy.
\end{equation}
From the definition of $g_{1, \delta}$, we have, for $0 \le j \le 2 m -1$,  
\begin{equation}\label{H-part1}
\Lambda_{\delta} (g_{1, \delta},  \R_j)  = \Lambda_{\delta} (g_{\delta},  \R). 
\end{equation}
Clearly,   for $0 \le j \le 2 m -1$,   
\begin{align*}
 \mathop{\iint}_{\R_j \times (\R(m) \setminus \R_j)}  \frac{ \varphi_{\delta} (|g_{1, \delta}(x) - g_{1, \delta}(y)|)}{|x-y|^{3}} \, dx \, dy = &  \mathop{\mathop{\iint}_{\R_j \times (\R(m) \setminus \R_j)}}_{|x_2 - y_2| < s}  \frac{ \varphi_{\delta} (|g_{1, \delta}(x) - g_{1, \delta}(y)|)}{|x-y|^{3}} \, dx \, dy \\[6pt] 
 & + \mathop{\mathop{\iint}_{\R_j \times (\R(m) \setminus \R_j)}}_{|x_2 - y_2| \ge s}  \frac{ \varphi_{\delta} (|g_{1, \delta}(x) - g_{1, \delta}(y)|)}{|x-y|^{3}} \, dx \, dy; 
\end{align*}
which yields, by the definition of $g_{1, \delta}$ and \eqref{cond-varphi-01},  
\begin{equation}\label{H-part2}
 \mathop{\iint}_{\R_j \times (\R(m) \setminus \R_j)}  \frac{ \varphi_{\delta} (|g_{1, \delta}(x) - g_{1, \delta}(y)|)}{|x-y|^{3}} \, dx \, dy  \le  \Lambda_{\delta}(g_\delta, \R^s) + C \delta m/ s^3. 
\end{equation}
Here and in what follows in this proof, $C$ denotes a positive constant independent of $\delta$ and  $m$. Combining \eqref{H-part0}, \eqref{H-part1},  and \eqref{H-part2} yields 
\begin{equation}\label{N-1}
\Lambda_{\delta} (g_{1, \delta}, \R(m)) \le 2 m \Lambda_{\delta}(g_\delta, \R) + 2m \Lambda_{\delta}(g_\delta, \R^s) + C \delta m^2 / s^3. 
\end{equation}
Define $g_{2, \delta}: \mR^2 \to \mR$ as follows 
\begin{equation*}
g_{2, \delta}(x) = \left\{ \begin{array}{cl} g_{1, \delta}(x) & \mbox{ if } x_1 \in (0, b_1), \\[6pt]
0 &  \mbox{ if } x_1 \le 0, \\[6pt]
1 & \mbox{ if } x_1 \ge b_1. 
\end{array}\right.
\end{equation*}
Since $g_{1, \delta}(x) = 0 $ for $x_1 <  \tau$ and $g_{1, \delta} (x)= 1$ for $x_1 > b_1 - \tau$,  by  \eqref{cond-varphi-01}, we have,   for $m \in \mN$,  
\begin{equation}\label{N-2}
\Lambda_{\delta} (g_{2, \delta}, (-mb_2, mb_2) \times (0, 2m b_2)) \le \Lambda_{\delta} (g_{1, \delta}, \R(m))  + C \delta m^4/ \tau^3. 
\end{equation}
Set, for $m \in \mN$ and $x \in Q$,  
$$
 g_{3, \delta, m}(x) = g_{2, \delta} \big(2 b_2 m (x_1 - 1/ 2, x_2) \big). 
$$
By a change of variables,   we have
\begin{equation}\label{N-3}
\Lambda_{\delta}(g_{2, \delta},  (-m b_2, mb_2) \times (0, 2 m b_2)) = 2 m b_2 \Lambda_{\delta}(g_{3, \delta, m}, Q). 
\end{equation}
Combining \eqref{N-1}, \eqref{N-2}, and \eqref{N-3} yields
\begin{equation}\label{N-3-1}
b_2 \Lambda_{\delta}(g_{3, \delta, m}, Q) \le \Lambda_{\delta}(g_\delta, \R) + \Lambda_{\delta}(g_\delta, \R^s) + C \delta m / s^3 + C \delta m^3 / \tau^3. 
\end{equation}
Since $0 \le g_{2, \delta} \le 1$,  it follows  from the definition of $g_{3, \delta, m}$ that  
\begin{equation*}
\|g_{3, \delta, m} - H_{1/2} \|_{L^1(Q)} \le C/m. 
\end{equation*}
By Lemma~\ref{lem-Def}, for every $\eps > 0$ there exists $m_\eps > 0$ such that if $m \ge m_\eps$ then 
\begin{equation}\label{N-4}
\liminf_{\delta \to 0} \Lambda_{\delta}(g_{3, \delta, m}, Q) \ge \C - \eps. 
\end{equation}
Taking $m = m_\eps$ in \eqref{N-3-1}, we derive from \eqref{N-4} that 
\begin{equation}
\liminf_{\delta \to 0} [ \Lambda_\delta(g_\delta, \R) + \Lambda_\delta(g_\delta, \R^s)] \ge (\C -\eps) b_2. 
\end{equation}
Since $\eps > 0$ is arbitrary, we obtain \eqref{M-1} by \eqref{SS}. 

We are now ready to prove 
\begin{equation}\label{M-3}
\liminf_{\delta \to 0}  \Lambda_\delta(g_\delta, \R)  \ge c \C |\R'|. 
\end{equation}
Without loss of generality, one may again assume \eqref{SS}. 
Fix $n \in \mN$ (arbitrary) and define $s = s(n)= b_2/ (4n^2)$. For $0 \le j \le n$,   by \eqref{M-1} (applied with $\R = \R \setminus \R^{j s}$), we have, 
$$
\liminf_{\delta \to 0} \big[\Lambda_{\delta}(g_\delta, \R \setminus \R^{j s}) + \Lambda_{\delta}(g_\delta, \R^{j s + s} \setminus \R^{j s}) \big] \ge \C (b_2 - 2 j s). 
$$
Summing these inequalities for $j$ from $0$ to $n-1$ and noting that 
$$
\Lambda_{\delta}(g_\delta, \R) \ge \sum_{j=0}^{n-1} \Lambda_{\delta}(g_\delta, \R^{j s + s} \setminus \R^{j s}) \quad 
\mbox{ and } \quad b_2 - 2 j s  \ge b_2 - b_2/ (2n), 
$$
we obtain 
$$
\liminf_{\delta \to 0 }(n+1)\Lambda_{\delta}(g_\delta, \R) \ge n \C b_2 [1 - 1/ (2n)]. 
$$
This implies 
$$
\liminf_{\delta \to 0} \Lambda_{\delta}(g_\delta, \R) \ge \frac{n}{n+1} \C b_2 [1 - 1/ (2n)]. 
$$
Since $n \in \mN$ is arbitrary, we obtain 
$$
\liminf_{\delta \to 0}  \Lambda_{\delta}(g_\delta, \R) \ge \C b_2. 
$$
The proof is complete. \proofend

\medskip

Here is a more general version of Lemma~\ref{lemp=1.1}. 
\begin{lemma}\label{lemp=1.1.1}
Let $c, \tau > 0$ and $(g_\delta) \subset 
L^1(\R)$ with $\R = (a_1,b_1) \times (a, b)^{d-1} $ for some $a_1 < b_1$ and $a < b$  be such that $\tau < (b_1 - a_1)/ 8$. Set
$$
A_\delta = \big\{x \in \R:  g_\delta(x) > 0 \mbox{ and } a_1 \le  x_1 \le a_1 + \tau\big\} \mbox{ and } B_\delta  = \big\{x \in \R: \;  g_\delta(x) <  c \mbox{ and }  b_1 -\tau \le x_1 \le  b_1 \big\}. 
$$
We have, with $\R' = (a, b)^{d-1}$, 
\begin{equation*}
\limf{\delta \lr 0}  \Lambda_{\delta} (g_\delta, \R)
\ge
c \C |\R'| - C_d c \limsup_{\delta \to 0} (|A_\delta| + |B_\delta|)/ \tau,  
\end{equation*}
\end{lemma}

\noindent{\bf Proof.} Define two continuous functions $f_1$ and $f_2$ in $\R$ which depend only on $x_1$ as follows 
\begin{equation*}
f_{1}(x) =   \left\{ \begin{array}{cl} \dsp  0 & \mbox{if } x_1 \le a_1 +  \tau/2,  \\[6pt]
c & \mbox{if } x_1 \ge a_1 +  \tau, \\[6pt]
\dsp \mbox{ affine w.r.t. $x_1$}  & \mbox{otherwise}, 
\end{array} \right.
\end{equation*}
and 
\begin{equation*}
f_{2}(x) =  \left\{ \begin{array}{cl} \dsp  0 & \mbox{if } x_1 \le b_1 -  \tau,  \\[6pt]
c & \mbox{if } x_1 \ge b_1- \tau/2, \\[6pt]
\dsp \mbox{ affine w.r.t. $x_1$}  & \mbox{otherwise}.
\end{array} \right.
\end{equation*}
Define 
$$
h_{1, \delta} = \max \Big( \min \big(g_\delta, c \big), 0 \Big),  \quad
h_{2, \delta} = \min( h_{1, \delta}, f_1), \quad \mbox{ and }
\quad 
h_{3, \delta} = \max(h_{2,\delta}, f_{2} ).
$$
By Corollary~\ref{corCor00}, we have 
\begin{equation*}
\Lambda_{\delta} (h_{1, \delta}, \R) \le \Lambda_{\delta} (g_\delta,  \R), 
\end{equation*}
and by Lemma~\ref{lemtechnical1}, we obtain 
\begin{equation*}
\Lambda_{\delta} (h_{2, \delta},  \R) \le \Lambda_{\delta} (h_{1, \delta},  \R)  + C_d c |\{h_{1, \delta} > f_1\}|/ \tau  \le  \Lambda_{\delta} (h_{2, \delta},  \R)  + C_dc  |A_\delta|/ \tau, 
\end{equation*} 
and 
\begin{equation*}
\Lambda_{\delta} (h_{3, \delta},  \R) \le \Lambda_{\delta} (h_{2, \delta},  \R)  + C_d c |\{h_{2, \delta} <  f_2\}| / \tau  \le \Lambda_{\delta} (h_{2, \delta},  \R)   + C_d c |B_\delta|/ \tau. 
\end{equation*} 
It follows that 
\begin{equation*}
\Lambda_{\delta} (h_{3, \delta},  \R)  \le \Lambda_{\delta} (g_\delta,  \R) + C_d c (|A_\delta| + |B_\delta|)/ \tau. 
\end{equation*}
One can easily check that  $h_{3, \delta}(x) = 0$ for $x_1 \le a_1 + \tau/ 2$,  $h_{3, \delta} (x)= c $ for $x_1 \ge b_1 -\tau/2$, and $0 \le h_{3, \delta} \le c$ in $\R$. 
Applying  Lemma~\ref{lemp=1.1} for $h_{3, \delta}$, we obtain the conclusion. \proofend

\medskip 
We next recall the definition of a Lebesgue surface (see \cite{NgGamma}): 

\begin{definition}\fontfamily{m} \selectfont
Let $g \in L^1 (\R)$ with $\R = \prod_{i=1}^d (a_i, b_i)$ for some $a_i < b_i$ ($1 \le i \le  d$) and $t \in
(a_1, b_1)$. Set $\R' = \prod_{i=2}^d (a_i, b_i)$. 
The surface $x_1 = t$ is said to be a Lebesgue
surface of $g$ if  for almost every $z' \in \R'$, $(t, z')$ is a Lebesgue point of $g$,
the restriction of $g$ to the surface $x_1 = t$ is integrable with
respect to $(d-1)$-Hausdorff measure, and
\begin{equation}\label{defBV}
\lim_{\eps \lr 0_+}\mathop{\mint}_{t - \eps }^{t + \eps}
\mathop{\int}_{\R'}| g(s, z') - g(t, z') | \,
dz' \, ds = 0.
\end{equation}
For $i=2, \dots, d$, we also define the notion of the Lebesgue
surface for surfaces $x_i = t$ with $t \in (a_i, b_i)$ in a 
similar manner.
\end{definition}

%
The following lemma plays a crucial role in our analysis; its proof relies on Lemma~\ref{lemp=1.1.1}.

\begin{lemma}\label{lemp=1liminf}
Let $g \in L^1 (\R)$ and $(g_\delta) \subset L^1(\R)$ with $\R =  \prod_{i=1}^d (a_i, b_i)$ for some $a_i < b_i$ ($1 \le i \le  d$)
such that $(g_\delta) \to g$ in  $L^1(\R)$. Set $\R' = \prod_{i=2}^d (a_i, b_i)$. 
Let  $a_1 < t_1 < t_2  <  b_1$  be such
that the surface $x_1 = t_j$ $(j=1, \, 2)$ is a Lebesgue surface
of $g$. We have
\begin{equation*} \liminf_{\delta \to 0} \Lambda_\delta \big(g_\delta, (t_1, t_2) \times
 \R' \big)
  \ge \C
\mathop{\int}_{\R' }|g(t_2, x') - g(t_1,x')|
\, d x'.
\end{equation*}
\end{lemma}

\noindent{\bf Proof.} Fix $\eps > 0$ (arbitrary). Let $A'$ be the
set of all elements $z' \in \R'$ such that, for $j=1, 2$, $(t_j, z')$ is a Lebesgue point of $g(t_j, \cdot)$ and $(t_j, z')$ is a Lebesgue point of $g$. Then $|A'| = |\R'|$. 
For each $z' \in A'$, there exists $l(z', \eps) > 0$ such that for every closed cube $Q'_l(z') \subset \mR^{d-1}$ centered at $z'$ with length $0 < l < l(z', \eps)$, we have
\begin{equation}\label{eq4.1p=1}
\left\{\begin{array}{c} \Big| \big\{(x_1, y') \in (t_1, t_1 + l) \times
Q'_l(z'); \, |g(x_1, y') - g(t_1, z')| \ge \eps/2 \big\} \Big| \le \eps l^d, \\[6pt]
\Big| \big\{(x_1, y') \in (t_2 - l, t_2) \times Q'_l(z'); \, |g(x_1, y') -
g(t_2, z')| \ge \eps/2 \big\} \Big|  \le \eps l^d,
\end{array} \right.
\end{equation}
and, for $j=1, 2$,  
\begin{equation}\label{totohaha}
\mint_{Q'_l(z')} |g(t_j, y') - g(t_j, z')| \, dy' < \eps. 
\end{equation}
Fix $z' \in A'$ and $0 < l < l(z', \eps)$. 
Since $(g_\delta) \to g$ in $L^1(\R)$, it follows
from~\eqref{eq4.1p=1} that, when $\delta$ is small,
\begin{equation} \label{equ2}
\left\{\begin{array}{c}  \Big| \big\{(x_1, y') \in (t_1, t_1 + l)
\times Q'_l(z'); \, |g_\delta(x_1, y') - g(t_1, z')| \ge \eps \big\} \Big|  \le 2 \eps l^d,\\[6pt]
 \Big| \big\{(x_1, y') \in (t_2 -  l, t_2) \times Q'_l(z'); \,
|g_\delta(x_1, y') - g(t_2, z')| \ge \eps \big\} \Big|
\le 2 \eps l^d.
\end{array} \right.
\end{equation}
We claim that 
\begin{equation}\label{jump-1}
\limf{\delta \lr 0} \Lambda_{\delta} \big(g_\delta, (t_1, t_2) \times
Q'_l(z') \big) \ge \C \big|g(t_2, z') - g(t_1,
z')\big|\, |Q'_l(z')| - C \eps|Q'_l(z')|, 
\end{equation}
for some positive constant $C$ depending only on $d$. Indeed, without loss of generality, one may assume that $g(t_2, z') \ge  g(t_1, z')$. It is clear that \eqref{jump-1} holds if $g(t_2, z') - g(t_1, z') \le 4 \eps$ by choosing $C = 10$. We now consider the case $g(t_2, z') - g(t_1, z') > 4 \eps$.  Applying Lemma~\ref{lemp=1.1.1} for $g_\delta - g(t_1, z') - \eps$ in the set $(t_1, t_2) \times
Q'_l(z')$,    $\tau = l$, and $c = g(t_2, z')  - g(t_1, z') - 2 \eps$,  we derive that
\begin{equation*}
\limf{\delta \lr 0} \Lambda_{\delta} \big(g_\delta, (t_1, t_2) \times
Q'_l(z') \big) \ge \C \big|g(t_2, z') - g(t_1,
z') -  2 \eps \big|\, |Q'_l(z')| - C \eps|Q'_l(z')| 
\end{equation*}
since,  by  \eqref{equ2},  
\begin{multline*}
|A_\delta|  =   \Big| \big\{(x_1, y') \in (t_1, t_1 + l) \times
Q'_l(z'); \,g_\delta (x_1, y') - g(t_1, z') - \eps > 0 \big\} \Big| \\[6pt]
 \le    \Big| \big\{(x_1, y') \in (t_1, t_1 + l) \times
Q'_l(z'); \, |g_\delta(x_1, y') - g(t_1, z')| \ge \eps \big\} \Big|  \le 2 \eps l^d = 2   \eps l |Q'_l(z')|
\end{multline*}   
and 
\begin{multline*}
|B_\delta|  =   \Big| \big\{(x_1, y') \in (t_2 - l, t_2) \times
Q'_l(z'); \,g_\delta (x_1, y') - g(t_1, z') - 2 \eps <  g(t_2, z')  - g(t_1, z') - \eps \big\} \Big| \\[6pt]
 \le    \Big| \big\{(x_1, y') \in (t_1, t_1 + l) \times
Q'_l(z'); \, |g(x_1, y') - g(t_1, z')| \ge \eps \big\} \Big|  \le 2 \eps l^d = 2  \eps l |Q'_l(z')|. 
\end{multline*}  
This implies Claim~\eqref{jump-1}. 

On the other hand, by Besicovitch's covering theorem (see e.g.,  \cite[Corollary 1 on page 35]{EGmeasure} \footnote{In \cite{EGmeasure}, this result is stated for balls but a similar argument works for cubes with arbitrary orientations. }), there exist a sequence   $(z_k')_{k \in \mN} \subset  A'$ and  disjoint cubes $\big(Q'_{l_k}(z_k') \big)_{k \in \mN} \subset \R'$ such that $0 < l_k < l(z_k', \eps)$ for every $k$, and 
\begin{equation}\label{coucou0}
|A'| = \sum_{ k \in \mN}  |Q'_{l_k}(z_k')|. 
\end{equation}
[For the convenience of the reader, we explain how to apply  \cite[Corollary 1]{EGmeasure} in our situation. We take $n = d-1$, $A$ is our $A'$, $
{\cal F} = \big\{ Q'_l(z'); \; z' \in A' \mbox{ and } 0 < l < l(z', \eps) \big\}$, and $U = \R'$. ]

It follows from \eqref{jump-1} that
\begin{equation}\label{coucou1}
\limf{\delta \lr 0}  \Lambda_\delta \big(g_\delta, (t_1, t_2) \times \R' \big)
 \ge \sum_k \Big( \C \big|g(t_2, z_k') - g(t_1,
z_k')\big|\, |Q'_{l_k}(z_k')| - C \eps|Q'_{l_k}(z_k')| \Big). 
\end{equation}
We claim that 
\begin{equation}\label{claim-claim}
 \big|g(t_2, z_k') - g(t_1,
z_k')\big|\, |Q'_{l_k}(z_k')|  \ge \int_{Q'_{l_k}(z_k')} |g(t_2, y') - g(t_1,
y')|  dy' - 2 \eps |Q'_{l_k}(z_k')|. 
\end{equation}
Indeed, we have
\begin{align*}
\int_{Q'_{l_k}(z_k')} |g(t_2, y') - g(t_1,
y')|  dy' 
\le &  \int_{Q'_{l_k}(z_k')} |g(t_2, y')  - g(t_2, z')| +  |g(t_1,
y') - g(t_1, z')|  dy' \\[6pt] 
& + \big|g(t_2, z_k') - g(t_1,
z_k')\big|\, |Q'_{l_k}(z_k')|; 
\end{align*}
which implies \eqref{claim-claim} by \eqref{totohaha}. 

Combining \eqref{coucou0}, \eqref{coucou1}, and \eqref{claim-claim} yields 
$$
\limf{\delta \lr 0}  \Lambda_\delta \big(g_\delta, (t_1, t_2) \times \R' \big) \ge \C \int_{\R'}|g(t_2, y') - g(t_1,
y')|  dy'  - C \eps |\R'|. 
$$
Since $\eps > 0$ is arbitrary, we obtain the conclusion. \proofend

\medskip

We next recall the notion of 
 essential variation in \cite{NgGamma} related to $BV$ functions. 

\begin{definition} \fontfamily{m} \selectfont
Let $g \in L^1 (\R)$ with $\R = \prod_{i=1}^d (a_i, b_i)$ for some $a_i < b_i$ ($1 \le i \le  d$). Set $\R' = \prod_{i=2}^d (a_i, b_i)$.  The essential variation
of $g$ in the first direction is defined as follows
\begin{equation*}
\ess V(g,1, \R) = \sup \left\{ \sum_{i =1}^{m}
\mathop{\int}_{\R'}|g(t_{i+1}, x') -
g(t_i,x')| \, d x' \right\},
\end{equation*}
where the supremum is taken over all finite partitions $\{a_1 <
t_1 < \dots < t_{m+1} < b_1 \}$ such that the surface $x_1 = t_k$
is a Lebesgue surface of $g$ for $1 \le k \le m+1$. For $2 \le j
\le d$, we also define $\ess V(g, j)$ the essential variation of
$g$ in the $j^{\mathrm{th}}$ direction  in a similar manner.
\end{definition}

The following result  provides a characterization of $BV$ functions (see e.g., \cite[Proposition 3]{NgGamma}).

\begin{proposition}\label{proBV} Let $g \in L^1 (\R)$ with $\R = \prod_{i=1}^d (a_i, b_i)$ for some $a_i < b_i$ ($1 \le i \le  d$). 
Then $g \in BV(\R)$ if and only if
\begin{equation*}
\ess V(g, j, \R) < + \infty, \quad \forall \, 1 \le j \le d.
\end{equation*}
Moreover, for $g \in BV(\R)$, 
\begin{equation*}
\ess V(g, j, \R) = \int_{\R} \left| \frac{\partial g}{\partial x_j} \right|\quad
\forall \, 1 \le j \le d.
\end{equation*}
\end{proposition}

As a consequence of Lemma~\ref{lemp=1liminf} and Proposition~\ref{proBV}, we have 
\begin{corollary}\label{cor-Gamma}
Let $g \in L^1 (\R)$ and $(g_\delta) \subset L^1(\R)$ with $\R =  \prod_{i=1}^d (a_i, b_i)$ for some $a_i < b_i$ ($1 \le i \le  d$)
such that $(g_\delta) \to g$ in  $L^1(\R)$.  We have
\begin{equation*} \liminf_{\delta \to 0} \Lambda_\delta \big(g_\delta,
 \R \big)
  \ge \C   \int_{\R} \left| \frac{\partial g}{\partial x_j} \right| \quad
\forall \, 1 \le j \le d.
\end{equation*}
\end{corollary}

\subsubsection{Proof of Property (G1) completed} 
Recall that for each $u \in BV(\Omega)$, $|\nabla u |$ is a Radon
measure on $\Omega$. By Radon-Nikodym's theorem, we may write 
$$
\nabla u = \sigma |\nabla u|, 
$$
for some $\sigma \in L^\infty(\Omega, |\nabla u|, \mR^d)$
 and $|\sigma | = 1$ $|\nabla u|$- a.e. 
(see e.g.,  \cite[Theorem 1 on page 167]{EGmeasure}). Then for $|\nabla u |$- a.e.  $x \in \Omega$, one has (see e.g., 
\cite[Theorem 1 on page 43]{EGmeasure})
\begin{equation*} 
 \lim_{r \lr 0}  \frac{1}{|\nabla u| (Q_r(x, \sigma(x)))} \int_{Q_r(x, \sigma(x)
)} \sigma(y)  |\nabla u(y)| \, dy =   \sigma(x) \; \footnote{In \cite{EGmeasure}, this result is stated for balls but a similar argument works for cubes with arbitrary orientations. }.  
\end{equation*}
Hereafter for any $(x, \sigma, r) \in \Omega \times \mS^{d-1}
\times (0, + \infty)$, $Q_r(x, \sigma)$ denotes the closed cube
centered at $x$ with edge length $2r$ such that
one of its faces is orthogonal to $\sigma$. It follows that, for $|\nabla u |$- a.e.  $x \in \Omega$, 
\begin{equation*}
\lim_{r \lr 0} \dsp \frac{1}{|\nabla u| (Q_r(x, \sigma(x)))}\int_{Q_r(x, \sigma(x)
)} \sigma (y) \cdot \sigma (x) |\nabla u(y)| \, dy  = 1. 
\end{equation*}
Since 
$$
\sigma (y) \cdot \sigma(x) \le |\sigma (y) \cdot \sigma(x)| \le 1, 
$$
we derive that, for $|\nabla u |$- a.e.  $x \in \Omega$, 
\begin{equation*}
\lim_{r \lr 0} \dsp \frac{1}{|\nabla u| (Q_r(x, \sigma(x)))}\int_{Q_r(x, \sigma(x)
)} |\sigma (y) \cdot \sigma (x)|  \, |\nabla u(y)| \, dy  = 1. 
\end{equation*}
In other words, for $|\nabla u |$- a.e.  $x \in \Omega$, 
\begin{equation}\label{condDgae}
\lim_{r \lr 0} \dsp \int_{Q_r(x, \sigma(x)
)} |\nabla u(y) \cdot \sigma (x)| \, dy  \Big/ \dsp \int_{Q_r(x, \sigma(x)
)} |\nabla u(y)| \, dy = 1. 
\end{equation}

Denote $A  = \{x \in \Omega;  \eqref{condDgae} \mbox{ holds} \}$. Fix $\eps > 0$ (arbitrary). For $x \in A$, there exists a sequence $ s_n = s_n(x, \eps) \to 0$ as $n \to + \infty$  such that, for all $n$,  
\begin{equation}\label{condDgae-11}
 \int_{Q_{s_n}(x, \sigma(x)
)} |\nabla u(y) \cdot \sigma (x)| \, dy  \Big/ \dsp \int_{Q_{s_n}(x, \sigma(x)
)} |\nabla u(y)| \, dy \ge 1 - \eps. 
\end{equation}
and
\begin{equation}\label{condDgae-1111}
\int_{\partial Q_{s_n} (x, \sigma(x))} |\nabla u (y)| \, dy = 0. 
\end{equation}
Set 
$$
{\cal F} = \big\{ Q_{s_n(x, \eps)}(x, \sigma(x)); x \in A \mbox{ and } n \in \mN\big\}. 
$$
By Besicovitch's covering theorem (see e.g.,  \cite[Corollary 1 on page 35]{EGmeasure} applied with $A$, ${\cal F}$, $U = \Omega$, and $\mu = |\nabla u|$), there exists a  collection  of
disjoint cubes  $\Big(Q_{r_k} \big(x_k, \sigma(x_k) \big) \Big)_{k \in \mN}$ with $x_k \in A$ and $r_k = s_{n_k}(x_k, \eps)$ such that
\begin{equation}\label{toT1}
| \nabla u | (\Omega) = |\nabla u | \Big(\bigcup_{k \in \mN} Q_{r_k}\big(x_k, \sigma(x_k) \big)\Big).
\end{equation}
From \eqref{condDgae-11} and \eqref{condDgae-1111}, we have 
\begin{equation}\label{condDgae-11-2}
\int_{Q_{r_k}(x_k, \sigma(x_k)
)} |\nabla u(y)| \, dy \le \frac{1 }{1 - \eps}  \int_{Q_{r_k}(x_k, \sigma(x_k)
)} |\nabla u(y) \cdot \sigma (x_k)| \, dy
\end{equation}
and
\begin{equation}\label{condDgae-1111-2}
\int_{\partial Q_{r_k} (x_k, \sigma(x_k))} |\nabla u (y)| \, dy = 0. 
\end{equation}
Combining  \eqref{toT1} and \eqref{condDgae-11-2} yields
\begin{equation}\label{toT2}
 | \nabla u |(\Omega) \le \frac{1}{1-\eps} \sum_{k \in \mN}  \int_{Q_{r_k}(x_k, \sigma(x_k)
)} |\nabla u(y) \cdot \sigma (x_k)| \, dy . 
\end{equation}
Applying Corollary~\ref{cor-Gamma} and using \eqref{condDgae-1111-2}, we obtain 
\begin{equation}\label{toT3}
\C \int_{Q_{r_k}(x_k, \sigma(x_k)
)} |\nabla u(y) \cdot \sigma (x_k)| \, dy \le
\limf{\delta \lr 0} \Lambda_{\delta} \big(u_\delta, Q_{r_k}(x_k, \sigma(x_k)) \big).
\end{equation}
From~\eqref{toT2} and \eqref{toT3}, we have
\begin{equation}\label{toT4}
\C | \nabla u |(\Omega) \le \frac{1}{ 1- \eps} \limf{\delta \lr 0} \Lambda_{\delta} \big(u_\delta, \Omega \big). 
\end{equation}
Since $\eps > 0$ is arbitrary, we have established that, for $u \in BV(\Omega)$, 
\begin{equation*}
\limf{\delta \lr 0} \Lambda_{\delta} \big(u_\delta, \Omega \big) \ge  \C | \nabla u |(\Omega) . 
\end{equation*}
Suppose now that $u \in BV_{\loc}(\Omega)$, we may apply the above for any $\omega \subset \subset \Omega$ and  therefore we conclude that 
$$
\limf{\delta \lr 0} \Lambda_{\delta} \big(u_\delta, \Omega \big) \ge  \C | \nabla u |(\Omega) . 
$$
Hence it now suffices to prove that if $\liminf_{\delta \to } \Lambda_{\delta}(u_\delta, \Omega) < + \infty$, then $u \in BV_{\loc}(\Omega)$. Indeed, this is a consequence of Corollary~\ref{cor-Gamma}. 

\medskip 
The proof is complete. 
\proofend

\subsection{Further properties of $K(\varphi)$} \label{sect-k}

This section deals with  properties of  $\C (\varphi)$ defined in \eqref{k-Gamma}. Our main result is:

\begin{theorem}\label{thm-K} We have
\begin{equation}\label{main-K}
\C (\varphi) \ge \C (c_1 \bvarphi_1) \quad \mbox{ for all } \varphi \in {\cal A};
\end{equation}
in particular,
$$
\inf_{\varphi \in {\cal A}} \C (\varphi) > 0.
$$
\end{theorem}

\noindent{\bf Proof.} The proof uses an idea in \cite[Section 2.3]{NgSob2}.
From the definition of $\C(c_1 \bvarphi_1)$, we have (see \cite[Lemma 8]{NgGamma})
\begin{equation}\label{main-K-1}
\forall \, \eps > 0, \; \exists \delta (\eps) > 0 \; \mbox{such that if } \| v - U \|_{L^1(Q)} < \delta(\eps) \mbox{ and } \delta < \delta(\eps),  \mbox{ then } \Lambda_\delta(v, c_1 \bvarphi_1) \ge \C(c_1 \bvarphi_1) - \eps.
\end{equation}
Next we {\bf fix} $\varphi \in {\cal A}$. 
{\bf Fix}  $(u_\delta) \subset L^1(Q)$ be such that $u_\delta \to U$ in $L^1(Q)$.
Our goal is to prove that
\begin{equation}\label{claim-k}
\liminf_{\delta \to 0 }\Lambda_\delta(u_\delta, \varphi) \ge \C(c_1 \bvarphi_1).
\end{equation}
Let $c > 1$ and $\eps > 0$. Since $\varphi$ is non-decreasing,
we have
\begin{equation}\label{p1-K}
 \int_{Q} \int_{Q} \frac{\varphi_\delta(|u_\delta(x) - u_\delta(y)|)}{|x - y|^{d+1}} \, dx dy
\ge  \sum_{k = - \infty}^{\infty} \mathop{\int_{Q} \int_{Q}}_{c^{-k-1} <  |u_\delta(x) - u_\delta(y)| \le c^{-k}} \frac{\varphi_\delta(c^{-k-1})}{|x - y|^{d+1}} \, dx dy.
\end{equation}
Using the fact that
\begin{multline*}
\mathop{\int_{Q} \int_{Q}}_{c^{-k-1} < |u_\delta(x) - u_\delta(y)| \le c^{-k}} \frac{1}{|x - y|^{d+1}} \, dx dy  \\[6pt]
= \mathop{\int_{Q} \int_{Q}}_{|u_\delta(x) - u_\delta(y)| >  c^{-k - 1}} \frac{1}{|x - y|^{d+1}} \, dx dy  - \mathop{\int_{Q} \int_{Q}}_{|u_\delta(x) - u_\delta(y)| >  c^{-k }} \frac{1}{|x - y|^{d+1}} \, dx dy,
\end{multline*}
we obtain
\begin{multline}\label{p2-K}
\sum_{k = - \infty}^{\infty}  \mathop{\int_{Q} \int_{Q}}_{c^{-k-1} <  |u_\delta(x) - u_\delta(y)| \le c^{-k}} \frac{\varphi_\delta(c^{-k-1})}{|x - y|^{d+1}} \, dx dy  \\[6pt]
=   \sum_{k = - \infty}^{\infty} \big[ \varphi_\delta (c^{-k}) - \varphi_\delta (c^{-k-1}) \big]  \mathop{\int_{Q} \int_{Q}}_{|u_\delta(x) - u_\delta(y)| >  c^{-k }} \frac{1}{|x - y|^{d+1}} \, dx dy.
\end{multline}
We have, for any $k_0 > 0$,
\begin{multline}\label{toto1-1}
  \sum_{k = - \infty}^{\infty} \big[ \varphi_\delta (c^{-k}) - \varphi_\delta (c^{-k-1}) \big]  \mathop{\int_{Q} \int_{Q}}_{|u_\delta(x) - u_\delta(y)| >  c^{-k }} \frac{1}{|x - y|^{d+1}} \, dx dy \\[6pt]
\ge \sum_{k = k_0}^{\infty} \big[ \varphi_\delta (c^{-k}) - \varphi_\delta (c^{-k-1}) \big]  c^k \mathop{\int_{Q} \int_{Q}}_{|u_\delta(x) - u_\delta(y)| >  c^{-k }} \frac{c^{-k}}{|x - y|^{d+1}} \, dx dy.
\end{multline}
Applying \eqref{main-K-1} with $v  = u_\delta$ and $\delta = c^{-k}$, we obtain
\begin{equation}\label{k-delta-eps}
c_1 \mathop{\int_{Q} \int_{Q}}_{|u_\delta(x) - u_\delta(y)| >  c^{-k }} \frac{c^{-k}}{|x - y|^{d+1}} \, dx dy \ge \C(c_1 \bvarphi_1) - \eps,
\end{equation}
provided $\|u_\delta - U \|_{L^1(Q)} < \delta(\eps)$ and $c^{-k} < \delta(\eps)$.
In particular, there exist $\tilde \delta(\eps) > 0$ and $k(\eps, c) \in \mN$ such that \eqref{k-delta-eps} holds for $\delta < \tilde \delta(\eps)$ and $k \ge k(\eps, c)$.
Combining \eqref{toto1-1} and \eqref{k-delta-eps} yields
\begin{multline}\label{p3-K}
  \sum_{k = - \infty}^{\infty} \big[ \varphi_\delta (c^{-k}) - \varphi_\delta (c^{-k-1}) \big]  \mathop{\int_{Q} \int_{Q}}_{|u_\delta(x) - u_\delta(y)| >  c^{-k }} \frac{1}{|x - y|^{d+1}} \, dx dy \\[6pt]
\ge  c_1^{-1}  \sum_{k_0}^{\infty} \big[\C(c_1 \bvarphi_1) - \eps \big] c^k \big[ \varphi_\delta (c^{-k}) - \varphi_\delta (c^{-k-1}) \big]  ,
\end{multline}
for $k_0 = k(\eps, c) $ and  $\delta < \tilde \delta(\eps)$.
We derive from \eqref{p1-K}, \eqref{p2-K}, and \eqref{p3-K} that, for $\delta < \tilde \delta(\eps)$,
\begin{equation}\label{p4-K}
\Lambda_\delta(u_\delta) \ge  c_1^{-1} \big[\C(c_1 \bvarphi_1) - \eps \big]   \sum_{k = k_0}^{\infty} c^k \big[ \varphi_\delta (c^{-k}) - \varphi_\delta (c^{-k-1}) \big].
\end{equation}
We have, since $\varphi_\delta \ge 0$,
\begin{equation}\label{AA-2}
\sum_{k=k_0}^{\infty} c^k  \big[ \varphi_\delta (c^{-k}) - \varphi_\delta (c^{-k-1}) \big] =    \sum_{k= k_0}^{\infty} c^k \varphi_\delta (c^{-k}) -  \sum_{k= k_0}^{\infty}  c^k\varphi_\delta (c^{-k-1})  \ge \frac{1}{c} \sum_{k=k_0+1}^{\infty} \varphi_\delta (c^{-k}) c^k (c - 1)
\end{equation}
and, since $\varphi_\delta$ is non-decreasing,
\begin{multline}\label{AA-3}
\sum_{k=k_0+1}^{\infty} \varphi_\delta (c^{-k}) c^k (c - 1)= \sum_{k=k_0+1}^{\infty} \varphi_\delta (c^{-k}) \int_{c^{-k-1}}^{c^{-k}} t^{-2} \, dt\\[6pt]
\ge    \sum_{k=k_0+1}^{\infty}\int_{c^{-k-1}}^{c^{-k}} \varphi_\delta(t) t^{-2} \, dt = \int_{0}^{c^{-k_0 - 1}} \varphi_\delta(t) t^{-2} \, dt.
\end{multline}
It follows from \eqref{p4-K}, \eqref{AA-2}, and \eqref{AA-3} that, for $\delta < \tilde \delta(\eps)$,
\begin{equation*}
\Lambda_\delta(u_\delta) \ge \frac{1}{c}  c_1^{-1} \big[\C(c_1 \bvarphi_1) - \eps \big]  \int_{0}^{c^{-k_0 - 1}} \varphi_\delta(t) t^{-2} \, dt.
\end{equation*}
Note that
\begin{equation*}
\lim_{\delta \to 0} \int_{0}^{c^{-k_0 - 1}} \varphi_\delta(t) t^{-2} \, dt = \lim_{\delta \to 0}  \int_0^{c^{-k_0 - 1}/\delta} \varphi(t) t^{-2} \, dt =  \int_0^{\infty} \varphi(t) t^{-2} \, dt =
 \gamma_d^{-1} \mbox{ by } \eqref{cond-varphi3}.
\end{equation*}
On the other hand, by  \eqref{cond-varphi3} applied with $c_1 \bvarphi_1$, we have
$$
\gamma_d c_1 \int_{1}^\infty t^{-2} \, d t = \gamma_d c_1 = 1.
$$
We derive that
\begin{equation*}
\liminf_{\delta \to 0} \Lambda_\delta(u_\delta) \ge \frac{1}{c}(\C(c_1 \bvarphi_1) - \eps).
\end{equation*}
Since $c>1$ and $\eps > 0$ are arbitrary, we obtain \eqref{main-K}.  \proofend

\medskip 
Theorem~\ref{thm-K} suggests the following question
\begin{question} Assume that $\varphi, \, \psi \in {\cal A}$ satisfy
\begin{equation}
\varphi \ge \psi \mbox{ near 0 } (\mbox{resp. } \varphi = \psi \mbox{ near } 0). 
\end{equation}
Is it true that 
$$
K(\varphi) \ge K (\psi) \; \big(\mbox{resp. } K(\varphi) = K(\psi) \big)? 
$$

\end{question}

\medskip
We conclude this section with a simple observation
\begin{proposition} The set ${\cal A}$ is convex and the function $\varphi \mapsto \C (\varphi)$ is concave on ${\cal A}$. Moreover, $t \mapsto \C (t \varphi + (1 -t) \psi)$ is continuous on $[0, 1]$ for all $\varphi, \psi \in {\cal A}$. In particular, $\C({\cal A})$ is an interval.
\end{proposition}

This is an immediate consequence of the fact that
$$
\C(\varphi) = \inf \liminf_{\delta \to 0} \Lambda_{\delta} (v_\delta, \varphi)
$$
and that $\varphi \mapsto \Lambda_{\delta}(v_\delta, \varphi)$ is linear for fixed $\delta$.

\subsection{Proof of the fact that $K(c_1 \bvarphi_1) < 1$ for every $d \ge 1$}
\label{sect3.6}

In view of Theorem~\ref{thm-gamma}, it suffices to construct a bounded domain $\Omega \subset \mR^d$, a function $u \in BV(\Omega)$ with $\int_{\Omega}|\nabla u| = 1$, a sequence $\delta_n \to 0$, and a sequence $(u_n) \subset L^1(\Omega)$ such that $u_n \to u$ in $L^1(\Omega)$ and 
\begin{equation}\label{HH0}
\limsup_{n \to + \infty} \Lambda_{\delta_n} (u_n, c_1 \bvarphi_1) < 1. 
\end{equation}
We take $\Omega = Q$, $\delta_n = 1/n$, $u(x) = x_1$ where $x = (x_1, x')$ with $x_1  \in (0, 1)$ and $x' \in Q' = (0, 1)^{d-1}$, and 
$$
u_n(x) = i/ n \mbox{ if } i/ n \le x_1 < (i+1)/ n \mbox{ for } 0 \le i \le n-1. 
$$
Clearly $u_n \to u$ in $L^1(Q)$ as $n \to + \infty$. 

It follows from the definition of $u_n$ and $u$ that for $(x, y) \in Q^2$, 
\begin{equation}
\mbox{ if } |u_n(x) - u_n(y)|  > 1/ n, \mbox{ then } |u(x) - u(y)| > 1/n, 
\end{equation}
which implies that 
\begin{equation}
A_n := \big\{(x, y) \in Q^2;  \; |u_n(x) - u_n(y)| > 1/n \big\} \subset B_n : = \big\{(x, y) \in Q^2;  \; |u(x) - u(y)| > 1/n \big\}. 
\end{equation}
Thus, by the definition of $\Lambda_{1/n}$, 
\begin{equation}\label{HH1}
\Lambda_{1/n} (u, c_1 \bvarphi_1) - \Lambda_{1/n} (u_n, c_1 \bvarphi_1) = \frac{c_1}{n} \mathop{\iint}_{B_n \setminus A_n} \frac{1}{|x - y|^{d+1}} \, dx \, dy. 
\end{equation}
For $0 \le i \le n-2$ and $n \ge 3$, set 
$$
Z_{i, n} = \big\{(x, y) \in Q^2; \; i/ n < x_1 < (i+1/2)/ n \mbox{ and } (i+3/2)/ n < y_1 < (i+2)/ n \big\}. 
$$
Note that if $(x, y) \in Z_{i, n}$ we have 
$$
u(y) - u(x) = y_1 - x_1 > 1/n \quad \mbox{ and } \quad u_n(y) - u_n(x) = (i+1)/ n - i/ n = 1/ n, 
$$
so that 
$$
Z_{i, n} \subset B_n \setminus A_n \mbox{ for } 0 \le i \le n-2. 
$$
On the other hand if $(x, y) \in Z_{i, n}$ we have 
$$
|x - y|^2 = |x_1 - y_1|^2 + |x' - y'|^2 \le 4/ n^2 + |x' - y'|^2, 
$$
and consequently
\begin{multline}\label{HH2}
\mathop{\iint}_{B_n \setminus A_n}  \frac{1}{|x-y|^{d+1}} \, dx \, dy \ge \sum_{i=0}^{n-2} \mathop{\iint}_{Z_{i, n}} \frac{1}{|x-y|^{d+1}} \, dx \, dy \\[6pt] \ge \sum_{i=0}^{n-2} \frac{1}{4n^2} \mathop{\iint}_{Q' \times Q'} \frac{1}{\big( (4/ n^2) + |x' - y'|^2\big)^{\frac{d+1}{2}}} \, dx' d y' \sim \frac{1}{n} \mathop{\iint}_{Q' \times Q'} \frac{1}{\big( (4/ n^2) + |x' - y'|^2\big)^{\frac{d+1}{2}}} \, dx' d y'. 
\end{multline}
Recall the (easy and) standard fact that 
\begin{equation}\label{HH3}
\mathop{\iint}_{Q' \times Q'} \frac{1}{\big(a^2+ |x' - y'|^2\big)^{\frac{d+1}{2}}} \, dx' d y' \sim 1/ a^2 \mbox{ for small } a. 
\end{equation}
Combining \eqref{HH1}, \eqref{HH2}, and \eqref{HH3} yields 
\begin{equation}\label{HH4}
\Lambda_{1/n} (u, c_1 \bvarphi_1) - \Lambda_{1/n} (u_n, c_1 \bvarphi_1)  \ge C_d  > 0. 
\end{equation}
From   Proposition~\ref{thm-Sobolev1p=1}, we get
\begin{equation}\label{HH5}
\lim_{n \to \infty} \Lambda_{1/n} (u) = \int_{Q} |\nabla u| = 1. 
\end{equation}
The desired result \eqref{HH0} follows from \eqref{HH4} and \eqref{HH5}. 
\proofend

\section{Compactness results. Proof of Theorems~\ref{thm-compact-p=1} and \ref{thm-compact2}}\label{sec-pro}

The following  subtle estimate from \cite[Theorem 1]{NgSob3} (with roots in \cite{BourNg}) plays a crucial role in the proof of Theorems~\ref{thm-compact-p=1} and \ref{thm-compact2}.

\begin{lemma}\label{proL2}
Let $d \ge 1$, $B_1$ be the unit ball (or cube),
and $u \in L^1(B_1)$.
There exists a positive constant $C_d$, depending only on $d$, such that
\begin{equation}\label{est-e-2}
\int_{B_1} \int_{B_1} |u(x) - u(y)| \, dx \, dy \le C_d \Big(
\mathop{\int_{B_1} \int_{B_1}}_{|u(x) - u(y)|
> 1} \frac{1}{|x-y|^{d+1}} \, dx \, dy +
1 \Big).
\end{equation}
\end{lemma}

By scaling, we obtain, for any ball or cube $B$,
\begin{equation}\label{est-e-2-1}
\int_{B} \int_{B} |u(x) - u(y)| \, dx \, dy \le C_d \Big( |B|^{1 + 1/d}
\mathop{\int_{B} \int_{B}}_{|u(x) - u(y)|
> 1} \frac{1}{|x-y|^{d+1}} \, dx \, dy +
|B|^2 \Big).
\end{equation}

The reader can find in \cite{BrezisNguyen} a connection between these inequalities and  the  $VMO / BMO$ spaces. 

\medskip Here is a question related to Lemma~\ref{proL2}: 

\begin{question} Is it true that 
\begin{equation}
\mathop{\int_{B_1} \int_{B_1}}_{|u(x) - u(y)| > 1} |u(x) - u(y)| \, dx \, dy  \le C_d \mathop{\int_{B_1} \int_{B_1}}_{|u(x) - u(y)| > 1} \frac{1}{|x - y|^{d+1}} \, dx \, dy \quad \forall \, u \in L^1(B_1) ? 
\end{equation}

\end{question}

\subsection{Proof of Theorem~\ref{thm-compact-p=1}}
In this subsection we fix $\delta = 1$. We recall the notation from \eqref{def-Lambda}
\begin{equation*}
\Lambda (u, \Omega) = \Lambda_{1} (u, \Omega) =  \int_{\Omega} \int_{\Omega}   \frac{\varphi(|u(x) - u(y)|)}{|x - y|^{d+ 1}} \, dx \, dy.
\end{equation*}

Here is an immediate consequence of Lemma~\ref{proL2}.

\begin{lemma}\label{lem-Gen} Let $B$ be a ball (or cube) and $\varphi$ be such that  \eqref{cond-varphi-decreasing} and
\eqref{cond-varphi-positive} hold, and let  $u \in L^1(B)$. We have
\begin{equation}\label{est1}
\int_B \int_B |u(x) - u(y)| \, dx \, dy \le C_d \Big\{ \frac{\lambda}{\varphi(\lambda)} |B|^{1 + 1/d}
\Lambda(u, B) +  \lambda |B|^{2} \Big\}, \quad \forall \, \lambda > 0.
\end{equation}
\end{lemma}

Assume $\Omega$ is bounded. Denote $\Gamma = \partial \Omega$, and set
\begin{equation*}
\Omega_t =\{x \in \mR^d; \; \mbox{dist }(x, \Omega) < t\}.
\end{equation*}
For $t$ small enough, every  $x \in \Omega_t \setminus \Omega$  can be uniquely  written as
\begin{equation}\label{small-t}
x = x_{\Gamma} + s \nu(x_\Gamma),
\end{equation}
where $x_\Gamma$ is the projection of $x$ onto $\Gamma$, $s = \mbox{dist }(x, \Gamma)$,  and $\nu(y)$ denotes the outward normal unit vector at $y \in \Gamma$.

\begin{lemma} \label{lem-extension} Assume that $\Omega$ is bounded. Fix $t > 0$  small enough such that \eqref{small-t} holds for any $x \in \Omega_t$. There exists an extension $U$ of $u$  in $\Omega_t$ such that
\begin{equation*}
\| U \|_{L^1(\Omega_t)} \le C \| u \|_{L^1(\Omega)} \quad \mbox{ and } \quad \Lambda (U, \Omega_t) \le C \Lambda(u,  \Omega),
\end{equation*}
for some positive constant $C$ depending only on $\Omega$.
\end{lemma}

\noindent{\bf Proof.} Define
\begin{equation*}
U(x) = \left\{\begin{array}{cl} u(x) & \mbox{ if } x \in \Omega, \\[6pt]
u\big(x_\Gamma- s \nu(x_\Gamma)\big) & \mbox{ if } x \in \Omega_t \setminus \Omega.
\end{array}\right.
\end{equation*}
It is clear that
\begin{equation*}
\| U \|_{L^1(\Omega_t)} \le C \| u \|_{L^1(\Omega)}.
\end{equation*}
In this proof $C$ denotes a positive constant depending only on $\Omega$. It remains to prove that
\begin{equation*}
 \Lambda(U,  \Omega_t) \le C \Lambda(u,  \Omega).
\end{equation*}
By the definition of $\Lambda$ in \eqref{def-Lambda}, it suffices to show that
\begin{equation}\label{main-part}
\int_{\Omega_t} \, dy  \int_{\Omega_t \setminus \Omega}  \frac{\varphi(|U(x) - U(y)|)}{|x - y|^{d+1}} \, dx  \le C \int_{\Omega} \int_{\Omega}  \frac{\varphi(|u(x) - u(y)|)}{|x - y|^{d+1}} \, dx \, dy. 
\end{equation}
If $x \in \Omega_t \setminus \Omega$ and $y \in \Omega_t \setminus \Omega$, then
\begin{equation*}
U\big(x_\Gamma + s_1 \nu(x_\Gamma)\big) - U\big(y_\Gamma + s_2 \nu(y_\Gamma)\big) = u\big(x_\Gamma - s_1 \nu(x_\Gamma)\big) - u\big(y_\Gamma - s_2 \nu(y_\Gamma)\big),
\end{equation*}
and
\begin{equation*}
\Big| \big(x_\Gamma + s_1 \nu(x_\Gamma)\big) - \big(y_\Gamma + s_2 \nu(y_\Gamma) \big) \Big| \ge C \Big| \big(x_\Gamma - s_1 \nu(x_\Gamma)\big) - \big(y_\Gamma - s_2 \nu(y_\Gamma) \big) \Big|,
\end{equation*}
and,  if $x \in \Omega_t \setminus \Omega$ and $y \in \Omega$, then
\begin{equation*}
U\big(x_\Gamma + s_1 \nu(x_\Gamma)\big) - U(y) = u\big(x_\Gamma - s_1 \nu(x_\Gamma)\big) - u(y),
\end{equation*}
and
\begin{equation*}
\Big| \big(x_\Gamma + s_1 \nu(x_\Gamma)\big) - y \Big| \ge C \Big| \big(x_\Gamma - s_1 \nu(x_\Gamma)\big) - y\big) \Big|.
\end{equation*}
Hence \eqref{main-part} holds. \proofend

\medskip

We are ready to present the

\medskip

\noindent {\bf Proof of Theorem~\ref{thm-compact-p=1}.} It suffices to consider the case where $\Omega$ is bounded. By Lemma~\ref{lem-extension}, one only needs to prove that up to a subsequence, $u_n \to u $ in $L^1_{\loc}(\Omega)$.   For a cube  $Q$ in $\Omega$, define
\begin{equation*}
F(u, Q) = \int_{Q} \int_{Q} \frac{ \varphi\big( |u(x) - u(y)|  \big)}{|x - y|^{d+1}} \, dx \, dy + |Q|.
\end{equation*}
Since,  by  \eqref{est1},
\begin{align*}   \frac{1
}{|Q|} \int_Q \int_Q |u(x) - u(y)| \, dx \, dy \le   C_d \Big\{ \frac{\lambda}{\varphi(\lambda)} |Q|^{1/d}
\int_{Q} \int_{Q} \frac{ \varphi \big( |u(x) - u(y)|  \big)}{|x - y|^{d+1}} \, dx \, dy + \lambda |Q| \Big\},  \end{align*}
it follows that
\begin{align}\label{basicE}
\frac{1}{|Q|} \int_Q \int_Q |u(x) - u(y)| \, dx \, dy    \le \rho(|Q|) F(u, Q),
\end{align}
where
\begin{equation*}
\rho(t) := C_d \inf_{\lambda > 0} \Big( \frac{\lambda t^{1/d}}{\varphi(\lambda)} + \lambda\Big).
\end{equation*}
It is clear that $\rho$ is non-decreasing and, by \eqref{cond-varphi-positive},
\begin{equation}\label{to}
\lim_{t \to 0} \rho(t) = 0.
\end{equation}
For $\eps > 0$ and $n \in \mN$, set
\begin{equation*}
u_{n, \eps}(x) = \frac{1}{\eps^d}\int_{Q_\eps(x)} u_n( y) \, dy,
\end{equation*}
where $Q_\eps(x)$ is the cube centered at $x$ of side $\eps$.
Fix an arbitrary cube $Q \subset \subset \Omega$. We claim that
\begin{equation}\label{compact-Thm2}
\int_{Q} |u_n(x) - u_{n, \eps} (x)| \, dx  \to 0 \mbox{ as } \eps \to 0, \mbox{ uniformly in } n.
\end{equation}
Let $\eps$ be small enough such that
$Q + \eps [-1, 1]^d \subset \Omega$.
Then there exists a finite family  $\Big(Q(j) \Big)_{i \in J} $ of disjoint open $\eps$-cubes such that
\begin{equation*}
 Q \subset \mbox{ interior } \Big( \bigcup_{j \in J}  \overline{ Q(j) } \Big) \subset \bigcup_{j \in J}  \overline{ 2 Q(j) } \subset \Omega,
\end{equation*}
and thus $\mbox{card} J \sim 1/ \eps^d$.  Here and in what follows $aQ(j)$ denotes the cube which has the same center as $Q(j)$ and of $a$ times its length.
Applying \eqref{basicE} we have
\begin{align}\label{sum-0}
\int_{Q(j)} |u_n(x) - u_{n, \eps} (x)| \, dx \le & \frac{C}{|2Q(j)|} \int_{2 Q(j)} \int_{2Q(j)} |u_n(x) - u_{n} (y)| \, dx \, dy \nonumber
\\[6pt] \le & C  \rho(2^d \eps^d) F \big(u_n, 2 Q(j) \big),
\end{align}
since $Q(j) + \eps [-1/2, 1/2]^d \subset 2 Q(j)$.
Note that the family $2 Q(j)$ is not disjoint, however, they have a finite number of  overlaps (depending only on $d$).
Therefore, for any $f \ge 0$,
\begin{equation}\label{sum-1}
\sum_{j} \int_{2 Q(j)} \int_{2 Q(j)} f \le C\int_{\Omega} \int_\Omega f,
\end{equation}
Summing with respect to $j$ in \eqref{sum-0}, we derive from \eqref{sum-1} that
\begin{equation}\label{sum-2}
 \int_{Q} |u_n(x) - u_{n, \eps} (x)| \, dx \le C \rho(2^d \eps^d) F\big(u_n, \Omega\big).
\end{equation}
Using  \eqref{compactness-p=1}, \eqref{to},  and \eqref{sum-2}, we obtain  \eqref{compact-Thm2}.
It follows from \eqref{compact-Thm2} and  a standard argument (see, e.g., the proof of the theorem of Riesz-Frechet-Kolmogorov in \cite[Theorem 4.26]{BrAnalyseEnglish}) that   there exists a subsequence
$(u_{n_k})$  of $(u_n)$ and $u \in L^1_{\loc}(\Omega)$ such that $(u_{n_k})$ converges to $u$ in $L^{1}_{\loc}(\Omega)$.\proofend

\begin{remark}  \fontfamily{m} \selectfont
Using Theorem~\ref{thm-compact-p=1}, one can prove that $\Lambda$ is lower semi-continuous with respect to weak convergence in $L^q$ for any $q \ge 1$.
\end{remark}

\subsection{Proof of Theorem~\ref{thm-compact2}}

It suffices to consider the case where $\Omega$ is bounded. By Lemma~\ref{lem-extension}, one only needs to prove that up to a subsequence, $u_n \to u $ in $L^1_{\loc}(\Omega)$.
Fix $\lambda_0  > 0$ such that $\varphi(\lambda_0) > 0$. Without loss of generality, one may assume that $\lambda_0 =1$. From \eqref{bd}, we have
\begin{equation*}
\sup_{n} \mathop{\int_{\Omega} \int_{\Omega}}_{|u_n(x) - u_n(y)| > \delta_n} \frac{\delta_n}{|x - y|^{d+1}} \, dx \, dy \le C.
\end{equation*}
We now follow the same strategy as in the proof of Theorem~\ref{thm-compact-p=1}.
Define
\begin{equation}\label{def-u-eps-Thm3}
u_{n, \eps}(x) = \frac{1}{\eps^d } \int_{Q_\eps(x)} u_n(y) \, dy,
\end{equation}
Here $Q_\eps(x)$ is the cube centered at $x$ of side $\eps$.
Fix an arbitrary cube $Q \subset \subset \Omega$. We claim that
\begin{equation}\label{compact-Thm3}
\lim_{\eps \to 0} \limsup_{n} \int_{Q} |u_n(x) - u_{n, \eps} (x)| \, dx = 0.
\end{equation}
Let $\eps$ be small enough such that
$Q + \eps [-1, 1]^d \subset \Omega$.
Then there exists a finite family  $\Big(Q(j) \Big)_{i \in J} $ of disjoint open $\eps$-cubes such that
\begin{equation*}
 Q \subset \mbox{ interior } \Big( \bigcup_{j \in J}  \overline{ Q(j) } \Big) \subset \bigcup_{j \in J}  \overline{ 2 Q(j) } \subset \Omega.
\end{equation*}
We have
\begin{align}\label{sum-0-Thm3}
\int_{Q(j)} |u_n(x) - u_{n, \eps} (x)| \, dx \le & \frac{C}{|2Q(j)|} \int_{2 Q(j)} \int_{2Q(j)} |u_n(x) - u_{n} (y)| \, dx \, dy, 
\end{align}
since $Q(j) + \eps [-1/2, 1/2]^d \subset 2 Q_j$. By \eqref{sum-0-Thm3} and  \eqref{est-e-2-1} with $B = 2 Q(j)$, we have
\begin{equation}\label{sum-1-Thm3}
\int_{Q(j)} |u_n(x) - u_{n, \eps} (x)| \, dx  \le  \C \Big(\eps
\mathop{\int_{2 Q(j)} \int_{2Q(j)}}_{|u_n(x) - u_n(y)|
> \delta_n} \frac{\delta_n}{|x-y|^{d+1}} \, dx \, dy + \delta_n  \eps^d \Big).
\end{equation}
Summing with respect to $j$ in \eqref{sum-1-Thm3}, we obtain
\begin{equation*}
 \int_{Q} |u_n(x) - u_{n, \eps} (x)| \, dx \le C (\eps + \delta_n).
\end{equation*}
Clearly, for fixed $n$,
\begin{equation*}
\lim_{\eps \to 0}\int_{Q} |u_n(x) - u_{n, \eps} (x)| \, dx = 0.
\end{equation*}
Therefore \eqref{compact-Thm3} holds and we conclude as in the proof of Theorem~\ref{thm-compact-p=1}. \proofend


\section{Some functionals related to Image Processing}\label{sec-image}

Given $q \ge 1$, $\lambda > 0$, $\delta > 0$, $d \ge 1$, $\Omega$  a
smooth bounded open subset of $\mR^d$, and $ f \in L^q(\Omega)$, consider the non-local, non-convex functional defined on $L^{q}(\Omega)$ by
\begin{equation}\label{defE-S}
E_{\delta}(u) = \lambda \int_{\Omega}| u - f|^q + \Lambda_{\delta}(u): = \lambda \int_{\Omega}| u- f|^q +
\delta \int_\Omega \int_\Omega \frac{\varphi(|u(x) - u(y)| / \delta) }{|x - y|^{ d  +  1}} \, dx \, dy.
\end{equation}
Our goal in this section is twofold. In the first part, we investigate the existence of a minimizer for $E_{\delta}$ ($\delta$ is fixed)  and then we study the behavior of these minimizers (or almost minimizers) as $\delta \to 0$. In the second part, we explain how these results are connected to Image Processing.

\subsection{Variational problems  associated with $E_{\delta}$}

We start with an immediate consequence of Theorem~\ref{thm-compact-p=1}.

\begin{corollary}\label{cor-image1p=1} Let $\delta > 0$ be fixed. Assume that  $\varphi \in {\cal A}$ satisfies \eqref{cond-varphi-positive}. There exists $u \in L^q(\Omega)$ such that
\begin{equation*}
E_{\delta}(u) = \inf_{w \in L^q(\Omega)} E_{\delta}(w).
\end{equation*}
\end{corollary}

\medskip
As we know from Theorem~\ref{thm-gamma}, under  assumptions \eqref{cond-varphi-0}-\eqref{cond-varphi3},  $\Lambda_\delta$ $\Gamma$-converges to $K \int_{\Omega} |\nabla \cdot |$ as $\delta \to 0$,  for some constant $ 0 < K \le 1$. Therefore, one may expect that the minimizers of $E_\delta$ converge to the unique minimizer in $L^{q}(\Omega)$ of $E_{0}$, where
\begin{equation}
E_{0}(w) = \lambda \int_{\Omega} |w - f|^q + K \int_{\Omega} |\nabla w|.
\end{equation}
If one does not assume  \eqref{cond-varphi-positive} one can not apply Corollary~\ref{cor-image1p=1}, and minimizers of $E_\delta$ might not exist; however, one can always consider almost minimizers. Here is slight generalization of Theorem~\ref{thm-image}.

\begin{theorem} \label{thm-4} Let $d \ge 1$, $q \ge 1$, $\Omega$ be a smooth  bounded open subset of $\mR^{d}$,  $f \in L^{q}(\Omega)$,  and $\varphi \in {\cal A}$. Let $(\delta_{n}), \; (\tau_{n})$ be two positive sequences converging to 0 as $n \to \infty$ and $u_{n} \in L^{q}(\Omega)$ be such that
\begin{equation}\label{defunk-image}
E_{\delta_{n}} (u_{n}) \le \inf_{u \in L^{q}(\Omega)} E_{\delta_{n}} (u) + \tau_{n}.
\end{equation}
Then $u_{n} \to u_{0}$ in $L^{q}(\Omega)$ where $u_{0}$ is the unique minimizer of the functional $E_{0}$ defined on $L^{q}(\Omega)\cap BV(\Omega)$ by
\begin{equation*}
E_{0}(u) : = \lambda \int_{\Omega} |u - f|^{q} + K \int_{\Omega} |\nabla u|,
\end{equation*}
and $0 < K \le 1$ is the constant in Theorem~\ref{thm-gamma}.
\end{theorem}

\noindent {\bf Proof.} It is clear that $(u_n)$ is bounded in $L^1(\Omega)$.  By Theorem~\ref{thm-compact2}, there exists a subsequence $(u_{n_k})$ which converges to some  $u_0$ a.e. and in $L^1(\Omega)$.  It follows from  Fatou's lemma and Property (G1) in Section~\ref{sec-gamma} that
\begin{equation}\label{est1-image}
E_{0}(u_0) \le \liminf_{k \to \infty} E_{\delta_{n_k}}(u_{n_k}).
\end{equation}

We will prove that $u_0$ is the unique minimizer of $E_0$ in $L^q(\Omega) \cap BV(\Omega)$.  Let $v \in L^{q}(\Omega) \cap BV(\Omega)$ be the unique minimizer of $E_0$. Applying Theorem~\ref{thm-gamma}, there exists $v_{n} \in L^{1}(\Omega)$ such that $v_{n} \to v$ in $L^{1}$ (without loss of generality,  one may assume that $v_{n} \to v$ a.e.) and
\begin{equation*}
\limsup_{n \to \infty} \Lambda_{\delta_n}(v_n) \le K \int_{\Omega} |\nabla v|.
\end{equation*}
For $A > 0$, recall the notation $T_A$ defined in \eqref{def-TA}. From \eqref{cond-varphi-decreasing}, we have
\begin{equation*}
\Lambda_{\delta_{n}}(T_A v_{n}) \le \Lambda_{\delta_{n}}(v_n).
\end{equation*}
By definition of $u_n$, we obtain 
\begin{equation*}
E_{\delta_{n}} (u_{n})   \le  \lambda
\int_{\Omega} |T_A v_{n} - f|^q + \Lambda_{\delta_{n}} (T_A v_{n}) + \tau_n \le \lambda \int_{\Omega} |T_A v_{n} - f|^q + \Lambda_{\delta_{n}} (v_{n}) + \tau_n.
\end{equation*}
Letting $n \to \infty$ yields 
\begin{equation*}
E_0(u_0) \le \liminf_{n \to \infty} E_{\delta_{n}} (u_{n})   \le \lambda \int_{\Omega} |T_A v - f|^q + K \int_{\Omega} |\nabla v|.
\end{equation*}
As  $A \to \infty$, we find 
\begin{equation}\label{important}
E_0(u_0) \le  \lambda \int_{\Omega} |v - f|^q + K \int_{\Omega} |\nabla v|.
\end{equation}

 This implies that $u_0$ is the unique minimizer of $E_0$.

\medskip
We next prove that $u_n \to u_0$ in $L^q$. Since
\begin{equation*}
E_0(u_0) \ge  \limsup_{n \to \infty} E_{\delta_{n}} (u_{n}) \ge E_0(u_0)
\end{equation*}
by \eqref{important},  and
\begin{equation*}
\liminf_{n \to \infty} \Lambda_{\delta_n} (u_n) \ge K \int_{\Omega} |\nabla u_0|,
\end{equation*}
by Theorem~\ref{thm-gamma},  we have
\begin{equation*}
\lim_{n \to \infty} \int_\Omega |u_{n} - f|^q = \int_{\Omega} |u_0 - f|^q.
\end{equation*}
In addition we know that $u_n - f \to u_0- f$ a.e. in $\Omega$. Therefore $u_n - f \to u_0-f $ in $L^q(\Omega)$; thus $u_n \to u_0$ in $L^q(\Omega)$. The proof is complete. \proofend

\begin{remark} \fontfamily{m} \selectfont   In case a Lavrentiev - type gap does occur (see Open problem~\ref{OP3} and the subsequent comments) it would be interesting to investigate what happens in Theorem~\ref{thm-4} if ${E_\delta}_{|L^q(\Omega)}$ is replaced by ${E_\delta}_{|C^0(\bar \Omega)}$ (with numerous possible variants). 
\end{remark}

%

\subsection{Connections with Image Processing}

 A fundamental challenge in Image Processing is to improve images of poor quality. Denoising is an immense subject, see,  e.g.,  the excellent survey by A.~Buades, B.~Coll and J.~M.~Morel \cite{BuadesCollMorelReview1}. One possible strategy is to introduce a filter $F$ and use  a variational formulation
 \begin{equation}\label{image}
 \min_{u} \Big\{ \lambda \int_{\Omega} |u - f|^{2} + F(u) \Big\},
 \end{equation}
or, alternatively, the associated Euler equation
\begin{equation}\label{Euler}
2 \lambda (u - f) + F'(u) = 0.
\end{equation}
Here $f$ is the given image of poor quality, $\lambda > 0$ is the fidelity parameter (fixed by experts) which governs how much filtering is desirable. Minimizers of \eqref{image} (or solutions to \eqref{Euler}) are the denoised images.

\medskip
Many types of filters are used in Image Processing. Here are three popular ones. The first one is the celebrated (ROF) filter due L.~Rudin, S.~Osher, and E.~Fatemi \cite{RudinOsherFatemi}:
\begin{equation*}
F(u)= \int_{\Omega} |\nabla u|
\end{equation*}
(see also \cite{ChambolleLions, Chan, ChoksiFonseca, HaddadMeyer}). The corresponding minimization problem is
\begin{equation*}
(ROF) \quad \quad \min_{u \in L^{2}(\Omega)}  \Big\{ \lambda \int_{\Omega} |u - f|^{2} + \int_{\Omega} |\nabla u| \Big\}.
\end{equation*}
The functional in $(ROF)$ is strictly convex. It follows from standard Functional Analysis that, given $f \in L^{2}(\Omega)$, there exists a unique minimizer $u_{0} \in BV(\Omega) \cap L^{2}(\Omega)$.

The second filter,  due to  G.~Gilboa and S.~Osher \cite{GilboaOsher1} (see also \cite{GilboaOsher2}),  is
\begin{equation*}
 F(u) = \int_{\Omega} \Big( \int_{\Omega} \frac{|u(x) - u(y)|^{2}}{|x - y|^{2}} w(x, y)\, dy\Big)^{1/2} \, dx,
\end{equation*}
where $w$ is a given weight function.
The corresponding minimization problem is
\begin{equation*}
(GO) \quad \quad  \min_{u \in L^{2}(\Omega)}  \Big\{ \lambda \int_{\Omega} |u - f|^{2} +  \int_{\Omega} \Big( \int_{\Omega} \frac{|u(x) - u(y)|^{2}}{|x - y|^{2}} w (x, y) \, dy\Big)^{1/2} \, dx  \Big\}.
\end{equation*}
The functional  in $(GO)$ is strictly convex.  Again by standard Functional Analysis, there exists a unique minimizer $u_0$ of $(GO)$.   One can prove (see \cite{BN-Two}) that if $w(x, y) = \rho_\eps(|x -y|)$, where $(\rho_\eps)$ is a sequence of mollifiers as in Remark~\ref{rem-BBM}, then the corresponding minimizers $(u_\eps)$ of $(GO_\eps)$ (i.e., (GO) with $w(x, y) = \rho_\eps(|x -y|)$) converge, as $\eps \to 0$,  to the unique solution of the $(ROF_k)$ problem


\begin{equation*}
(ROF_k) \quad \min_{u \in L^2(\Omega) \cap BV(\Omega)}  \Big\{ \lambda \int_{\Omega}|u-f|^2 + k \int_{\Omega} |\nabla u| \Big\},
\end{equation*}
where
\begin{equation*}
k= \Big( \int_{\mS^{d-1}} |\sigma \cdot e|^2 \, d \sigma \Big)^{1/2},
\end{equation*}
for some $e \in \mS^{d-1}$. The proof in \cite{BN-Two} is strongly inspired by the results of J.~Bourgain, H.~Brezis, and P.~Mironescu in \cite{BBM},  A. Ponce \cite{Po}, and
G.~Leoni and D.~Spector \cite{LeoniSpector1}.

In a similar spirit,  G.~Aubert and P.~Kornprobst  in \cite{AubertKornprobst} have proposed to use the filter
\begin{equation*}
 F(u) = I_{\eps} (u) = \int_{\Omega} \int_{\Omega} \frac{|u(x) - u(y)|}{|x - y|} \rho_{\eps}(|x -y|) \, dx \, dy,
\end{equation*}
and the corresponding minimization problem is
\begin{equation*}
(AK_{\eps}) \quad \quad \min_{u \in L^{2}}  \Big\{ \lambda \int_{\Omega} |u - f|^{2} +  \int_{\Omega} \int_{\Omega} \frac{|u(x) - u(y)|}{|x - y|} \rho_{\eps}(|x -y|) \, dx \, dy \Big\}.
\end{equation*}
As above, $(AK_\eps)$ admits a unique minimizer $u_\eps$ and, as $\eps \to 0$, $(u_\eps)$ converges to the solution of $(ROF_{\gamma_d})$ where $\gamma_d$ is the constant defined in \eqref{cond-varphi3}.

\medskip

The third type of filter was introduced in the pioneering works of L.~S.~Lee~\cite{Lee} and  L. P.~ Yaroslavsky  (see \cite{Yaroslavsky1, YaroslavskyEden});
more details can be found in  the expository paper by A. ~Buades, B. Coll, and J.~M.~Morel
\cite{BuadesCollMorelReview1}; see also \cite{BuadesCollMorel2,  Morel1,PKTD, SmithBrady}) where the terms ``neighbourhood filters", ``non-local means" and ``bilateral filters" are used.  Originally, they were not formulated as variational problems. In an important paper K.~Kindermann, S.~Osher and P.~W.~Jones \cite{KindermannOsherJones} showed that some of these filters come from the Euler-Lagrange equation of a minimization problem where the  functional $F$ has the form
\begin{equation*}
F(u) =  \int_{\Omega} \int_{\Omega} \varphi \big( |u(x) - u(y)|/ \delta \big) w (|x -y|) \, dx \, dy,
\end{equation*}
$\delta> 0$ is a fixed small parameter, $\varphi$ is a given {\bf non-convex} function, and $w\ge 0$ is a weight function.  The corresponding minimization problem is
\begin{equation*}
(YNF_{\delta}) \quad \quad \min_{u \in L^{2}}  \Big\{ \lambda \int_{\Omega} |u - f|^{2} +  \int_{\Omega} \int_{\Omega} \varphi \big( |u(x) - u(y)|/ \delta \big) w (|x -y|) \, dx \, dy \Big\}.
\end{equation*}
Here are some examples of $\varphi's$ and $w's$ used in Image Processing see,  e.g.,  \cite[Section 3]{KindermannOsherJones}:

\begin{itemize}

\item[i)]$\varphi = \bvarphi_2$ or $\varphi = \bvarphi_3$ (from the list of examples in the Introduction).

\item[ii)] $w =1$ or
\begin{equation*} w(t)  = \left\{
\begin{array}{cl}
1 & \mbox{ if }  t< \rho, \\[6pt]
0 & \mbox{ otherwise},
\end{array}\right.
\end{equation*}
for some $\rho > 0$.
\end{itemize}

In this paper,  we suggest a new example for $w$:
\begin{equation}\label{new-choice}
w(|x - y|) = \frac{1}{|x-y|^{d+1}}.
\end{equation}
Taking $\lambda \sim 1/ \delta$, more precisely $\lambda = \gamma/ \delta$, we are led to
the minimization problem:
\begin{equation}\label{new-F}
\min_{u \in L^{2}}  \Big\{ \gamma \int_{\Omega} |u - f|^{2} + \Lambda_{\delta}(u) \Big\}.
\end{equation}
Up to now, there was {\bf no} rigorous analysis whatsoever for problems of the form $(YNF_\delta)$. Even the existence of minimizers in $(YNF_\delta)$, for fixed $\delta$,  was lacking. Our contributions for the new choice of $w$ in \eqref{new-choice} are twofold:
\begin{enumerate}
\item Existence of minimizers for \eqref{new-F} under fairly general assumptions on $\varphi$ (Theorem~\ref{thm-compact-p=1}).

\item Asymptotic analysis as $\delta \to 0$:  $(YNF_{\delta}) \to (ROF_K)$ (Theorem~\ref{thm-image}).
\end{enumerate}


\appendix

\section{Appendix: Proof of  Pathology~\ref{pathology-2}} \label{A}
\renewcommand{\theequation}{A\arabic{equation}}
\renewcommand{\thelemma}{A\arabic{lemma}}
  \setcounter{equation}{0}  
  \setcounter{lemma}{0}  

We construct a function $u \in W^{1, 1}(0, 1)$ such that,  for $\varphi = c_1 \bvarphi_1$, 
\begin{equation}\label{ine-Fatou1-A}
\liminf_{\delta \to 0} \Lambda_{\delta} (u) = \int_{\Omega} |\nabla u| \quad \mbox{ and } \quad  \limsup_{\delta \to 0} \Lambda_{\delta} (u) = + \infty. 
\end{equation}

Set $x_n = 1 - 1/ n$ for $n \ge 1$.
Set $\delta_1 = 1/100$ and   $T_1 = e^{-\delta_1^{-1}}$.
Let $y_1$ be the middle point of the interval $(x_1, x_1 + T_1)$
and fix $0< t_1< T_1/ 4$  such that
$$
 \int_{x_1}^{y_1 - t_1} \, dx \int_{y_1 + t_1}^{x_1 + T_1} \frac{\delta_1}{|x - y|^2}  \, dy \ge 1.
$$
Since
\begin{equation}\label{A-infty}
 \int_{\alpha}^{\beta} \, dx \int_{\beta}^{\gamma} \frac{\delta_1}{|x - y|^2}  \, dy = + \infty,
\end{equation}
for all $\alpha < \beta < \gamma$, such a $t_1$ exists.  Define $u_1 \in W^{1, 1}(0, 1)$ by
$$
u_1(x) = \left\{ \begin{array}{cl}
\mbox{ constant in }  & [0, x_1], \\[6pt]
 \mbox{ affine in }  &  [x_1, y_1 - t_1 ],  \\[6pt]
\mbox{ affine in } & [y_1 - t_1 , y_1 + t_1],  \\[6pt]
\mbox{ affine in }  & [y_1 + t_1 , x_1 + T_1], \\[6pt]
\mbox{ constant in }  & [x_1 + T_1,  1], 
\end{array} \right.
\mbox{ and }
\left\{ \begin{array}{c}
u_1 (x_1)= 0, \\[6pt]
u_1(y_1 - t_1) =  \delta_1 / 3, \\[6pt]
u_1(y_1  + t_1) = 2 \delta_1 / 3, \\[6pt]
u_1(x_1  + T_1) =  \delta_1.
\end{array}\right.
$$
Assuming that $\delta_{k}$, $T_{k}$, $t_{k}$, and $u_k$ are constructed for $1 \le k \le n-1$ and for $n \ge 2$ such  that $u_k$ is Lipschitz. We then obtain  $\delta_n$, $T_n$, $t_n$, and $u_n$ as follows.
Fix  $0< \delta_{n}<  \delta_{n-1}/ 8$ sufficiently small such that
\begin{equation}\label{A-p1}
\Lambda_{2\delta_{n}} (u_{n-1}) + \int_{0}^{x_{n-1} + T_{n-1}} \, dx \int_{x_n}^1 \frac{2 \delta_n c_1}{|x-y|^2}  \, dy  \le \int_0^1 |u_{n-1}'| + 1/ n.
\end{equation}
Such a constant $\delta_n$ exists by Proposition~\ref{thm-Sobolev1p=1} (in fact $u_{n-1}$ is only Lipschitz; however Proposition~\ref{thm-Sobolev1p=1}  holds as well for Lipschitz functions, see also Proposition~\ref{pro-C}). Set $T_n = e^{- \delta_n^{-1}}$ and let $y_n$ be the middle point of the interval $(x_n, x_n + T_n)$
and fix $0< t_n< T_n/ 4$  such that
\begin{equation}\label{A-p2}
 \int_{x_n}^{y_n - t_n}  \, dx  \int_{y_n + t_n}^{x_n + T_n} \frac{\delta_n}{|x - y|^2} \, dy \ge n.
\end{equation}
Such a $t_n$ exists by \eqref{A-infty}.  Define a continuous function $w_n: [0, 1] \mapsto [0, 1]$, $n \ge 2$,  as follows 
$$
w_n(x) = \left\{ \begin{array}{cl}
\mbox{ constant in }  & [0, x_n], \\[6pt]
\mbox{ affine in }  &  [x_n, y_n - t_n ],  \\[6pt]
\mbox{ affine in } & [y_n - t_n , y_n + t_n],  \\[6pt]
\mbox{ affine in }  & [y_n + t_n , x_n + T_n], \\[6pt]
\mbox{ constant in }  & [x_n + T_n,  1], 
\end{array} \right.
\mbox{ and }
\left\{ \begin{array}{c}
w_n(0) = 0, \\[6pt]
w_n(y_n - t_n) =  \delta_n / 3, \\[6pt]
w_n(y_n  + t_n) = 2 \delta_n / 3, \\[6pt]
w_n (x_n  + T_n) =  \delta_n.
\end{array}\right.
$$
Set
$$
u_n = u_{n - 1} + w_n \mbox{ in } (0, 1).
$$
Since $w_n$ and $u_{n-1}$ are Lipschitz, it follows that $u_n$ is Lipschitz. Moreover, one can verify that $(u_n)$ converges in $W^{1, 1}(0, 1)$ by noting that
$$
\| w_n \|_{W^{1,1}(0,1)} \le 2 \delta_n \le 2 \delta_1/ 8^{n-1}.
$$
Let $u$ be the limit of $(u_n)$ in $W^{1, 1}(0, 1)$. We derive  from the construction of $u_n$ that $u$ is non-decreasing, and for $n \ge 1$, 
\begin{equation}\label{A-part1}
u (x) = u_n(x) \mbox{ for } x \le x_{n+1},
\end{equation}
\begin{equation}\label{A-part2}
u \mbox{ is constant in } (x_n + T_n, x_{n+1}),
\end{equation}
\begin{equation}\label{A-part3}
u(1) - u(x_n) \le \sum_{k \ge n} \delta_k < 2 \delta_n,
\end{equation}
since $\delta_k < \delta_{k-1}/ 8$. We have 
\begin{align*}
\Lambda_{2\delta_{n}}(u) = &\;  \int_0^1 \int_0^1 \frac{\varphi_{2\delta_n}(|u(x) - u(y)|)}{|x-y|^2} \, dx \, dy\\[6pt]
= & \;  \int_0^{x_n} \int_0^{x_n} \cdots  + \int_{x_n}^1 \int_{x_n}^1 \cdots  + 2 \int_0^{x_n} \int_{x_n}^{1} \cdots  \quad \mbox{ where } \cdots = \frac{\varphi_{2\delta_n}(|u(x) - u(y)|)}{|x-y|^2}.  
\end{align*}
It is clear that 
$$
 \int_0^{x_n} \int_0^{x_n} \cdots  \; \;  \mathop{\le}^{\mathrm{ by } \;  \eqref{A-part1}}  \; \;  \Lambda_{2\delta_n}(u_{n-1}),  
$$
$$
\int_{x_n}^1 \int_{x_n}^1 \cdots  \; \; \mathop{=}^{\mathrm{ by } \; \eqref{A-part3}}  \; \;  0, 
$$
and 
$$
2 \int_0^{x_n} \int_{x_n}^{1}  \cdots  \; \; \mathop{\le}^{\mathrm{ by } \;  \eqref{A-part3}}  \; \;  \int_{0}^{x_{n-1} + T_{n-1}} \, dx \int_{x_n}^1 \frac{2\delta_n c_1}{|x-y|^2}  \, dy, 
$$
since $u$ is constant in $[x_{n-1} + T_{n-1}, x_n]$. It follows from \eqref{A-p1} that 
\begin{equation}\label{A-part4}
\Lambda_{2\delta_{n}}(u)   \le \int_0^1 |u_{n-1}'| + 1/ n.
\end{equation}
On the other hand,  from  \eqref{A-p2}, \eqref{A-part1}, and the definition of $w_n$, we have, for $n \ge 1$,  
\begin{multline}\label{A-part5}
\Lambda_{\delta_n/ 3}(u) \ge  \int_{x_n}^{y_n - t_n} \, dx  \int_{y_n + t_n}^{x_n + T_n} \frac{\varphi_{\delta_n/3}(|u(x) - u(y)|)}{|x-y|^2} \, dy 
\\[6pt]= \int_{x_n}^{y_n - t_n} \, dx  \int_{y_n + t_n}^{x_n + T_n} \frac{\varphi_{\delta_n/3}(|w_n(x) - w_n(y)|)}{|x-y|^2} \, dy \\[6pt] 
\ge \int_{x_n}^{y_n - t_n} \, dx  \int_{y_n + t_n}^{x_n + T_n} \frac{c_1 \delta_n / 3}{|x - y|^2} \, dy \ge c_1 n / 3.
\end{multline}
Combining \eqref{A-part4} and \eqref{A-part5} and noting that $u_n \to u$ in $W^{1, 1}(0, 1)$, we obtain the conclusion.
\proofend

\section{Appendix: Proof of  Pathology~\ref{pathology-3}} \label{B}
\renewcommand{\theequation}{B\arabic{equation}}
\renewcommand{\thelemma}{B\arabic{lemma}}
  \setcounter{equation}{0}  
  \setcounter{lemma}{0}  

We first establish \eqref{ine-Fatou2} for $\varphi = c_1 \bvarphi_1$ where $c_1 = 1/ 2$ is the normalization constant.

Let $c \ge 5$ and for each $k \in \mN$ ($k \ge 4$) define a non-decreasing continuous function $v_{k}: [0, 1] \mapsto [0, 1]$ with $v_k(1) = 1$  as follows
\begin{equation}\label{def-vkm-Pro2}
v_{k} (x) =  \left\{\begin{array}{cl} i/k   & \mbox{ for }  x \in [i/ k, (i+1)/ k - 1/(ck)] \quad \forall \, i = 0, \cdots, k-1, \\[6pt]
\mbox{affine} & \mbox{ for }  x \in [(i+1)/ k - 1/(ck), (i+1)/ k] \quad \forall \, i = 0, \cdots, k-1.
\end{array}\right.
\end{equation}
Clearly,  
\begin{equation}\label{vk-1}
\mbox{if } |v_k(x) - v_k(y)| > 1/k \mbox{ then } |x - y | > 1/k. 
\end{equation}
%
Define 
$$
V_k(x): = \lim_{c \to + \infty} v_k(x) \quad \mbox{ for } x \in [0, 1].  
$$
Since $c_1 = 1/2$, one can show that (see \cite[page 683]{NgGammaCRAS})
\begin{equation*}
A_0: = \limsup_{k \to \infty} \Lambda_{1/k} (V_k) =  \limsup_{k \to \infty} \frac{1}{k} \sum_{i=0}^{k - 2} \int_{i/k}^{(i+1)/k} \, dx \int_{(i+2)/k}^1 \frac{1}{|x - y|^2}  \, dy  
< 1 =\lim_{\delta \to 0} \Lambda_\delta (x, [0, 1]).
\end{equation*}
Since, for $c \ge 2$,  
$$
\frac{1}{k} \sum_{i=0}^{k - 2} \int_{i/k}^{(i+1)/k} \, dx \int_{(i+2)/k -1/(ck)}^{(i+2)/k} \frac{1}{|x - y|^2}  \, dy \le 
\frac{k-1}{k} \frac{1}{k}\frac{1}{ck} \Big( \frac{1}{k}  - \frac{1}{ck}\Big)^{-2} \le  \frac{1}{c} \Big(1  - \frac{1}{c}\Big)^{-2}   \le \frac{4}{c},
$$
it follows that, for sufficiently large $c$,  
\begin{equation}\label{depart}
\limsup_{k \to \infty} \frac{1}{k} \sum_{i=0}^{k - 2} \int_{i/k}^{(i+1)/k} \, dx \int_{(i+2)/k -1/(ck)}^1 \frac{1}{|x - y|^2}  \, dy \le \frac{A_0 + 1}{2} <1. 
\end{equation}
{\bf Fix} such a constant  $c$. We are now going to define by induction a sequence of $u_n: [0, 1] \mapsto [0, 1]$.
Set
\begin{equation*}
u_0 = v_{4}.
\end{equation*}
Assume that $u_{n-1}$ ($n \ge 1$) is defined  and satisfies the following properties:
\begin{equation}\label{recurence0-Pro2}
\mbox{$u_{n-1}$ is non-decreasing, continuous, and  piecewise affine,} \quad u_{n-1}(0) = 0,  
\end{equation}
and  there exists a partition $0  = t_{0, n-1} < t_{1, n-1} < \cdots < t_{2 l_{n-1}, n-1} = 1$ such that, with the notation $J_{i, n-1} = [t_{i, n-1}, t_{i+1, n-1}]$,  the following four properties hold: 
\begin{equation}\label{recurence1-Pro2}
u_{n-1} \mbox{ is  constant on } J_{2i, n-1}  \quad \mbox{ for } i=0, \cdots, l_{n-1}-1, 
\end{equation}
\begin{equation}\label{recurence2-Pro2}
u_{n-1} \mbox{ is  affine and not constant on } J_{2i+1, n-1} \quad \mbox{ for } i =0, \cdots, l_{n-1}-1,
\end{equation}
the total variation of $u_{n-1}$ on the interval $J_{i, n-1}$ with $i$ odd (where $u_{n-1}$ is not constant) is always $1/ l_{n-1}$, i.e.,
\begin{equation}\label{recurence3-Pro2} u_{n-1}(t_{2i + 2, n-1}) - u_{n-1}( t_{2i+1, n-1} ) = 1/ l_{n-1}  \quad \mbox{ for } i =0, \cdots, l_{n-1}-1,
\end{equation}
and  the intervals $J_{i, n-1}$ with $i$ odd have the same length which is less than the one of any interval $J_{i, n-1}$ with $i$ even, i.e., 
\begin{equation}\label{recurence4-Pro2}
|J_{1, n-1}| = |J_{3, n-1}|  = \cdots =   |J_{2l_{n-1}-1, n-1}|  < |J_{2i, n-1}| \quad  \mbox{ for } i = 0, \cdots, l_{n-1}-1. 
\end{equation}
Since $u_{n-1} (0) = 0$,  it follows from the properties of $u_{n-1}$ in \eqref{recurence1-Pro2} and \eqref{recurence2-Pro2} that 
\begin{equation}\label{u-n-1-1} 
 u_{n-1}(t) = s/ l_{n-1} +  i/ l_{n-1}  \mbox{ for }  t \in J_{2i+1, n-1} \mbox{ where } s  = (t - t_{2i +1, n-1}) / |J_{2i + 1, n-1}|. 
\end{equation}
Set
\begin{equation}\label{def-Bn}
B_{n-1} = \bigcup_{i=0}^{l_{n-1}-1} J_{2i, n-1}
\end{equation}
($B_{n-1}$ is the union of all intervals on which $u_{n-1}$ is constant).
For $n \in \mN$, let $k_n$  be a sufficient large integer such that
\begin{equation}\label{def-k1}
\frac{1}{k_n}  \int_{-1}^0 \, dx \int_{\tau_n}^1 \frac{1}{|x-y|^2}  \, dy <  \frac{1}{n} \quad \mbox{ where } \quad  \tau_n = |J_{1, n-1}| / k_n
\end{equation}
and
\begin{equation}\label{def-k2-1}
\frac{1}{k_n} \sum_{i=0}^{k_n - 2} \int_{i/k_n}^{(i+1)/k_n} \, dx \int_{(i+2)/k_n -1/(ck_n)}^1 \frac{1}{|x - y|^2}  \, dy  \le  \frac{A_0+1}{2}.
\end{equation}
Since, for a small positive number $\tau$,  
\begin{equation*}
 \int_{-1}^0 \, dx \int_{\tau}^1 \frac{1}{|x-y|^2}  \, dy \le |\ln \tau|, 
\end{equation*}
such a constant $k_n$ exists by \eqref{depart}. 
Define
\begin{equation}\label{def-un-Pro2}
u_n (t) = \left\{\begin{array}{cl} u_{n-1}(t) & \mbox{ if } t \in   B_{n-1},  \\[6pt]
\dsp \frac{1}{l_{n-1}}v_{k_n} (s) + \frac{i}{l_{n-1}}& \mbox{ if } t \in J_{2i + 1, n-1} \mbox{ for some } 0 \le i \le l_{n-1} -1,
\end{array}\right.
\end{equation}
where $s  = (t - t_{2i +1, n-1}) / |J_{2i + 1, n-1}|$. Then $u_n$  satisfies \eqref{recurence0-Pro2}-\eqref{recurence4-Pro2} for some  $l_n$ and $t_{i, n}$.
Since  $0 \le v_k(x) \le x$ for $x \in [0, 1]$, we deduce  from \eqref{u-n-1-1} and  the definition of $u_{n}$ that $u_n \le u_{n-1}$. On the other hand, we derive from \eqref{u-n-1-1} and  \eqref{def-un-Pro2} that, for $m\ge n$,  
$$
\|u_m - u_n \|_{L^\infty(0, 1)} \le 1/ l_n. 
$$
Hence the sequence $(u_n)$ is Cauchy in $C([0, 1])$. 
Let $u$ be the limit and set
\begin{equation}\label{def-dn}
\delta_n =  1/ (l_{n-1} k_n).
\end{equation}
It follows from the definition of $u_n$ and $u$ that 
\begin{equation}\label{u=u_n}
u(t) = u_n(t)  \mbox{ for }  t = t_{i, n-1} \mbox{ with } 0 \le i \le 2 l_{n-1}.
\end{equation}
From the construction of $u_n$ in \eqref{def-un-Pro2}, the property of $v_k$ in  \eqref{vk-1}, and \eqref{recurence4-Pro2}, we derive that 
\begin{equation}\label{lip-Pro2}
\mbox{ if }| u(x) -  u(y)| > \delta_n, \mbox{ then }  |x  - y| > \tau_n,
\end{equation}
where $\tau_n$ is defined in \eqref{def-k1}. Since $u_{n-1}$ is constant in $J_{2i, n-1}$ for $0 \le i \le l_{n-1} - 1$ by \eqref{recurence1-Pro2},  it follows from \eqref{def-un-Pro2} that $u$ is constant   in  $J_{2i, n-1}$ for $0 \le i \le l_{n-1} - 1$.  We derive that   
\begin{multline}\label{p1-Pro2}
\mathop{\int_0^1 \int_0^1}_{|u(x) - u(y)| > \delta_n} \frac{\delta_n}{|x - y|^2} \, dx \, dy \le \sum_{i=0}^{l_{n-1} - 1} \mathop{\iint_{J_{2i+1, n-1}^2}}_{|u(x) - u(y)| > \delta_n} \frac{\delta_n}{|x - y|^2} dx \, dy  \\[6pt]
 +   \sum_{i=0}^{2l_{n-1}-1} \mathop{\int_{J_{i, n-1}}  \, dx \int_{[0, 1] \setminus J_{i, n-1}}}_{|u(x) - u(y)| > \delta_n} \frac{\delta_n}{|x - y|^2 }  \, dy.
\end{multline}
Using \eqref{lip-Pro2}, we have, by \eqref{def-k1},
\begin{equation}\label{p2-Pro2}
 \sum_{i=0}^{2l_{n-1}-1} \mathop{\int_{J_{i, n-1}}  \, dx \int_{[0, 1] \setminus J_{i, n-1}}}_{|u(x) - u(y)| > \delta_n} \frac{\delta_n}{|x - y|^2 }  \, dy
\le 4 l_{n-1} \int_{-1}^0  \, dx \int_{\tau_n}^1 \frac{ \delta_n }{|x-y|^2} \, dy \le 4/ n.
\end{equation}
We now estimate, for $0 \le i \le l_{n-1} - 1$, 
$$
\mathop{\iint_{J_{2i+1, n-1}^2}}_{|u(x) - u(y)| > \delta_n} \frac{\delta_n}{|x - y|^2} dx \, dy.
$$
Define $g_i : J_{2i + 1, n-1} \to [0,  1]$,  for $0 \le i \le l_{n-1} - 1$,  as follows 
$$
g_i(x) =  (x - t_{2i +1, n-1}) / |J_{2i + 1, n-1}| \quad \mbox{ for } x \in J_{2i + 1, n-1}. 
$$
We claim that, for $0 \le i \le l_{n-1} - 1$,
$$
\mbox{ if } (x, y) \in J_{2i +1, n-1}^2,  \quad |u(x) - u(y)| > \delta_n,   \quad g_i(x) \in [i/ k_{n}, (i+1)/ k_n], \quad \mbox{ and } \quad  x < y, 
$$
then  
$$
g_i(y) \in \big[(i+2)/ k_n - 1/ (c k_n), 1\big].
$$ 
In fact, if $g_i(z) \in [i/ k_{n},  (i+2)/ k_n - 1 / (c k_n)]$ then 
$$
u_n \left(g_i^{-1}\Big(\frac{i}{k_n} \Big) \right) =u \left(g_i^{-1}\Big(\frac{i}{k_n} \Big) \right)  \le u(z) \le u \left(g_i^{-1}\Big( \frac{i+2}{k_n} -  \frac{1}{c k_n}  \Big) \right) =  u \left(g_i^{-1}\Big( \frac{i+2}{k_n} -  \frac{1}{c k_n}  \Big) \right). 
$$
Here we used \eqref{u=u_n}  and the fact that $u$ is non-decreasing. 
It follows from the definition of $u_n$ that,   if $g_i(x), g_i(y) \in [i/ k_{n},  (i+2)/ k_n - 1 / (c k_n)]$ then 
$$
|u(y) - u(x)| \le \frac{1}{l_{n-1}} \left|v_{k_n} \Big( \frac{i+2}{k_n} - \frac{1}{c k_n}  \Big) - v_{k_{n}} \Big( \frac{i}{k_n} \Big) \right| \le \frac{1}{k_n l_{n-1}} = \delta_n. 
$$
The claim is proved. 

By a change of variables, for $i=0, \cdots, 2 l_{n-1} -1$,   
$$
(x, y) \mapsto \Big(g_i(x), g_i(y) \Big)  \mbox{ for } (x , y) \in J_{2i+1, n-1}^2, 
$$ 
we deduce from the claim  that 
\begin{align*}
\sum_{i=0}^{l_{n-1} - 1} \mathop{\iint_{J_{2i+1, n-1}^2}}_{|u(x) - u(y)| > \delta_n} \frac{\delta_n}{|x - y|^2} dx \, dy 
& \le 2 l_{n-1} \delta_n \sum_{j=0}^{k_n - 2} \int_{j/k_n}^{(j+1)/k_n} \, dx \int_{(j+2)/k_n -1/(ck_n)}^1 \frac{1}{|x - y|^2}  \, dy. 
\end{align*}
It follows from  \eqref{def-k2-1} and \eqref{def-dn} that  
\begin{equation}\label{p3-Pro2}
\sum_{i=0}^{l_{n-1} - 1} \mathop{\iint_{J_{2i+1, n-1}^2}}_{|u(x) - u(y)| > \delta_n} \frac{\delta_n}{|x - y|^2} dx \, dy 
 \le  A_0 + 1.  
\end{equation}
Combining  \eqref{p1-Pro2}, \eqref{p2-Pro2}, and \eqref{p3-Pro2} yields
\begin{equation*}
\limsup_{n \to \infty} \mathop{\int_0^1 \int_0^1}_{|u(x) - u(y)| > \delta_n} \frac{\delta_n}{|x - y|^2} \, dx \, dy  \le A_0 + 1.
\end{equation*}
Since $c_1 = 1/2$, we have, for $\varphi= c_1 \bvarphi_1$,
\begin{equation}\label{lim-u-Pro2}
\limsup_{n \to \infty} \Lambda_{\delta_n}(u) \le (A_0 + 1)/ 2 < 1.
\end{equation}
Note that $u \in C([0, 1])$ is non-decreasing and $u(0) =0$ and $u(1) = 1$. This implies 
$$
\int_0^1|u'| = 1. 
$$
 Therefore \eqref{ine-Fatou2}  holds  for $\varphi = c_1 \bvarphi_1$ and $u$.

We next construct  a continuous function $\varphi_\ell$ which is ``close" to $c_1 \bvarphi_1$ such that  \eqref{ine-Fatou2} holds for $\varphi_\ell$ and the function $u$ constructed above. For $\ell \ge 1$, define a continuous function   $\varphi_\ell: [0 , + \infty) \mapsto \mR$ by
\begin{equation*}
\varphi_\ell (t) = \left\{\begin{array}{cl} \alpha_\ell  & \mbox{ if } t \ge 1 + 1/ \ell, \\[6pt]
0 & \mbox{ if } t \le 1, \\[6pt]
\mbox{affine} & \mbox{ if } t \in [1, 1 + 1/ \ell], 
\end{array}\right.
\end{equation*}
where $\alpha_\ell$ is the constant such that
\begin{equation*}
\gamma_{1} \int_0^\infty \varphi_\ell(t) t^{-2} \, dt  = 2 \int_0^\infty \varphi_\ell(t) t^{-2} \, dt = 1.
\end{equation*}
Then $\varphi_{\ell} \in {\cal A}$. Moreover,  $\varphi_\ell (t) \le \alpha_\ell \bvarphi_1 (\beta_\ell t)$ where $\beta_\ell = 1 + 1/\ell$.  It follows  from \eqref{lim-u-Pro2} that
\begin{multline*}
\liminf_{\delta \to 0}  \int_{0}^1 \int_{0}^1  \frac{\beta_\ell \delta \varphi_\ell (|u(x) - u(y)| / (\beta_\ell \delta) ) }{|x-y|^{2}} \, dx \, dy  \\[6pt]
\le \liminf_{\delta \to 0}  \int_{0}^1 \int_{0}^1  \frac{\alpha_\ell  \beta_\ell \delta\bvarphi_1 \Big(|u(x) - u(y)| \big/ \delta \Big) }{|x-y|^{2}} \, dx \, dy
\le  a_\ell \beta_\ell (A_0 + 1).
\end{multline*}
Since $a_\ell \to c_1 =  1/2$  and $\beta_\ell \to 1$ as $\ell \to +\infty $,  the conclusion holds for $\varphi_\ell$ when  $\ell$ is large. The proof is complete. 
\proofend

\section{Appendix: Pointwise convergence of $\Lambda_{\delta}(u)$ when $u \in W^{1, p}(\Omega)$} \label{C}
\renewcommand{\theequation}{C\arabic{equation}}
\renewcommand{\thelemma}{C\arabic{lemma}}
\renewcommand{\theproposition}{C\arabic{proposition}}

  \setcounter{equation}{0}  
  \setcounter{lemma}{0}  
 \setcounter{proposition}{0} 

In this section, we prove the following result

\begin{proposition} \label{pro-C} Let $d \ge 1$, $\Omega$ be a smooth bounded open subset  of $\mR^d$, and $\varphi \in {\cal A}$. We have
$$
\lim_{\delta \to 0} \Lambda_{\delta} (u) = \int_{\Omega} |\nabla u|  \mbox{ for } u \in \bigcup_{p> 1} W^{1, p}(\Omega)
$$
\end{proposition}

\noindent{\bf Proof.} We already know by Proposition~\ref{thm-Sobolev1p=1} that 
\begin{equation}
\liminf_{\delta \to 0} \Lambda_{\delta} (u) \ge \int_{\Omega} |\nabla u| \quad \forall \, u \in W^{1, 1}(\Omega). 
\end{equation}
Assume now that $u \in W^{1, p}(\Omega)$ for some $p>1$. We are going to prove that 
\begin{equation}
\limsup_{\delta \to 0} \Lambda_{\delta} (u) \le \int_{\Omega} |\nabla u|. 
\end{equation}
Consider an extension of $u$ to $\mR^d$ which belongs to $W^{1, p}(\mR^d)$. For simplicity, we still denote the extension by  $u$. 

Clearly
$$
\Lambda_{\delta}(u) \le \int_{\Omega} dx \int_{\mR^d} \frac{\varphi_\delta(|u(x) - u(y)|)}{|x - y|^{d+1}} \, dy, 
$$
and thus it suffices to establish that 
\begin{equation}\label{C-estimate-*}
\lim_{\delta \to 0} \int_{\Omega} dx \int_{\mR^d} \frac{\varphi_\delta(|u(x) - u(y)|)}{|x - y|^{d+1}} \, dy = \int_{\Omega} |\nabla u| \, dx. 
\end{equation}
Using polar coordinates and  a change of variables, we have, as in \eqref{tt4}, 
\begin{equation}\label{C-part0}
 \int_{\Omega} \, d x \int_{\mR^d}  \frac{\varphi_\delta  (|u(x) - u(y)|) }{|x - y|^{ d  +  1}} \,
dy = \int_{\Omega} \, dx \int_{0}^{\infty}  \, d h  \int_{\mS^{d-1}}    \frac{1}{h^2}\varphi\Big(|u(x + \delta h \sigma) - u(x)| \big/ \delta \Big)  \, d \sigma. 
\end{equation}
As in \eqref{mono1}, we also obtain
\begin{multline}\label{C-part0-1}
\lim_{\delta \to 0}\frac{1}{h^{2}} \varphi \Big(|u(x + \delta h \sigma) - u(x)| \big/ \delta \Big) =
\frac{1}{h^{2}}\varphi \Big(|\nabla u(x) \cdot \sigma| h \Big) \\[6pt]
\mbox{ for  a.e. } (x, \, h, \, \sigma) \in \Omega \times (0, + \infty) \times \mS^{d-1}. 
\end{multline}
As in \eqref{limit}, we have 
\begin{equation}\label{C-part0-2}
\int_{\Omega} \, dx  \int_{0}^\infty  \, d h \int_{\mS^{d-1}} \frac{1}{h^{2}} \varphi\Big( |\nabla u (x) \cdot \sigma| h\Big)  \, d \sigma
 =  \int_{\Omega} |\nabla u| \, dx. 
\end{equation}
On the other hand, since $\varphi$ is non-decreasing, it follows that, for $\delta >0$,  
\begin{equation}\label{C-part1}
 \frac{1}{h^2}\varphi \big(|u(x + \delta h \sigma) - u(x)| / \delta \big) \le \frac{1}{h^2} \varphi \big(M_{\sigma} (\nabla u) (x) h \big) \mbox{ for  a.e. } (x, \, h, \, \sigma) \in \mR^d \times (0, + \infty) \times \mS^{d-1},   \end{equation}
where 
\begin{equation}\label{C-def-M}
M_{\sigma} (\nabla u)(x): = \sup_{\tau > 0} \int_{0}^1 |\nabla u(x + s \tau  \sigma) \cdot \sigma| \, ds \quad  \mbox{ for }   x \in \mR^d. 
\end{equation}
Indeed, we have 
\begin{equation*}
|u(x  + \delta h \sigma) - u(x)| / \delta \le \int_0^1 h  |\nabla u(x + s \delta h \sigma) \cdot \sigma | \, ds
\le  h \sup_{\tau > 0}  \int_0^1 h  |\nabla u(x + s \tau \sigma) \cdot \sigma | \, ds. 
\end{equation*}
We claim that 
\begin{equation}\label{C-claim2}
 \frac{1}{h^2} \varphi \big( M_\sigma(\nabla u)(x) h \big) \,d  x  \in L^1 \big(\Omega \times (0, + \infty) \times \mS^{d-1} \big). 
\end{equation}
Assuming \eqref{C-claim2}, we may then apply the dominated convergence theorem using \eqref{C-part0}, \eqref{C-part0-1}, \eqref{C-part0-2}, \eqref{C-part1}, and \eqref{C-claim2},  and conclude that \eqref{C-estimate-*} holds. 

To show \eqref{C-claim2}, it suffices to prove that, for all $\sigma \in \mS^{d-1}$, 
\begin{equation}\label{C-claim1}
 \int_{\Omega} dx \int_{0}^{\infty}   \frac{1}{h^2} \varphi \big( M_\sigma(\nabla u)(x) h \big) \,d  h   \le C \Big(\int_{\mR^d} |\nabla u|^p \Big)^{1/p}. 
\end{equation}
Here and in what follows $C$ denotes a positive constant independent of $u$ and $\delta$; it depends only on $\Omega$ and $p$. 
For simplicity of notation, we  assume that $\sigma = e_{d}: = (0, \cdots, 0, 1 )$.
By a change of variables, we have 
\begin{align}\label{C-part2} 
  \int_{\Omega}  dx  \int_{0}^{\infty}  \frac{1}{h^2} \varphi\big( M_{e_d}(\nabla u)(x) h \big) \,d  h 
= &  \int_{\Omega} \big|M_{e_d} (\nabla u) (x)\big|  \, d  x  \int_0^\infty \varphi(t) t^{-2} \, d t   \nonumber \\[6pt] 
= &   \gamma_d^{-1} \int_{\Omega} \big|M_{e_d} (\nabla u) (x)\big|  \, d  x \le C\Big( \int_{\Omega} \big|M_{e_d} (\nabla u) (x)\big|^p  \, d  x \Big)^{1/p}.
\end{align}
Note that 
$$
M_{e_d}(\nabla u)(x) =  \sup_{\tau > 0} \int_{0}^1 |\partial_{x_d} u(x', x_d + s \tau)| \, ds  = \sup_{\tau > 0} \mint_{x_d}^{x_d + \tau}  |\partial_{x_d} u(x', s)| \, ds. 
$$
We have 
\begin{multline}\label{C-haha}
\int_{\Omega} \big|M_{e_d} (\nabla u) (x)\big|^p \, dx \le \int_{\mR^d} \big|M_{e_d} (\nabla  u) (x)|^p \, dx
=  \int_{\mR^{d-1}} \, d x' \int_{\mR} \big|M_{e_d} (\nabla u) (x', x_d) \big|^p \, dx_d. 
\end{multline}
Since, by the theory of maximal functions in one dimension,
\begin{equation*} 
\int_{\mR} \big|M_{e_d} (\nabla u) (x', x_d) \big|^p \, dx_d  \le C
 \int_{\mR} |\partial_{x_d} u (x', x_d)|^p \, dx_d,
\end{equation*}
it follows from \eqref{C-haha} that 
\begin{equation}\label{C-haha-2}
\int_{\Omega} \big|M_{e_d} (\nabla u) (x)\big|^p \, dx \le C
 \int_{\mR^{d} }  |\nabla u (x)|^p \, dx. 
\end{equation}
Combining \eqref{C-part2} and \eqref{C-haha-2} implies \eqref{C-claim1}  for $\sigma = e_d$.  The proof is complete. \proofend

\begin{remark}  \fontfamily{m} \selectfont 
The above proof shows that 
$$
\Lambda_\delta(u) \le C_p \| \nabla u\|_{L^p(\Omega)} \quad \forall \, u \in W^{1, p}(\Omega). 
$$
Such an estimate is inspired by ideas due to H-M. Nguyen \cite{NgSob1} and A. Ponce and J. Van Schaftingen in \cite{Po2}.  
\end{remark}

\bigskip
\noindent{\bf Acknowledgement.} We are extremely grateful  to J. Bourgain for sharing fruitful ideas which led to the joint work \cite{BourNg} with the second author,  and served as an ``early  catalyst" for our works in this direction.  The first author (H.B.)  warmly thanks R. Kimmel and J. M. Morel  for useful discussions concerning Image Processing.

\providecommand{\bysame}{\leavevmode\hbox to3em{\hrulefill}\thinspace}
\providecommand{\MR}{\relax\ifhmode\unskip\space\fi MR }
\providecommand{\MRhref}[2]{%
  \href{http://www.ams.org/mathscinet-getitem?mr=#1}{#2}
}
\providecommand{\href}[2]{#2}

\end{document}